\pdfoutput=1
\RequirePackage{ifpdf}
\ifpdf 
\documentclass[pdftex]{sigma}
\else
\documentclass{sigma}
\fi

\usepackage[labelfont=bf,labelsep=period]{caption}

\usepackage{pgf,pstricks,tikz}
\usepackage[absolute]{textpos}
\usepackage{stmaryrd}
\usepackage{dsfont}
\usetikzlibrary{shapes,snakes}
\usetikzlibrary{positioning}
\usepackage{ytableau}
\usepackage{dynkin-diagrams}
\usepackage{genyoungtabtikz}
\usepackage{makecell}

\newtheorem{Theorem}{Theorem}[section]
\newtheorem{Corollary}[Theorem]{Corollary}
\newtheorem{Lemma}[Theorem]{Lemma}
\newtheorem{Proposition}[Theorem]{Proposition}
 { \theoremstyle{definition}
\newtheorem{Definition}[Theorem]{Definition}

\newtheorem{Example}[Theorem]{Example}
\newtheorem{Remark}[Theorem]{Remark} }

\numberwithin{equation}{section}

\newcommand{\al}{\alpha}

\newcommand{\ga}{\gamma}

\newcommand{\La}{\Lambda}

\newcommand{\fh}{\mathfrak{h}}
\newcommand{\fg}{\mathfrak{g}}
\newcommand{\cC}{\mathcal{C}}

\def\cch{\mathcal H}
\def\ccsc{\mathcal {SC}}
\def\ccdd{\mathcal {DD}}
\def\ccd{\mathcal {D}}
\def\ccp{\mathcal P}

\DeclareMathOperator{\core}{core}
\DeclareMathOperator{\quot}{quot}

\newcommand{\scal}[2]{\langle #1,#2\rangle}

\newcommand{\Z}{\mathbb{Z}}
\newcommand{\N}{\mathbb{N}}

\newcommand{\func}[1]{\operatorname{#1}}

\begin{document}

\newcommand{\arXivNumber}{2404.10532}

\renewcommand{\thefootnote}{}

\renewcommand{\PaperNumber}{023}

\FirstPageHeading

\ShortArticleName{Macdonald identities, Weyl--Kac denominator formulas and affine Grassmannian elements}

\ArticleName{Macdonald Identities, Weyl--Kac Denominator\\ Formulas and Affine Grassmannian Elements\footnote{This paper is a~contribution to the Special Issue on Basic Hypergeometric Series Associated with Root Systems and Applications in honor of Stephen C.~Milne's 75th birthday. The~full collection is available at \href{https://www.emis.de/journals/SIGMA/Milne.html}{https://www.emis.de/journals/SIGMA/Milne.html}}}

\Author{C\'{e}dric LECOUVEY and David WAHICHE}

\AuthorNameForHeading{C.~Lecouvey and D.~Wahiche}
\Address{Univ. de Tours, UMR CNRS 7013, Institut Denis Poisson, France}
\Email{\href{mailto:cedric.lecouvey@univ-tours.fr}{cedric.lecouvey@univ-tours.fr}, \href{mailto:wahiche@univ-tours.fr}{wahiche@univ-tours.fr}}
\URLaddress{\url{https://www.idpoisson.fr/lecouvey/}, \url{https://www.idpoisson.fr/wahiche/}}

\ArticleDates{Received June 20, 2024, in final form March 04, 2025; Published online April 10, 2025}

\Abstract{The Nekrasov--Okounkov formula gives an expression for the Fourier coefficients of the Euler functions as a sum of hook length products. This formula can be deduced from a specialization in a renormalization of the affine type $A$ Weyl denominator formula and the use of a polynomial argument. In this paper, we rephrase the renormalized Weyl--Kac denominator formula as a sum parametrized by affine Grassmannian elements. This naturally gives rise to the (dual) atomic length of the root system considered introduced by Chapelier-Laget and Gerber. We then provide an interpretation of this atomic length as the cardinality of some subsets of $n$-core partitions by using foldings of affine Dynkin diagrams. This interpretation does not permit the direct use of a polynomial argument for all affine root systems. We show that this obstruction can be overcome by computing the atomic length of certain families of integer partitions. Then we show how hook-length statistics on these partitions are connected with the Coxeter length on affine Grassmannian elements and Nekrasov--Okounkov type formulas.}


\Classification{05E05; 05E10; 11P81}

\renewcommand{\thefootnote}{\arabic{footnote}}
\setcounter{footnote}{0}

\begin{flushright}
\begin{minipage}{58mm}
\it To Stephen Milne on the occasion\\ of his 75th birthday
\end{minipage}
\end{flushright}

\section{Introduction}

The Weyl--Kac formula is a cornerstone in the representation theory of
infinite-dimensional Lie algebras. It
permits the computation of the character of a highest weight irreducible
representation and naturally generalizes the classical Weyl character
formula relevant for the finite-dimensional simple Lie algebras over
the complex number field. The Weyl denominator formula (see Section \ref{WKDF} for the notation) can be written%
\begin{equation*}
\prod_{\alpha \in R_{+}}(1-{\rm e}^{-\alpha })^{m_{\alpha }}=\sum_{w\in
	W_{a}}\varepsilon (w){\rm e}^{w(\rho )-\rho }
\end{equation*}%
and reflects the fact that the trivial representation has a character equal
to $1$. When applied to affine root systems and by putting $q={\rm e}^{-\delta}$, it yields a rich class of $q$-series. Note  the seminal work of Lepowsky and Milne in~\cite{MilneLepowsky}, exploring the connection between partition identities and the Weyl--Kac formula or a proof of the identity in affine type \smash{$A_{n-1}^{(1)}$} (see~\cite{Milne}), using a multiple series generalization of the $q$-binomial theorem. Moreover, as proved by Han~\cite{Ha} by using the affine
root systems of type $A$, it yields the Nekrasov--Okounkov
formula%
\begin{equation}\label{NOdebut}
\sum_{\lambda \in \mathcal{P}}q^{\vert \lambda \vert }\prod_{h\in
	\mathcal{H}(\lambda )}\biggl( 1-\frac{z}{h^{2}}\biggr) =\prod_{k\geq
	1}\bigl(1-q^{k}\bigr)^{z-1},
\end{equation}%
where $\mathcal{P}$ is the set of partitions, $\mathcal{H}(\lambda )$ the
multiset of hook lengths in the partition $\lambda $ and $z$ is any fixed
complex number. When $z=0$, one recovers the Euler
generating function for the set of partitions. By assuming that $z=n$, with $n
$ an integer greater than or equal to $2$, it is also possible to derive
interesting generating series identities for $n$-core partitions,
that is, the subset of $\mathcal{P}$ containing exactly the partitions with
no hook length equal to $n$.

As mentioned in~\cite{CRV,RW}, where a two parameter generalization of \eqref{NOdebut} is derived, there exists a~$z,u$-analogue of the Nekrasov--Okounkov formula. It can be regarded as a~reformulation of a~result due to Dehaye and Han~\cite{HD} based on a specialization of the Macdonald identity in type~\smash{$A^{(1)}_{n-1}$} or obtained by using the refined topological vertex as in~\cite{INRS}. This $z,u$-analogue can be written as follows:
\begin{equation}\label{Hande}
\sum_{\lambda\in\ccp}q^{\lvert \lambda\rvert}\prod_{h\in\cch(\lambda)}\frac{\bigl(1-zu^h\bigr)\bigl(1-z^{-1}u^h\bigr)}{\bigl(1-u^h\bigr)^2}=\prod_{k,r\geq 1}\frac{\bigl(1-zu^rq^k\bigr)^r\bigl(1-z^{-1}u^rq^k\bigr)^r}{\bigl(1-u^{r-1}q^k\bigr)^r\bigl(1-u^{r+1}q^k\bigr)^r}.
\end{equation}
\noindent Note that taking $z=u^n$ for a strictly positive integer $n$ and letting $q\rightarrow 1$ in \eqref{Hande} yields \eqref{NOdebut} for $z=n^2$. By remarking that for any positive integer $k$, the coefficients in $q^k$ on both sides of~\eqref{NOdebut} are polynomials in $z$, we can use a polynomial argument, which shows the identity for any complex number $z$.

In fact, one can deduce numerous generalizations of the
Nekrasov--Okounkov formula from the Weyl--Kac denominator formula, not only in
affine type $A$, but also for the other affine root~systems. This was done by P\'etr\'eolle for types \smash{$C_n^{(1)}$} and \smash{$D_{n+1}^{(2)}$} (see~\cite{MP}) and by the second~author for all seven infinite families of affine Lie algebras in~\cite[Section $5.3$]{Wmac} to which we refer the reader
for a more complete introduction to the history and developments concerning partition identities linked to the Weyl--Kac denominator formula. A central idea of~\cite{Wmac} is to get expansions of
the Weyl denominator formula associated to classical affine root systems in
terms of the irreducible Weyl characters, that is, the characters of the
irreducible highest weight representations corresponding to the underlying
finite-dimensional Lie algebras. Another important aspect of the results
proved in~\cite{Wmac} is the use, as index sets for the previous expansion,
of~families of self-conjugate and doubled distinct partitions (i.e., some concatenations of two copies of a partition with distinct parts) which can be regarded as natural analogues of core partitions for each classical affine root system. Nevertheless, the methods used are based on case-by-case
computations and combinatorial manipulations on binary codings of
partitions through the Littlewood decomposition (see, for instance,~\cite{HJ}).

In contrast, our goal here is to fully use the powerful machinery of
affine root systems, as developed by Macdonald and Kac and presented in~\cite{carter_affine,Mac}
or~\cite{KacB}, to first expand the affine Weyl--Kac denominator formula in
terms of the ordinary Weyl characters by using a summation over the
affine Grassmannian elements. Recall here these are the minimal-length
elements in the cosets of an affine Weyl group by its classical parabolic
subgroup. The set of affine Grassmannian elements (called the affine
Grassmannian in the literature) comes with a natural statistics called the
\textquotedblleft atomic length\textquotedblright\ introduced and studied
(in a more general context) by Chapelier-Laget and Gerber in~\cite{CGT}. In
affine type \smash{$A_{n-1}^{(1)}$}, the affine Grassmannian elements are in
one-to-one correspondence with the $n$-cores and the atomic length is just
the number of boxes of the associated Young diagram. For the classical
affine root systems (twisted or not) and in types~\smash{$G_{2}^{(1)}$} and~\smash{$D_{4}^{(3)}$}
(here we follow the classification of affine root systems in
Kac's book~\cite{KacB}), we show that it is possible to parametrize the
affine Grassmannian elements by using particular subsets of self-conjugate $2n
$-cores. In each case, we are able to write a simple combinatorial formula
for their atomic lengths just in terms of their number of boxes and the
number of their boxes (or nodes) of residues $0$ or $n$. In particular, in
the non-twisted cases, we recover the formulas obtained in~\cite{STW}.
However, our methods---based on folding of Dynkin diagrams---are different from those of~\cite{STW} and also well-suited to
consider the twisted affine cases. Alternatively, one can compute the atomic
lengths by counting boxes in the Young diagrams of some self-conjugate $2n$%
-cores with weights on boxes depending on their residue (in type $A$, all
the weights are equal to~$1$). This approach has the advantage of
homogeneity, but it is not~well-adapted for considering~at the same time a
given family \smash{$\bigl(X_{n}^{(a)}\bigr)_{n\geq 2}$} of Dynkin diagrams (i.e., when the
type is fixed and one considers all the possible ranks at once). Also in
order to get a natural parametrization of the solutions of certain Diophantine
equations, as proposed in~\cite{BCT}, it is important to label the affine
Grassmannian elements of the previous types by families of partitions for
which the atomic length equates the numbers of boxes. This is what we
explain in~Section~\ref{sec:combi}, where we show how the families of distinguished
self-conjugate $2n$-cores obtained in~Section~\ref{sec:atomic} are related to the families
of partitions introduced in~\cite{Wmac} by simple bijections sending the atomic length to the number of boxes. This also gives us the dictionary to connect the Nekrasov--Okounkov type formulas as stated in~\cite{Wmac} with the affine Grassmannian elements.

The paper is organized as follows. In Section \ref{sec:2}, we recall the background on affine root systems on which the affine Weyl denominator formulas are based. We then decompose this affine Weyl denominator in terms of the generalized Weyl characters in Section \ref{sec:3}. Observe that these expressions are essentially equivalent to Theorems 20.03 and 20.04 in Carter's book~\cite{carter_affine}. In Section~\ref{sec:4}, we simplify the previous decomposition in order to get a decomposition in terms of genuine characters, i.e., characters labelled by dominant weights. This also makes the appearance of affine Grassmannian elements more natural. In Section \ref{sec:atomic}, we use techniques of foldings of Dynkin diagram techniques to
realize each affine root system in an affine root system of type~$A$. This
permits in particular to get the desired formulas for the atomic length in
terms of weighted boxes in self-conjugate $2n$-cores. Section \ref{sec:combi} presents the combinatorics of partitions and the Littlewood decomposition which are the keys to understand the connections between the different combinatorial models that we employ. This also allows for the description of simple bijections between the model of self-conjugate cores and the model
of partitions obtained in~\cite{Wmac}. The~basic idea here is first to
identify each $2n$-core with its so-called $2n$-charge $\beta $ \big(a vector in~$\mathbb{Z}^{2n}$ whose sum of coordinates is equal to zero\big); next we
express the atomic lengths in each type in terms of~$\beta $; and finally, we
add a number $a$ of zero parts to $\beta $ in order to be able to identify the atomic length obtained as a number of boxes in a relevant partition.
When the Dynkin diagram considered has no sub-diagram of classical type $D$,
it suffices to use $(2n+a)$-core partitions or some of their
half-reductions regarded as distinct partitions. Otherwise, the situation becomes more complicated, but one
can still reduce the problem to simple families of distinct partitions, even if
they do not come from $(2n+a)$-cores in general. In Section~\ref{sec:abacus}, we
show how the dominant weights appearing in the Weyl character expansion of
the affine Weyl denominator formula of Section~\ref{sec:4} and obtained from the affine Grassmannian elements can be computed directly
from the previous combinatorics of partitions thanks to the notion of
$V_{(g,n)}$-coding introduced in~\cite{Wmac}. This yields in particular a more uniform presentation of the generalized Nekrasov--Okounkov formulas.
 Finally, in Proposition \ref{prop:inversion} we exhibit how, in affine type $A$, some refinements of the inversion sets containing roots with fixed heights have the same cardinality as some subsets of hook lengths in the corresponding core partition.

\section{Affine root systems and affine Lie algebras}\label{sec:2}


In this section, we recal some facts about affine root systems and affine Lie algebras. In particular, we detail how the affine Weyl group $W_a$ acts on the affine weight lattice. In fact, the affine Weyl group is the semi-direct product of the underlying finite Weyl group $W$ by a~group of translations. Roughly speaking, this explains why non-affine characters arise in our rewriting of the affine Weyl denominator formulas.

Set $I=\{0,1,\dots,n\}$, $I^\ast=\{1,2,\dots,n\}$ and let \smash{$A=(a_{i,j})_{(i,j)\in I^{2}}$} be a generalized Cartan matrix of a classical affine root system. These affine root systems
are classified in~\cite[p.~54, Table~1]{KacB}. The
matrix $A$ has rank $n$ and there exists a unique vector $v = (a_i)_{i\in I}\in\mathbb{Z}^{n+1}$ with $(a_i)_{i\in I}$ relatively prime and a unique vector \smash{$v^{\vee} =\bigl(a_i^\vee\bigr)_{i\in I}\in \Z^{n+1}$} with \smash{$\bigl(a^\vee_i\bigr)_{i\in I}$} relatively prime~such that \smash{$v^{\vee}\cdot A=A\cdot {^t}\hspace{-.1mm}v=0$}. Note that we have $a_0=a_0^\vee = 1$ except in the \smash{$A_{2n}^{(2)}$} case for~which $a_0=2$ and \smash{$a_0^\vee = 1$}. We refer to~\cite{carter_affine,KacB} for details and proofs of the results presented in this section.

Let $\bigl(\fh,\Pi,\Pi^\vee\bigr)$ be a realization of $A$, that is:
\begin{enumerate}\itemsep=0pt
\item[(1)] $\fh$ is a complex vector space of dimension $n+2$,
\item[(2)] $\Pi=\{\al_0,\dots,\al_n\}\subset \fh^\ast$ and $\Pi^\vee = \big\{\al_0^\vee,\dots,\al_n^\vee\big\}\subset \fh$ are linearly independent subsets,
\item[(3)] $a_{i,j} = \bigl\langle\al_j,\al^\vee_i\bigr\rangle$ for $0\leq i,j\leq n$.
\end{enumerate}
Here, $\scal{\cdot}{\cdot}\colon \fh^\ast\times \fh\to \mathbb{C}$ denotes the pairing $\scal{\al}{h} = \al(h)$.
Fix an element $d\in \fh$ such that $\scal{\al_i}{d} = \delta_{0,i}$ for all $i\in I$ so that $\Pi^\vee\cup \{d\}$ is a basis of $\fh$.
We denote by $\fg$ the affine Lie algebra associated to this datum and refer to~\cite{KacB} for its definition.

Let $\La_0,\dots,\La_{n}\subset \fh^\ast$ be such that
\[
\La_i\bigl(\al^\vee_j\bigr) = \delta_{i,j} \qquad\text{and}\qquad \La_i(d) = 0\quad \text{ for all $i,j\in I$.}
\]
We set $\delta = \sum_{j=0}^{n} a_j \al_j$ so that
\[
\delta\bigl(\al_i^\vee\bigr) = \sum_{j=0}^n a_j \al_j\bigl(\al_i^\vee\bigr)=\sum_{j=0}^n a_j a_{i,j} = [Av]_{i} =0\qquad \text{and} \qquad \delta(d) = a_0.
\]
Then the family $(\La_0,\dots,\La_n,\delta)$ is the dual basis of $\bigl(\al^\vee_0,\dots,\al^\vee_n,d\bigr)\subset \fh$ and we have
\[
\al_j = \sum_{i\in I} \bigl\langle\al_j,\al_i^\vee\bigr\rangle \La_i + \scal{\al_j}{d}\delta = \sum_{i\in I} a_{i,j} \La_i \quad \text{for all $j\in I^\ast$.}
\]

\smallskip

For any $i\in I^{\ast}$, set%
\[
\omega_{i}=\Lambda_{i}-\frac{a_{i}^{\vee}}{a_{0}^{\vee}}\Lambda_{0}.
\]

Recall the following classical lemma.

\begin{Lemma}
	For any $j\in I^{\ast}$, we have
	\[
	\alpha_{j}=\sum_{i\in I^{\ast}}a_{i,j}\omega_{i},%
	\]
	that is, the weights $\omega_{i}$, $i\in I^{\ast}$ and the roots $\alpha_{i}$, $i\in
	I^{\ast}$ can be regarded as the dominant weights and the simple roots of the
	finite root system with Cartan matrix $(a_{i,j})_{(i,j)\in I^{\ast}\times
		I^{\ast}}$.
\end{Lemma}

\begin{proof}
For any $j\in I^{\ast}$, we have%
\begin{align*}
\alpha_{j}=\sum_{i\in I}a_{i,j}\Lambda_{i}&{}=\sum_{i\in I^{\ast}}a_{i,j}%
\omega_{i}+\Biggl( a_{0,j}+\sum_{i\in I^{\ast}}\frac{a_{i,j}a_{i}^{\vee}%
}{a_{0}^{\vee}}\Biggr) \Lambda_{0}\\
&{}=
\sum_{i\in I^{\ast}}a_{i,j}\omega_{i}+\frac{1}{a_{0}^{\vee}}\Biggl(
a_{0,j}a_{0}^{\vee}+\sum_{i\in I^{\ast}}a_{i,j}a_{i}^{\vee}\Biggr)
\Lambda_{0}.
\end{align*}
But now, by definition of $v^{\vee}$, we have%
\[
\sum_{i\in I}a_{i}^{\vee}a_{i,j}=0
\]
for any $j\in J^{\ast}$ hence the expected result.
\end{proof}

Let us set%
\[
\eta^{\vee}:=\sum_{i\in I}\frac{a_{i}^{\vee}}{a_{0}^{\vee}}.
\]

There exists an invariant non degenerate symmetric bilinear form $(\cdot,\cdot)$ on $\fh$, uniquely defined~by
\begin{gather*}
\bigl(\al_i^\vee,\al^\vee_j\bigr) = \bigl\langle\al_i,\al_j^\vee\bigr\rangle a_i{a_i^\vee}^{-1}=a_j{a_j^\vee}^{-1}a_{j,i}\quad\mbox{for $i,j\in I$},\\
\bigl(\al_i^\vee,d\bigr) = 0 \quad \mbox{for $i\in I^\ast$},\qquad
\bigl(\al_0^\vee,d\bigr) = a_0 \quad \mbox{for $i\in I^\ast$},\qquad
(d,d) = 0.
\end{gather*}
It can be checked that $\bigl(\al_i^\vee,\al_j^\vee\bigr) = \bigl(\al_j^\vee,\al_i^\vee\bigr)$ for all $i,j\in I$. Let $\nu$ be the associated map from $\fh$ to its dual:
\begin{align*}
\nu \colon \ \fh  \to  \fh^\ast,\qquad
 h \mapsto  (h,\cdot).
\end{align*}
The form $(\cdot,\cdot)$ on $\fh$ then induces a form on $\fh^\ast$ via $\nu$. We still denote this form $(\cdot,\cdot)$. Then we have $\nu(d) =\La_0$, \smash{$a_{i,j} = \frac{2(\al_i,\al_j)}{(\al_i,\al_i)}$} and \smash{$\nu\bigl(\al_i^{\vee}\bigr) = \frac{2\al_i}{(\al_i,\al_i)}$} and
\begin{gather*}
(\al_i,\al_j) = a_i^\vee{a_i}^{-1}a_{i,j}\quad\mbox{for $i,j\in I$},\qquad
(\al_i,\La_0) = 0 \quad \mbox{for $i\in I^\ast$},\\
(\al_0,\La_0) = a_0^{-1} ,\qquad
(\La_0,\La_0) = 0,\qquad
(\La_0,\delta)=1.
\end{gather*}

\textbf{Warning.} The previous scalar product does not yield in
general the Euclidean norm on the real part of $\fh^{\ast}$ as it appears in
many references (see, for example,~\cite{Bourbaki}). For example, in type~$C_{n}$
the long roots have length equal to $2$ and the short roots length $1$ which
is not the most common convention. We indicate in the table below the relation
between the previous norm $ \Vert \beta \Vert ^{2}=(\beta,\beta)$ and
the Euclidean norm $ \Vert \beta \Vert _{2}^{2}$ used in~\cite{Bourbaki}.

\begin{table}[h]\centering \renewcommand{\arraystretch}{1.32}\caption{Affine types and associated datum.}
\begin{tabular}{|c|c|c|c|c|c|c|c|}
\hline
type & $v=(a_{0},\dots ,a_{n})$ & $v^{\vee }=\bigl(a_{0}^{\vee },\dots
,a_{n}^{\vee }\bigr)$ & $Q$ & $Q^{\vee }$ & $\eta $ & $\eta ^{\vee }$ & $%
 \Vert \beta  \Vert ^{2}$ \\ \hline
$A_{n}^{(1)}$ & $1^{n+1}$ & $1^{n+1}$ & $A_{n}$ & $A_{n}$ & $n+1$ & $n+1$ & $%
 \Vert \beta  \Vert _{2}^{2}$ \\ \hline
$B_{n}^{(1)}$ & $1^{2}2^{n-1}$ & $1^{2}2^{n-2}1$ & $B_{n}$ & $C_{n}$ & $2n$
& $2n-1$ & $ \Vert \beta  \Vert _{2}^{2}$ \\ \hline
$C_{n}^{(1)}$ & $12^{n-1}1$ & $1^{n+1}$ & $C_{n}$ & $B_{n}$ & $2n$ & $n+1$ &
$\frac{1}{2} \Vert \beta  \Vert _{2}^{2}$ \\ \hline
$D_{n}^{(1)}$ & $1^{2}2^{n-3}1^{2}$ & $1^{2}2^{n-3}1^{2}$ & $D_{n}$ & $D_{n}$
& $2n-2$ & $2n-2$ & $ \Vert \beta  \Vert _{2}^{2}$ \\ \hline
$A_{2n-1}^{(2)}$ & $1^{2}2^{n-2}1$ & $1^{2}2^{n-1}$ & $C_{n}$ & $B_{n}$ & $%
2n-1$ & $2n$ & $ \Vert \beta  \Vert _{2}^{2}$ \\ \hline
$A_{2n}^{(2)}$ & $2^{n}1$ & $12^{n}$ & $C_{n}$ & $B_{n}$ & $2n+1$ & $2n+1$ &
$ \Vert \beta  \Vert _{2}^{2}$ \\ \hline
$D_{n+1}^{(2)}$ & $1^{n+1}$ & $12^{n-1}1$ & $B_{n}$ & $C_{n}$ & $n+1$ & $2n$
& $2 \Vert \beta  \Vert _{2}^{2}$ \\ \hline
$G_{2}^{(1)}$ & $123$ & $121$ & $G_{2}$ & $G_{2}^{t}$ & $6$ & $4$ & $\frac{1%
}{3} \Vert \beta  \Vert _{2}^{2}$ \\ \hline
$D_{4}^{(3)}$ & $121$ & $123$ & $G_{2}^{t}$ & $G_{2}$ & $4$ & $6$ & $%
 \Vert \beta  \Vert _{2}^{2}$ \\ \hline
\end{tabular}
\end{table}

There are different objects associated to the datum $\bigl(\fh,\Pi,\Pi^\vee\bigr)$, some of them lives in $\fh$ other in $\fh^\ast$. We will as much as possible use the following convention: we will add the suffix ``co'' to the name of the object to indicate that it naturally lives in $\fh$ (even though we sometime think of it as an element of $\fh^\ast$) and we will add a superscript $^\vee$ to the notation. For instance, $\al_i^\vee$ is called a coroot as it is an element of~$\fh$. We now introduce various lattices:
\begin{itemize}
\item  the coweight lattice: $P_a^\vee := \bigoplus_{i\in I} \Z \al_i^\vee +\frac{1}{a_0}\Z d\subset \fh$,
\item  the weight lattice $P_a = \bigl\{\ga\in \fh^\ast\mid \ga\bigl(P_a^\vee\bigr) \subset \Z\bigr\} = \bigoplus_{i\in I} \Z\La_i +\frac{1}{a_0} \Z\delta$,
 \item  the dominant weights in $\fh^\ast$ are the elements of $P^+_a = \bigoplus_{i\in I} \Z_{\geq 0}\La_i +\frac{1}{a_0} \Z\delta$,
\item  the root lattice $Q_a=\bigoplus_{i\in I} \Z \al_i$,
\item  the coroot lattice $Q_a^\vee=\bigoplus_{i\in I} \Z \al^\vee_i$.
\end{itemize}
For any $i\in I$, we define the simple reflection $s_{i}$ on $\fh^\ast$ by
\begin{equation*}
s_{i}(x)=x-\bigl\langle\alpha_{i}^{\vee},x\bigr\rangle\alpha_{i}\quad \text{for any
}x\in \fh^\ast . \label{defSi}
\end{equation*}
The affine Weyl group $W_{a}$ is the subgroup of $\mathrm{GL}(\mathfrak{h}^{\ast})$ generated
by the reflections $s_{i}$. Since $(\La_0,\dots,\La_n,\delta)$ is the dual basis of $\bigl(\al_0^\vee,\dots,\al_n^\vee,d\bigr)$, we have for all $i\in I$
\begin{align*}
s_i(\delta)=\delta - \bigl\langle\alpha_{i}^{\vee},\delta\bigr\rangle \al_i = \delta,\qquad
s_i(\La_j)= \La_j \quad\text{if $i\neq j$},\qquad
s_i(\al_i)= -\al_i.
\end{align*}
The Weyl group $W_a$ is acting on the weight
lattice $P_{a}$.

Let $\mathring{A}$
be the matrix obtained from $A$ by deleting the row and the column corresponding to~$0$. Then it is well-known that $\mathring{A}$ is a Cartan matrix of finite type. Let \smash{$\mathring{\fh}{}^\ast$} and $\mathring{\fh}$ be the vector spaces spanned by the subset $\mathring{\Pi} = \Pi\setminus \{\al_0\}$ and $\mathring{\Pi}{}^\vee = \Pi\setminus \bigl\{\al^\vee_0\bigr\}$. The root and weight lattices associated to $\mathring{A}$ are $Q=\bigoplus_{i\in I^\ast}\mathbb{Z}\alpha_{i}$ and $P=\bigoplus_{i\in I^\ast}\mathbb{Z}\omega_{i}$ where we have set $\omega_{i}=\Lambda_{i}-a_{i}^{\vee}\Lambda_{0}$ for all $i\in I^\ast$. We denote by $W$ the finite Weyl group associated to $\mathring{A}$: it is generated by the orthogonal reflections $s_i$ with respect to the hyperplane orthogonal to $\al_i$ in $\mathring{\fh}{}^\ast$. Finally, we set $Q^\vee=\bigoplus_{i\in I^\ast}\mathbb{Z}\alpha^\vee_{i}\subset \mathring{\fh}{}^\ast$. Note that the reflection $s_i\in W_a$ for $i\in I^\ast$ stabilizes~$\mathring{\fh}{}^\ast$ so that the Weyl group $W$ can be seen as the subgroup of $W_a$ generated by $(s_i)_{i\in I^\ast}$.

Let
$
\theta=\delta-a_0\alpha_0$.

\textbf{Warning.} The root $\theta$ does not always coincide with the highest root of the finite root system $\mathring{\Pi} = \{\al_1,\dots,\al_n\}$. This is, for example, the case in type \smash{$A_{2n-1}^{(2)}$}. We refer the reader to~\cite[Proposition~17.18]{carter_affine} for more details.

Let $s_\theta$ be the orthogonal reflection with respect to $\theta$ defined by $s_\theta(\ga) = \ga - \scal{\theta^\vee}{\ga}\theta$. For all $\ga\in \fh^\ast$, we have
\[
s_0s_{\theta}(\ga)= \gamma +(\ga,\delta)\theta - \left((\ga,\theta)+\frac{1}{2}(\theta,\theta)(\ga,\delta)\right)\delta.
\]
Consider $\beta\in \mathring{\fh}{}^\ast$. We define the map $t_\beta$ on $\fh^\ast$ by
\[
t_{\beta}(\gamma) = \gamma +(\gamma,\delta)\beta - \left((\gamma,\beta)+\frac{1}{2}(\beta,\beta)(\gamma,\delta)\right)\delta
\]
for all $\gamma$ in $\fh^\ast$.
Note that if $(\ga,\delta) = 0$ (as it is the case when $\ga\in \mathring{\fh}{}^\ast$), we get a simpler formula
\[
t_{\beta}(\gamma) = \gamma - (\gamma,\beta)\delta.
\]
It is important to observe that $t_{\beta}$ does not act as a translation on $\mathring{\fh}{}^\ast$.
Nevertheless, it can be shown that
\begin{itemize}\itemsep=0pt
\item  $t_\beta\circ t_{\beta'}=t_{\beta+\beta'}$ for all $\beta,\beta'\in \fh^\ast$,
\item  $w\circ t_{\beta} \circ w^{-1} = t_{w(\beta)}$ for all $\beta\in \fh^\ast$ and $w\in W$,
\item  $s_{0}=t_{\theta}s_{\theta}=s_{\theta}t_{-\theta}$.
\end{itemize}
These relations tell us that $W_a$ is the semi-direct product of the finite Weyl group $W$ with the lattice that is generated by the $W$-orbit of $\frac{1}{a_0}\theta$. Let us denote by $M^{*}$ this lattice. So that we have
$W_a\simeq W\ltimes M^{*}$.
The table below gives the lattice $M^{*}$ expressed as a sublattice of the underlying finite root system $Q$ thanks to the simple roots (observe that $M^{*}$ is a sublattice of~$Q$ by definition).
\begin{table}[th]\centering\renewcommand{\arraystretch}{1.32}
\caption{Table of affine types with their corresponding underlying type of finite root system an~$M^{\ast}$.}\label{table2}
\begin{tabular}{|c|c|c|}
\hline
type & $Q$ & $M^{\ast }$ \\ \hline
$A_{n}^{(1)}$ & $A_{n}$ & ${\textstyle\bigoplus\limits_{i=1}^{n}}\mathbb{%
	Z\alpha }_{i}$ \\ \hline
$B_{n}^{(1)}$ & $B_{n}$ & ${\textstyle\bigoplus\limits_{i=1}^{n-1}}\mathbb{%
	Z\alpha }_{i}\oplus 2\mathbb{Z}\alpha _{n}$ \\ \hline
$C_{n}^{(1)}$ & $C_{n}$ & ${\textstyle\bigoplus\limits_{i=1}^{n-1}}2\mathbb{%
	Z\alpha }_{i}\oplus \mathbb{Z}\alpha _{n}$ \\ \hline
$D_{n}^{(1)}$ & $D_{n}$ & ${\textstyle\bigoplus\limits_{i=1}^{n}}\mathbb{%
	Z\alpha }_{i}$ \\ \hline
$A_{2n-1}^{(2)}$ & $C_{n}$ & ${\textstyle\bigoplus\limits_{i=1}^{n}}\mathbb{%
	Z\alpha }_{i}$ \\ \hline
$A_{2n}^{(2)}$ & $C_{n}$ & ${\textstyle\bigoplus\limits_{i=1}^{n-1}}\mathbb{%
	Z\alpha }_{i}\oplus \frac{1}{2}\mathbb{Z}\alpha _{n}$ \\ \hline
$D_{n+1}^{(2)}$ & $B_{n}$ & ${\textstyle\bigoplus\limits_{i=1}^{n}}\mathbb{%
	Z\alpha }_{i}$ \\ \hline
$G_{2}^{(1)}$ & $G_{2}$ & ${\textstyle}\mathbb{Z\alpha }_{1}\oplus 3\mathbb{%
	Z\alpha }_{2}$ \\ \hline
$D_{4}^{(3)}$ & $G_{2}^{t}$ & ${\textstyle}\mathbb{Z\alpha }_{1}\oplus
\mathbb{Z\alpha }_{2}$ \\ \hline
\end{tabular}
\end{table}

\section{Macdonald formula types}\label{sec:3}

In this section, we examine how the decomposition of $W_a$ as a semi-direct product provides a~connection with the virtual (non-affine) characters.

We now consider the Macdonald identity from the point of view of the Weyl--Kac
denominator formula. The computations of this section are essentially the same
as in Section 20 of Carter's book~\cite{carter_affine}. For any affine root
system, the Weyl denominator formula can be written
\begin{equation}
\prod_{\alpha\in R_{+}}(1-{\rm e}^{-\alpha})^{m_{\alpha}}=\sum_{w\in W_{a}%
}\varepsilon(w){\rm e}^{w(\rho)-\rho},\label{WKDF}%
\end{equation}
where $\rho$ is any element of $P_{a}$ such that $\bigl\langle\rho,\al_i^\vee\bigr\rangle=1$
for any $i\in I$. Here the numbers $m_{\alpha}$ are the multiplicity of the
positive roots of the affine root system considered. Contrary to the classical
Weyl denominator formula (for a non affine finite root system), $\rho$ is not
unique here and cannot be defined as the half sum of the positive roots.
Nevertheless, with the previous notation, we can take%
\[
\rho=\sum_{i\in I}\Lambda_{i}=\eta^{\vee}\Lambda_{0}+\mathring{\rho},%
\]
where
\[
\mathring{\rho}=\frac{1}{2}\sum_{\alpha\in\mathring{R}_{+}}\alpha=\sum_{i \in
	I^{*}}\omega_{i}%
\]
is the half sum of the positive roots associated to the underlying finite root
system. Recall the following notation: for any classical weight \smash{$\gamma
\in\mathring{P}$}%
\[
a_{\gamma}=\sum_{u\in W}\varepsilon(w){\rm e}^{w(\gamma)}.
\]
In particular,
\[
a_{\mathring{\rho}}=\sum_{u\in W}\varepsilon(w){\rm e}^{w(\mathring{\rho})}={\rm e}^{
	\mathring{\rho}}\prod_{\alpha\in\mathring{R}_{+}}(1-{\rm e}^{-\alpha})
\]
from the Weyl denominator formula for the non affine root system. From the
(usual, that is, non affine) Weyl Kac character formula, we can define the
virtual characters%
\[
s_{\gamma}=\frac{a_{\gamma+\mathring{\rho}}}{a_{\mathring{\rho}}},\qquad\gamma
\in\mathring{P}.
\]
Then $s_{\gamma}=0$ or there exists a unique dominant weight $\lambda\in
P_{+}$ and an element $u\in W$ such that
\begin{equation}
s_{\gamma}=\varepsilon(u)s_{\lambda}\label{SR}%
\end{equation}
with $\lambda=u(\gamma+\mathring{\rho})-\mathring{\rho}=u\circ\gamma$ (the
so-called ``dot''-action). Recall also that each element $w$ of the affine Weyl
group $W_{a}$ admits a unique decomposition of the form%
\[
w=t_{\gamma}u \quad \text{with }\gamma\in M^{\ast} \ \text{and}\ u\in W.
\]

We will now compute the right-hand side of (\ref{WKDF}) by using the previous
decomposition of the elements in $W_{a}$. In fact, we will compute the
convenient renormalization below:%
\[
\Delta=\prod_{\alpha\in R_{+}\setminus\mathring{R}_{+}}(1-{\rm e}^{-\alpha
})^{m_{\alpha}}=\frac{1}{{\rm e}^{-\mathring{\rho}}a_{\mathring{\rho}}}\sum_{w\in
	W_{a}}\varepsilon(w){\rm e}^{w(\rho)-\rho}.
\]
First, fix $w=t_{\gamma}u$ in $W_{a}$. We have%
\begin{align*}
w(\rho)-\rho&{}=t_{\gamma}u\bigl(\eta^{\vee}\Lambda_{0}+\mathring{\rho}\bigr)-\eta^{\vee}\Lambda_{0}-\mathring{\rho}
 =t_{\gamma}\bigl( \eta^{\vee}\Lambda_{0}+u(\mathring{\rho})\bigr) -\eta^{\vee}\Lambda_{0}-\mathring{\rho}\\
&{}=\eta^{\vee}( t_{\gamma}(\Lambda_{0})-\Lambda_{0}) +t_{\gamma}u(\mathring{\rho})-\mathring{\rho},
\end{align*}
where the first equality uses the fact that $u(\Lambda_{0})=\Lambda_{0}$ for
any element $u$ in $W$ the finite Weyl group. We can now apply the formula
recalled in the previous section giving the action of a~translation
$t_{\gamma}$ on the element of $P_{a}$. We get%
\[
t_{\gamma}(\Lambda_{0})-\Lambda_{0}=\Lambda_{0}+\gamma-\frac{1}{2} \Vert
\gamma \Vert ^{2}\delta-\Lambda_{0}=\gamma-\frac{1}{2} \Vert
\gamma \Vert ^{2}\delta
\]
and also%
\[
t_{\gamma}u(\mathring{\rho})-\mathring{\rho}=u(\mathring{\rho})-(\gamma
,u(\mathring{\rho}))\delta-\mathring{\rho}=u(\mathring{\rho})-\mathring{\rho
}-\bigl(u^{-1}(\gamma),\mathring{\rho}\bigr)\delta.
\]
This yields
\[
w(\rho)-\rho=\eta^{\vee}\gamma+u(\mathring{\rho})-\mathring{\rho}-\biggl(
\frac{\eta^{\vee}}{2}\Vert \gamma\Vert ^{2}+\bigl(u^{-1}(\gamma
),\mathring{\rho}\bigr)\biggr) \delta
\]
and by setting $q={\rm e}^{-\delta}$%
\[
\Delta=\frac{1}{{\rm e}^{-\mathring{\rho}}a_{\mathring{\rho}}}\sum_{\gamma\in
	M^{\ast}}\sum_{u\in W}\varepsilon(t_{\gamma})\varepsilon(u)q^{\frac{\eta
		^{\vee}}{2}\Vert \gamma\Vert ^{2}+(u^{-1}(\gamma),\mathring{\rho}%
	)}{\rm e}^{\eta^{\vee}\gamma+u(\mathring{\rho})-\mathring{\rho}}.
\]
Now, since the classical root lattice $M^{\ast}$ is $W$-invariant, we can set
$\beta=u^{-1}(\gamma)\in M^{\ast}$ in the previous expression and obtain%
\[
\Delta=\frac{1}{{\rm e}^{-\mathring{\rho}}a_{\mathring{\rho}}}\sum_{\beta\in
	M^{\ast}}\sum_{u\in W}\varepsilon(u)q^{\frac{\eta^{\vee}}{2} \Vert
	\beta \Vert ^{2}+(\beta,\mathring{\rho})}{\rm e}^{u(\eta^{\vee}\beta
	)+u(\mathring{\rho})-\mathring{\rho}}%
\]
because $\varepsilon(t_{\beta})=\varepsilon(t_{\gamma})=1$ (the signature of
any translation in $W_{a}$ is equal to $1$) and $ \Vert \beta \Vert
^{2}=\Vert \gamma\Vert ^{2}$. This can be rewritten%
\begin{align}
\Delta &{}=\frac{1}{{\rm e}^{-\mathring{\rho}}a_{\mathring{\rho}}}\sum_{\beta\in
	M^{\ast}}q^{\frac{\eta^{\vee}}{2} \Vert \beta \Vert ^{2}%
	+(\beta,\mathring{\rho})}\sum_{u\in W}\varepsilon(u){\rm e}^{u(\eta^{\vee}%
	\beta+\mathring{\rho})-\mathring{\rho}}\nonumber\\
&{}=\sum_{\beta\in M^{\ast}}q^{\frac{\eta^{\vee}}{2} \Vert \beta
	 \Vert ^{2}+(\beta,\mathring{\rho})}\frac{a_{\eta^{\vee}\beta}%
}{a_{\mathring{\rho}}}=\sum_{\beta\in M^{\ast}}q^{\frac{\eta^{\vee}}%
	{2} \Vert \beta \Vert ^{2}+(\beta,\mathring{\rho})}s_{\eta^{\vee
	}\beta}.\label{Detal1}
\end{align}
In this last formula, we use an indexation by the affine translations instead
of the affine Grassmannians elements (which are the minimal length
representatives of the left cosets in $W_{a}/W$). Observe that each coset in
$W_{a}/W$ also contains a unique translation $t_{\beta}$ but it is not of
minimal length in general. In fact, we will go further in the following
section by applying the straightening rules for the virtual characters
(\ref{SR}) and also associating to each translation $\beta$, the unique affine Grassmannian element $c(\beta)$ such that $c(\beta)$ is of minimal length in the coset $t_{\beta}W$.
Observe that $t_{\beta}(\Lambda_{0})=c(\beta)(\Lambda_{0})$. In the
following, it will be convenient to rewrite our formula (\ref{Detal1}) by
changing $\beta$ into $-\beta$ (which is clearly possible since $M^{\ast}$ is
a lattice). Since $ \Vert \beta \Vert ^{2}= \Vert -\beta
 \Vert ^{2}$, this gives the expression%
\begin{equation}
\Delta=\sum_{\beta\in M^{\ast}}q^{\frac{\eta^{\vee}}{2} \Vert
	\beta \Vert ^{2}-(\beta,\mathring{\rho})}s_{-\eta^{\vee}\beta
}.\label{Delta_Beta}%
\end{equation}

\begin{Remark}
	The previous expansion of $\Delta$ can be regarded as an analogue in affine type of some classical formulas due to Littlewood (see~\cite[p.~79]{Macbook}). We will study this analogy more deeply in a future work (see~\cite{LW2}).
\end{Remark}

\section{Connection with the affine Grassmannian}\label{sec:4}

Given the results of the previous section, we now examine how the Macdonald identity can be rewritten as a signed $q$-series of
	ordinary Weyl characters running over the affine Grassmannian elements and how it relates to the so-called atomic length as introduced by Chapelier-Laget and Gerber in~\cite{CGT}.
	
As already mentioned, the affine Grassmannian is the set of minimal length
elements in the left cosets of $W_{a}/W$. We shall denote by $W_{a}^{0}$ these
affine Grassmannian elements. Let us recall that any $w=t_{\beta}u$ in $W_{a}$
with $\beta\in M^{\ast}$ and $u\in W$ can be also written
\[
w=t_{\beta}u=u\bigl(u^{-1}t_{\beta}u\bigr)=ut_{u^{-1}(\beta)}.
\]
Therefore, we can use decompositions of $w$ of both forms $t_{\beta}u$ or
$ut_{\gamma}$. The second one is particularly well-adapted for computing the
length of $w$ since we have the formula%
\[
\ell(ut_{\gamma})=\sum_{\alpha\in\mathring{R}_{+}}\vert (\gamma
,\alpha)+\chi(u(\alpha)) \vert \text{ with }\chi(\alpha)=
\begin{cases}
0&\text{on }\mathring{R}_{+},\\
1&\text{on } {-}\mathring{R}_{+}.
\end{cases}
\]
From this, it is not difficult to check that a translation $t_{\gamma}$
belongs to $W_{a}^{0}$ if and only if ${\gamma\in M^{\ast}\cap(-P_{+})}$, that
is an antidominant weight for the finite root system. Indeed, we then have~${(\gamma,\alpha)\leq0}$ whereas~${\chi(w(\alpha))\geq0}$. We can in fact
completely characterize the elements in~$W_{a}^{0}$ by the following lemma
\big(which generalizes the previous observation on the translations~in~$W_{a}^{0}$\big).

\begin{Lemma}
	The element $w$ belongs to $W_{a}^{0}$ if and only if its admits a
	decomposition $w=ut_{\nu}$ such that $\nu\in M^{\ast}\cap(-P_{+})$ and $ u\in
	W$ is of minimal length in a left coset of $W/W_{\nu}$, where $W_{\nu}$ is the
	stabilizer of $\nu$ under the action of the finite Weyl group.
\end{Lemma}

Now let us consider a translation $t_{\beta}$ with $\beta\in M^{\ast}$ as in
(\ref{Delta_Beta}). In general, we do not have~$t_{\beta}$ in $W_{a}^{0}$ but
this can be corrected. To do this, observe that $W\cdot\beta$, the orbit of
$\beta$ under the~action of the finite Weyl group $W$, intersects $-P_{+}$ in
a unique antidominant weight $\nu$ with~${-\nu\in P_{+}}$. In general, one
rather uses the intersection with $P_{+}$, this works similarly with $-P_{+}$
just by composing with $w_{0}$, the maximal length element in $W$. Write as
usual $W^{\nu}$ for the set of minimal length elements in the cosets of
$W/W_{\nu}$. Then, there exists a unique $u\in W^{\nu}$ such that~${\beta=u(\nu)}$. We then have
\[
t_{\beta}=(ut_{\nu})u^{-1}=cu^{-1}%
\]
with $c=ut_{\nu}\in W_{a}^{0}$ and $u^{-1}\in W$. Since $\beta=u(\nu)$ and
$t_{\beta}=(ut_{\nu})u^{-1}$, we get
\[
 \Vert \beta \Vert ^{2}= \Vert \nu \Vert ^{2}\qquad\text{and}\qquad \varepsilon(t_{\beta})=\varepsilon(t_{\nu})=1.
\]
We can also compute the analogue of the number of boxes in a core from%
\[
\Lambda_{0}-c(\Lambda_{0})=\Lambda_{0}-u\left( \Lambda_{0}+\nu-\frac{1}%
{2} \Vert \nu \Vert ^{2}\delta\right) =\Lambda_{0}-\Lambda_{0}%
-u(\nu)+\frac{1}{2} \Vert \nu \Vert ^{2}\delta=\frac{1}{2} \Vert
\beta \Vert ^{2}\delta-\beta
\]
by setting
\begin{align*}
L^{^{\vee}}(c)&=(\Lambda_{0}-c(\Lambda_{0}),\rho) =\eta^{\vee}(
\Lambda_{0}-ut_{\nu}(\Lambda_{0}),\Lambda_{0}) +( \Lambda
_{0}-ut_{\nu}(\Lambda_{0}),\mathring{\rho}) \\
& =\eta^{\vee}\left( \frac{1}{2} \Vert \beta \Vert ^{2}\delta
-\beta,\Lambda_{0}\right) +\left( \frac{1}{2} \Vert \beta \Vert
^{2}\delta-\beta,\mathring{\rho}\right)
\end{align*}
and by using the equalities $(\Lambda_{0},\delta)=1$, $(\delta,\mathring{\rho
})=0$ and $(\Lambda_{0},\beta)=0$, we obtain
\[
L^{^{\vee}}(c)=\frac{\eta^{\vee}}{2} \Vert \beta \Vert ^{2}%
-(\beta,\mathring{\rho})\quad \text{with }c=ut_{\nu}=t_{\beta}u.
\]

\begin{Remark}
	\label{rem:reseau} Observe that the definition of $L^{^{\vee}}(c)$ is very
	close to the so-called ``Atomic Length'' introduced by Chapelier-Laget and
	Gerber and generalizing the number of boxes in a core partition (see
	\cite[Corollary~8.2]{CGT}) which satisfies
	\[
	L(c)=\frac{\eta}{2} \Vert \beta \Vert ^{2}-\bigl(\beta,\mathring{\rho}^{\vee}\bigr).
	\]
	In particular, $L(c)$ counts the number of simple roots appearing in the
	decomposition of the weight $\Lambda_{0}-c(\Lambda_{0})$ on the basis of
	simple roots. Therefore, we have $L^{^{\vee}}(c)=L(c)$ in the simply laced
	cases. More generally, one can observe that the lattice $M^{\ast}$ and the
	norm $ \Vert \cdot \Vert ^{2}$ are~the~same in type~\smash{$B_{n}^{(1)}$} and
	\smash{$A_{2n-1}^{(2)}=\bigl(B_{n}^{(1)}\bigr)^{t}$}. For types \smash{$C_{n}^{(1)}$} and \smash{$D_{n+1}^{(2)}=\bigl(C_{n}^{(1)}\bigr)^{t}$},
    we only have
\[
M_{C_{n}^{(1)}}^{\ast}=2M_{D_{n+1}^{(2)}}^{\ast} \qquad \text{but} \quad \Vert \cdot \Vert _{C_{n}^{(1)}}^{2}=\frac{1}{4} \Vert \cdot \Vert _{D_{n+1}^{(2)}}^{2}.
\]
 This
	allows us to conclude that the statistics $L$ and~$L^{^{\vee}}$ take exactly
	the same values up to transposition of the Cartan matrices.
\end{Remark}

To rewrite $\Delta$ in terms of the elements in $W_{a}^{0}$, observe that the
previous construction gives a~bijection
\begin{equation}
M^{\ast}\rightarrow W_{a}^{0},\qquad
\beta\longmapsto c=ut_{\nu}=t_{\beta}u
 \label{betac}
\end{equation}
such that $t_{\beta}=ut_{\nu}u^{-1}=cu^{-1}$. We so have $1=\varepsilon
(t_{\beta})=\varepsilon(c)\varepsilon\bigl(u^{-1}\bigr)$. Now for each term in the sum
(\ref{Delta_Beta}), we get by the previous computations and remarks the
equality%
\begin{align*}
q^{\frac{\eta^{\vee}}{2} \Vert \beta \Vert ^{2}-(\beta,\mathring
	{\rho})}s_{-\eta^{\vee}\beta}&=\varepsilon(c)\varepsilon\bigl(u^{-1}\bigr)q^{L^{^{\vee}%
	}(c)}s_{-\eta^{\vee}u(\nu)}=\varepsilon(c)q^{L^{^{\vee}}(c)}s_{u^{-1}(
	-\eta^{\vee}u(\nu)+\mathring{\rho}) -\mathring{\rho}}
\\	&=\varepsilon
(c)q^{L^{^{\vee}}(c)}s_{-\eta^{\vee}\nu+u^{-1}(\mathring{\rho})-\mathring
	{\rho}},%
\end{align*}
where $-\eta^{\vee}\nu\in P_{+}$. This gives the theorem below.

\begin{Theorem}
	\label{thm:macdo} With the previous notation, the renormalized Weyl--Kac
	denominator formula can be written
	\[
	\Delta=\sum_{c\in W_{a}^{0}}\varepsilon(c)q^{L^{^{\vee}}(c)}s_{-\eta^{\vee}%
		\nu+u^{-1}(\mathring{\rho})-\mathring{\rho}},%
	\]
	where we set $c=ut_{\nu}$ with $\nu\in M^{\ast}\cap(-P_{+})$ and $u\in W^{\nu
	}$ for any element of the affine Grassmannian $W_{a}^{0}$. Moreover, each
	weight $-\eta^{\vee}\nu+u^{-1}(\mathring{\rho})-\mathring{\rho}$ belongs to
	$P_{+}$, the set of dominant weights for $\mathring{A}$.
\end{Theorem}

It remains to prove the last claim of the theorem. First observe that we
always have $M^{\ast}\subset P$ since $M^{\ast}\subset Q\subset P$. Thus
$-\eta^{\vee}\nu+u^{-1}(\mathring{\rho})-\mathring{\rho}$ belongs to $P$ and
it suffices to prove that it is dominant. This is done in the lemma below.

\begin{Lemma}
	For any dominant weight $\lambda\in M^{\ast}$ and any $u\in W^{\lambda}$
	$($i.e., of minimal length in the cosets of $W/W_{\lambda})$, the weight
	\[
	\eta^{\vee}\lambda+\bigl(u^{-1}(\mathring{\rho})-\mathring{\rho}\bigr)
	\]
	is dominant.
\end{Lemma}

\begin{proof}
	We need to prove that
	\[
	\bigl(\eta^{\vee}\lambda+\bigl(u^{-1}(\mathring{\rho})-\mathring{\rho}\bigr),\alpha_{i}%
	\bigr)=\bigl(\eta^{\vee}(\lambda,\alpha_{i})+\bigl(u^{-1}(\mathring{\rho})-\mathring{\rho
	}\bigr), \alpha_{i}\bigr)\geq0
	\]
	for any $i\in I^{\ast}$. Observe that we have
	\[
	\bigl(\bigl(u^{-1}(\mathring{\rho})-\mathring{\rho}\bigr), \alpha_{i}\bigr) =\bigl(u^{-1}(\mathring{\rho
	}),\alpha_{i}\bigr)-(\mathring{\rho},\alpha_{i})=(\mathring{\rho},u(\alpha
	_{i}))-\frac{1}{2} \Vert \alpha_{i} \Vert ^{2}
	\]
	since \smash{$(\mathring{\rho},\alpha_{i})=\frac{ \Vert \alpha_{i} \Vert
		^{2}}{2}\bigl(\mathring{\rho},\alpha_{i}^{\vee}\bigr)=\frac{ \Vert \alpha
		_{i} \Vert ^{2}}{2}$}. If $u(\alpha_{i})\in R_{+}$, then $(\mathring{\rho
	},u(\alpha_{i}))\geq\frac{1}{2} \Vert \alpha_{i} \Vert ^{2}$ and
	therefore $\bigl(u^{-1}(\mathring{\rho})-\mathring{\rho}\bigr),\alpha_{i})\geq0$. We are
	done because $\lambda$ is dominant and thus $\eta^{\vee}(\lambda,\alpha
	_{i})\geq0$. Now assume $u(\alpha_{i})\in-R_{+}$ is a negative root. Then, we
	know that there exists a reduced expression of $u$ ending by $s_{i}$, that is
	of the form $u=u^{\prime}s_{i}$ with $\ell(u)=\ell(u^{\prime})+1$. Since $u\in
	W^{\lambda}$, this implies that $(\lambda,\alpha_{i})\geq p$ where $p=1$ in
	all affine types except in types \smash{$C_{n}^{(1)}$} and \smash{$G_{2}^{(1)}$} where $p=2$
	and $p=3$, respectively (because $\lambda\in M^{\ast}\cap P_{+}$). Indeed, we
	would have otherwise $u(\lambda)=u^{\prime}s_{i}(\lambda)=u^{\prime}$ and~$u$~would
    not be of minimal of length. We thus have
	\[
	\bigl(\eta^{\vee}\lambda+\bigl(u^{-1}(\mathring{\rho})-\mathring{\rho}\bigr),\alpha_{i}\bigr)\geq
	p\eta^{\vee}-(\mathring{\rho},-u(\alpha_{i}))-\frac{1}{2} \Vert \alpha
	_{i} \Vert ^{2}\geq0,
	\]
	where $\alpha\in R_{+}$. One can then check that for any $\alpha\in R_{+}$,
	we have $(\mathring{\rho},-u(\alpha_{i}))+\frac{1}{2} \Vert \alpha
	_{i} \Vert ^{2}\leq p\eta^{\vee}$ and thus the desired inequality
	\[
	\bigl(\eta^{\vee}\lambda+\bigl(u^{-1}(\mathring{\rho})-\mathring{\rho}\bigr),\alpha_{i}%
	\bigr)\geq0.\qedhere
\tag*{\qed}
\]
\renewcommand{\qed}{}
\end{proof}

\begin{Remark}
\
\begin{itemize}\itemsep=0pt
\item Recall that the Euler product is $\prod_{k\geq1} \bigl(
		1-q^{k}\bigr) $. The Nekrasov--Okounkov formula \eqref{NOdebut} gives an
		expansion of powers of the Euler product for any complex $z$. A crucial step in
		Han's proof to derive \eqref{NOdebut} is a specialization of the Macdonald
		identity. For every semisimple Lie algebra $\mathfrak{g}$ of rank $n$,
		Macdonald~\cite{Mac} proves that by setting $q={\rm e}^{-\delta}$ and~${{\rm e}^{\alpha
			_{i}}\mapsto\pm1}$ for~${1\leq i\leq n}$
		, the left-hand side of Theorem \ref{thm:macdo} is equal to
		\[
		\delta(q)=\prod_{k\geq1} \bigl( 1-q^{k}\bigr) ^{\dim\mathfrak{g}}.
		\]
		Moreover, a variation of the Euler product, called the Dedekind $\eta$-function,
		is defined by
		\[
		\eta(q)=q^{1/24}\prod_{k\geq1} \bigl( 1-q^{k}\bigr),
		\]
		so that $q^{\dim\mathfrak{g}/24}\delta(q)=\eta(q)^{\dim\mathfrak{g}}$. The
		term $q^{\dim\mathfrak{g}/24}$ arises from the ``strange formula'' \mbox{\cite[p.~243]{FVries}}:
        $\Phi_{R}(\rho,\rho)=\dim\mathfrak{g}/24$, where
		$\Phi_{R}$ is the scalar product on $R$ induced by the Killing form on
		$\mathfrak{g}$. Therefore, one so gets an explicit connection between the powers
		of the Dedekind $\eta$-function and weights $-\eta^{\vee}\nu+u^{-1}(\mathring
		{\rho})-\mathring{\rho}$ associated with elements of the affine Grassmannian,
		which we do not detail here.

\item Theorem \ref{thm:macdo} is in fact a particular case of a more general result expressing the Weyl denominator formula associated to a root system in terms of the characters corresponding to one of its parabolic subroot systems that will be detailed and exploited elsewhere.
\end{itemize}	
\end{Remark}

\section{Atomic length and cores}\label{sec:atomic}

As seen in the previous section, the Weyl--Kac denominator formula relates to the affine Grassmannian elements, defined algebraically. In this section, we explain how they can be described in a purely combinatorial way in terms of core partitions. More precisely, we give a simple description of the affine
Grassmannian elements of the previous classical affine root systems (twisted
or not) with underlying finite root system of rank $n$ in terms of families
of $2n$-core partitions. We will also consider the affine root system of
type \smash{$G_{2}^{(1)}$} and \smash{$D_{4}^{(3)}$} where our description will use $6$%
-core partitions. For each of the previous affine types, we provide an
embedding of its associated root and weight lattices in a weight lattice of
affine type $A$. Incidentally, this yields an embedding of the
corresponding affine Weyl group in a group of affine permutations (i.e., a
Weyl group of affine type $A$). Similar results for the previous non-twisted affine types also appeared in~\cite{STW} where they are obtained by different techniques. We give in Figures~\ref{Dynkin-diag} and~\ref{Dynkin-diag2} the labelling of the affine
Dynkin diagrams that we shall use in the following.

\begin{figure}[th]\centering \begin{tabular}{rc}
$A_{1}^{(1)}$ &\dynkin[labels={0,1},
edge length=1.5cm,
extended] A[1]{1}\\
$A_{n}^{(1)}(n\geq2)$ &\dynkin[labels={0,1,2,n-1,n},
edge length=1.5cm,
extended] A[1]{}\\
$B_{n}^{(1)}(n\geq3)$ &\dynkin[labels={0,1,2,3,n-2,n-1,n},
edge length=1.5cm,
extended] B[1]{}\\
$C_{n}^{(1)}(n\geq2)$ &\dynkin[labels={0,1,2,n-2,n-1,n},
edge length=1.5cm,
extended] C[1]{}\\
$D_{n}^{(1)}(n\geq4)$ &\dynkin [labels={0,1,2,3,n-3,n-2,n-1,n},
edge length=1.5cm,
extended] D[1]{} \\

\end{tabular}\caption{The Dynkin diagrams of the extended simple root systems.}\label{Dynkin-diag}
\end{figure}

\begin{figure}[th]\begin{tabular}{rc}
$A_{2}^{(2)}$ &\dynkin[labels={0,1},
edge length=1.5cm,
extended] A[2]{2}\\
$A_{2n}^{(2)}(n\geq2)$ &\dynkin[labels={0,1,2,3,n-2,n-1,n},
edge length=1.5cm,
extended] A[2]{even}\\
$A_{2n-1}^{(2)}(n\geq3)$ &\dynkin[labels={0,1,2,3,4,n-2,n-1,n},
edge length=1.5cm,
extended] A[2]{odd}\\
$D_{n+1}^{(2)}(n\geq2)$ &\dynkin[labels={0,1,2,3,n-3,n-2,n-1,n},
edge length=1.5cm,
extended] D[2]{} \\

\end{tabular}\caption{The Dynkin diagrams of the twisted simple root systems.}\label{Dynkin-diag2}\end{figure}

The results presented in this section are mainly based on foldings of the affine Dynkin diagram of \smash{$A_{2n-1}^{(1)}$}. Figure~\ref{aaaaa} shows how one can get the Dynkin diagram of type \smash{$C_n^{(1)}$} by gluing nodes of the Dynkin diagram of type \smash{$A_{2n-1}^{(1)}$}.

\begin{figure}[th!]\centering
\begin{tabular}{ccc}
\dynkin[fold,edge length=1.4cm]A[1]{**.*****.**} & $\longrightarrow$ & \dynkin[extended,edge length=1.4cm] C[1]{}
\end{tabular}
\caption{The Dynkin diagram of type \smash{$C_n^{(1)}$} by gluing nodes of the Dynkin diagram of type~\smash{$A_{2n-1}^{(1)}$}.}\label{aaaaa}
\end{figure}

\subsection{Affine type A and core partitions}

A \textit{partition} $\lambda$ of a positive integer $n$ is a non-increasing sequence of positive integers $\lambda=(\lambda_1,\lambda_2,\dots,\lambda_\ell)$ such that $\lvert \lambda \rvert := \lambda_1+\lambda_2+\dots+\lambda_\ell = n$. The $\lambda_i$'s are the \textit{parts} of $\lambda$, the number~$\ell$ of parts being the \textit{length} of $\lambda$, denoted by $\ell(\lambda)$. The \textit{weight} of $\lambda$ is $\lvert \lambda \rvert$. For convenience, set~${\lambda_i=0}$ for~all~$i>\ell(\lambda)$.

Each partition can be represented by its Ferrers diagram, which consists of a finite collection of boxes arranged in left-justified rows, with the row lengths in non-increasing order. The \textit{Durfee square} of $\lambda$ is the maximal square fitting in the Ferrers diagram. Its diagonal, denoted by $D(\lambda)$, will be called the main diagonal of $\lambda$. Its size will be denoted $d=d_\lambda:=\max(s \mid \lambda_s\geq s)$.
 The~partition $\lambda^{\mathrm{tr}}=\bigl(\lambda_1^{\mathrm{tr}},\lambda_2^{\mathrm{tr}},\dots,\lambda_{\lambda_1}^{\mathrm{tr}}\bigr)$ is the \textit{conjugate} of $\lambda$, where $\lambda_j^{\mathrm{tr}}$ denotes the number of boxes in the column $j$.

 For each box $s$ in the Ferrers diagram of a partition $\lambda$ (for short we will say for each box~$s$ in~$\lambda$), one defines the \textit{arm-length} (respectively \textit{leg-length}) as the number of boxes in the same row (respectively in the same column) as $s$ strictly to the right of (respectively strictly below) the box~$s$. The \textit{hook length} of $s$, denoted by $h_s(\lambda)$ or $h_s$, is the number of boxes $u$ such~that either~${u=s}$, or $u$ lies strictly below (respectively to the right) of $s$ in the same column (respectively row).

 The \textit{hook lengths multiset} of $\lambda$, denoted by $\mathcal{H}(\lambda)$, is the multiset of all hook lengths of $\lambda$. For any positive integer $n$, the multiset of all hook lengths that are congruent to $0 \pmod n$ is denoted by $\mathcal{H}_n(\lambda)$. Note that $\mathcal{H}(\lambda)=\mathcal{H}_1(\lambda)$. A partition $\omega$ is an \textit{$n$-core} if $n\not\in \cch(\omega)$ or equivalent, thanks to the Littlewood decomposition introduced in Section \ref{sec:combi}, if $\cch_n(\omega)=\varnothing$. For example, the only $2$-cores are the ``staircase'' partitions $(k,k-1,\dots,1)$, where $k$ is any positive integer.

Given a partition $\lambda $, an addable (resp.\ removable) node is a
node $b$ such that $\lambda \sqcup \{b\}$ (resp.~${\lambda \setminus \{b\}}$)
is again the diagram of a partition. The content of a node $b$ appearing in $%
\lambda $ at the intersection of its $j$-th column and its $i$-th row is
defined as $c(b)=j-i$. The residue of~$b$ is the value of $c(b)$ modulo $n$%
, that is, $r(b)=c(b)\bmod n$. The nodes with the same residue~$i$ are
called the $i$-nodes of $\lambda $. It is then easy to check that for any $%
i=0,\dots ,n-1$, an $n$-core cannot contain a mix of addable ($A$) and
removable ($R$) $i$-nodes: they are either $i$-nodes $A$, or $i$-nodes~$R$.

Let us denote by $\mathcal{C}_{n}$ the set of $n$-cores and write $%
\widetilde{\mathfrak{S}}_{n}$ for the affine Weyl group of type \smash{$A_{n-1}^{(1)}$}.
The group \smash{$\widetilde{\mathfrak{S}}_{n}$} is a Coxeter group
with simple generating reflections $s_{0},\dots ,s_{n-1}$. There is
a~classical action of $\widetilde{\mathfrak{S}}_{n}$ on the set $\mathcal{C}%
_{n}$: for any $i=0,\dots ,n-1$ and any $c\in \mathcal{C}_{n}$, the core $%
s_{i}\cdot c$ is obtained by removing all the addable $i$-nodes of $c$ when $%
c$ contains only removable $i$-nodes and adding all the possible addable $i$%
-nodes in $c$ when $c$ does not contain any addable $i$-node. Note that if $c$ contains neither removable $i$-nodes nor addable ones, then $%
s_{i}\cdot c=c$. Moreover, observe that $%
\mathcal{C}_{n}$ is stable by the transposition (or conjugation) $\mathrm{tr}
$ operation on partitions (exchanging the rows and the columns in the Young
diagrams). It is in fact easy to see that for any core~${c\in \mathcal{C}_{n}}$
such that $c=s_{i_{k}}\cdots s_{i_{1}}\cdot \varnothing $, we have $c^{\mathrm{%
tr}}=s_{n-i_{k}}\cdots s_{n-i_{1}}\cdot \varnothing $.

\begin{Example}
Assume $n=3$ and consider the $3$-core $c=(4,2,1,1)$. We give below the
action of the simple reflections $s_{0}$ and $s_{1}$ on $c$ (the number
indicated are the residues of the nodes):
\begin{equation*}
\begin{tabular}{lllll}
\cline{1-3}
\multicolumn{1}{|l}{$0$} & \multicolumn{1}{|l}{$1$} & \multicolumn{1}{|l}{$2$%
} & \multicolumn{1}{|l}{$0$} & $1$ \\ \cline{1-1}\cline{1-3}
\multicolumn{1}{|l}{$2$} & \multicolumn{1}{|l}{$0$} & $1$ & & \\
\cline{1-1}
\multicolumn{1}{|l}{$1$} & \multicolumn{1}{|l}{$2$} & & & \\ \cline{1-1}
$0$ & & & & \\
$2$ & & & &
\end{tabular}%
\overset{0}{\longleftarrow }%
\begin{tabular}{lllll}
\cline{1-4}
\multicolumn{1}{|l}{$0$} & \multicolumn{1}{|l}{$1$} & \multicolumn{1}{|l}{$2$%
} & \multicolumn{1}{|l}{$0$} & \multicolumn{1}{|l}{$1$} \\
\cline{1-2}\cline{1-4}
\multicolumn{1}{|l}{$2$} & \multicolumn{1}{|l}{$0$} & \multicolumn{1}{|l}{$1$%
} & & \\ \cline{1-1}\cline{1-2}
\multicolumn{1}{|l}{$1$} & \multicolumn{1}{|l}{$2$} & & & \\ \cline{1-1}
\multicolumn{1}{|l}{$0$} & \multicolumn{1}{|l}{} & & & \\ \cline{1-1}
$2$ & & & &
\end{tabular}%
\overset{1}{\rightarrow }%
\begin{tabular}{lllll}
\hline
\multicolumn{1}{|l}{$0$} & \multicolumn{1}{|l}{$1$} & \multicolumn{1}{|l}{$2$%
} & \multicolumn{1}{|l}{$0$} & \multicolumn{1}{|l|}{$1$} \\ \hline
\multicolumn{1}{|l}{$2$} & \multicolumn{1}{|l}{$0$} & \multicolumn{1}{|l}{$1$%
} & \multicolumn{1}{|l}{} & \\ \cline{1-1}\cline{1-3}
\multicolumn{1}{|l}{$1$} & \multicolumn{1}{|l}{$2$} & & & \\ \cline{1-1}
\multicolumn{1}{|l}{$0$} & \multicolumn{1}{|l}{} & & & \\ \cline{1-1}
$1$ & & & &
\end{tabular}.%
\end{equation*}
\end{Example}

A partition $\mu$ is \textit{self-conjugate} if $\mu=\mu^{\mathrm{tr}}$ or equivalently its Ferrers diagram is symmetric along the main diagonal. Let $\mathcal{SC}$ be the set of self-conjugate partitions and $\cC_{n}^s$ denote the set of self-conjugate $n$-cores.

\begin{Proposition}
We have $\mathcal{C}_{n}=\widetilde{\mathfrak{S}}_{n}\cdot \varnothing $, that is, the previous action is transitive on $\mathcal{C}_{n}$. Moreover, the
stabilizer of the empty node $\varnothing $ is the symmetric group $\mathfrak{S%
}_{n}=\langle s_{1},\dots ,s_{n-1}\rangle $ and the set $\mathcal{C}_{n}$
gives a parametrization of the affine Grassmannian elements of type \smash{$A_{n-1}^{(1)}$}.
\end{Proposition}

\begin{proof}
	It follows easily from the fact that any $n$-core partition has at least on removable $i$-node for an integer $i \in I$.
	Also it is clear that the stabilizer of the fundamental weight $\Lambda_0$ is the symmetric group $\mathfrak{S}_{n}$.
\end{proof}

For simplicity, we will slightly abuse the notation and identify each affine
Grassmannian element with its corresponding $n$-core. In this case, for
any $c\in W_{a}^{0}=\mathcal{C}_{n}$, we get
\begin{equation}
\Lambda _{0}-c(\Lambda _{0})=\sum_{i=0}^{n-1}a_{i}(c)\alpha _{i},
\label{DelaAffineA}
\end{equation}
where $a_{i}(c)$ is the number of $i$-nodes in the core $c$. Therefore, the
atomic length of $c$ satisfies
\begin{equation*}
L(c)=\sum_{i=0}^{n-1}a_{i}(c)=\vert c \vert,
\end{equation*}%
where $\vert c \vert $ is just the number of nodes in the core $c$. In the following paragraphs, we will see that it is possible to get a
similar combinatorial interpretation for any of the previous affine root
systems. For any integer $N\geq 2$, we will denote by $\alpha _{0},\dots,\alpha _{N-1}$ and $\Lambda _{0},\dots ,\Lambda _{N}$ (without
superscript) the simple roots and the fundamental weights of the affine root
system of type~\smash{$A_{N-1}^{(1)}$}. We will then express the simple roots
\smash{$\alpha _{i}^{X_{n}^{(a)}}$} of the affine root system \smash{$X_{n}^{(a)}$} in terms
of the~$\alpha _{j}$'s, and similarly the fundamental weights \smash{$\Lambda
_{i}^{X_{n}^{(a)}}$} in terms of the $\Lambda _{j}$'s. We will also decompose
the simple reflections generating each affine Weyl group \smash{$W_{X_{n}^{(a)}}$}
in terms of $s_{0},\dots ,s_{N-1}$, the simple reflections generating $%
\widetilde{\mathfrak{S}}_{N}$. The same convention will be used for the
atomic lengths: we will express \smash{$L_{X_{n}^{(a)}}$} in terms of $L$ the atomic
length of affine type $A$.

\subsection[Type C\_n\^{}(1)]{Type $\boldsymbol{C_{n}^{(1)}}$}

We can realize the affine root system of type~\smash{$C_{n}^{(1)}$} as the subsystem
of the root system of type~\smash{$A_{2n-1}^{(1)}$} such that
\begin{gather}
\alpha _{i}^{C_{n}^{(1)}}=\alpha _{i}+\alpha _{2n-i}, \quad i=1,\dots ,n-1 ,\qquad
\alpha _{n}^{C_{n}^{(1)}}=2\alpha _{n}, \qquad
\alpha _{0}^{C_{n}^{(1)}}=2\alpha _{0},\nonumber\\
\Lambda _{i}^{C_{n}^{(1)}}=\Lambda _{i}+\Lambda _{2n-i}, \quad i=1,\dots ,n-1, \qquad
\Lambda _{n}^{C_{n}^{(1)}}=2\Lambda _{n} ,\qquad
\Lambda _{0}^{C_{n}^{(1)}}=2\Lambda _{0}.
  \label{RelC}
\end{gather}
Indeed, one can then check that the previous relations are compatible with
the Dynkin diagram of type \smash{$C_{n}^{(1)}$}. The affine Weyl group \smash{$W_{^{C_{n}^{(1)}}}$} can then be seen as the subgroup of $\widetilde{
	\mathfrak{S}}_{2n}$ such that
\[
W_{^{C_{n}^{(1)}}}=\bigl\langle s_{i}^{C_{n}^{(1)}},\, i=0,\dots ,n\bigr\rangle
\]
with
\[
s_{i}^{C_{n}^{(1)}}=s_{i}s_{2n-i},\quad i=1,\dots ,n-1, \qquad
s_{n}^{C_{n}^{(1)}}=s_{n}, \qquad
s_{0}^{C_{n}^{(1)}}=s_{0}.
\]
The following Lemma is easy to prove. Set
\begin{equation*}
\mathcal{C}_{2n}^{s}=\bigl\{c\in \mathcal{C}_{2n}\mid c^{\mathrm{tr}}=c \bigr\}.
\end{equation*}

\begin{Lemma}
	\label{Lemma_selfconjugate}We have $W_{^{C_{n}^{(1)}}}\cdot \varnothing =
	\mathcal{C}_{2n}^{s}$, that is, the affine Grassmannian elements of type $
	C_{n}^{(1)}$ are parametrized by the self-conjugate $2n$-cores.
\end{Lemma}

\begin{proof}
	One first checks that \smash{$\mathcal{C}_{2n}^{s}$} is stable under the action of \smash{$W_{^{C_{n}^{(1)}}}$}.
    Next, it suffices to ob\-serve~that for any nonempty $%
	c\in \mathcal{C}_{2n}^{s}$, there exists at least an integer $i=0,\dots ,n$
	such that all~$i$~and $2n-i$ nodes of $c$ are removable because $c$
	is self-conjugate. This implies that \smash{$ s_{i}^{C_{n}^{(1)}}\cdot c=c^{\prime }$}
    where $c^{\prime }\neq c$ belongs to \smash{$\mathcal{C}_{2n}^{s}$} and has a
	number of nodes strictly less than $c$.
\end{proof}

We will identify the elements of the affine Grassmannian \smash{$W_{^{C_{n}^{(1)}}}^{0}$}
with the self-conjugate $2n$-cores in $\mathcal{C}%
_{2n}^{s}$. Recall also that the atomic length $L(c)$ in affine type $A$
then counts the number of boxes in the Young diagram of $c$. Then for each $c
$ in \smash{$\mathcal{C}_{2n}^{s}$}, we can write%
\begin{equation*}
\Lambda _{0}^{C_{n}^{(1)}}-c(\Lambda
_{0}^{C_{n}^{(1)}})=\sum_{i=0}^{n}a_{i}^{_{C_{n}^{(1)}}}(c)\alpha
_{i}^{C_{n}^{(1)}},
\end{equation*}%
which gives by using (\ref{RelC})%
\begin{align*}
	&2\Lambda _{0}-c(2\Lambda _{0})=2a_{0}^{C_{n}^{(1)}}(c)\alpha
	_{0}+2a_{n}^{C_{n}^{(1)}}(c)\alpha
	_{n}+\sum_{i=1}^{n-1}a_{i}^{C_{n}^{(1)}}(c)(\alpha _{i}+\alpha _{2n-i}), \qquad\text{i.e.,} \\
	&\Lambda _{0}-c(\Lambda _{0})=a_{0}^{C_{n}^{(1)}}(c)\alpha
	_{0}+a_{n}^{C_{n}^{(1)}}(c)\alpha _{n}+\sum_{i=1}^{n-1}\frac{1}{2}
	a_{i}^{C_{n}^{(1)}}(c)(\alpha _{i}+\alpha _{2n-i}).
\end{align*}
By comparing with (\ref{DelaAffineA}), we get
\begin{gather*}
a_{i}^{C_{n}^{(1)}}(c)=2a_{i}(c)=2a_{2n-i}(c)\quad \text{for }i=1,\dots ,n-1, \\
a_{n}^{C_{n}^{(1)}}(c)=a_{n}(c) \qquad \text{and} \qquad a_{0}^{C_{n}^{(1)}}(c)=a_{0}(c).
\end{gather*}
The following proposition is obtained by interpreting \smash{$L_{C_{n}^{(1)}}(c)$}
as the number of simple roots of type \smash{$C_{n}^{(1)}$} in the previous
decomposition of $\Lambda _{0}-c(\Lambda _{0})$.

\begin{Proposition}
	\label{Prop_C(n,1)}For any \smash{$c\in W_{^{C_{n}^{(1)}}}^{0}=\mathcal{C}_{2n}^{s}$}, we have
	\begin{equation*}
	L_{C_{n}^{(1)}}(c)=L(c), 
	\end{equation*}%
	that is, \smash{$L_{C_{n}^{(1)}}(c)$} is equal to the number of boxes in the self
	conjugate $2n$-core $c$.
\end{Proposition}

\subsection[Type D\_\{n+1\}\^{}(2)]{Type $\boldsymbol{D_{n+1}^{(2)}}$}

\label{sucsecD^2}This time, we can realize the affine root system of type $%
D_{n+1}^{(2)}$ as the subsystem of the root system of type \smash{$A_{2n-1}^{(1)}$}
such that
\begin{gather}
\alpha _{i}^{D_{n+1}^{(2)}}=\alpha _{i}+\alpha _{2n-i},\quad i=1,\dots ,n-1 ,\qquad
\alpha _{n}^{D_{n+1}^{(2)}}=\alpha _{n} ,\qquad
\alpha _{0}^{D_{n+1}^{(2)}}=\alpha _{0},\nonumber
\\
\Lambda _{i}^{D_{n+1}^{(2)}}=\Lambda _{i}+\Lambda _{2n-i},\quad i=1,\dots ,n-1 ,\qquad
\Lambda _{n}^{D_{n+1}^{(2)}}=\Lambda _{n}, \qquad
\Lambda _{0}^{D_{n+1}^{(2)}}=\Lambda _{0}.
 \label{RelCt}
\end{gather}
The affine Weyl group \smash{$W_{^{D_{n+1}^{(2)}}}$} is equal to \smash{$W_{^{C_{n}^{(1)}}}$}
and we yet have \smash{$W_{^{D_{n+1}^{(2)}}}\cdot \varnothing =\mathcal{C}_{2n}^{s}$}.
For each $c$ in~\smash{$\mathcal{C}_{2n}^{s}$}, we can write
\begin{equation*}
\Lambda _{0}^{D_{n+1}^{(2)}}-c\bigl(\Lambda
_{0}^{D_{n+1}^{(2)}}\bigr)=\sum_{i=0}^{n}a_{i}^{_{D_{n+1}^{(2)}}}(c)\alpha
_{i}^{D_{n+1}^{(2)}},
\end{equation*}%
which gives by using (\ref{RelCt})%
\begin{equation*}
\Lambda _{0}-c(\Lambda _{0})=a_{0}^{D_{n+1}^{(2)}}(c)\alpha
_{0}+a_{n}^{D_{n+1}^{(2)}}(c)\alpha
_{n}+\sum_{i=1}^{n-1}a_{i}^{D_{n+1}^{(2)}}(c)(\alpha _{i}+\alpha _{2n-i}).
\end{equation*}
By comparing with (\ref{DelaAffineA}), we get
\begin{gather}
a_{i}^{D_{n+1}^{(2)}}(c)=a_{i}(c)=a_{2n-i}(c)\quad \text{for }i=1,\dots ,n-1 ,\nonumber\\
a_{n}^{D_{n+1}^{(2)}}(c)=a_{n}(c)\qquad \text{and}\qquad
a_{0}^{D_{n+1}^{(2)}}(c)=a_{0}(c).
\label{A_i(Dn+1,2)}
\end{gather}

\begin{Proposition}
	\label{Prop_D(n+1,2)}For any \smash{$c\in W_{^{D_{n+1}^{(2)}}}^{0}=\mathcal{C}
	_{2n}^{s}$}, we have
	\begin{equation*}
	L_{D_{n+1}^{(2)}}(c)=\frac{1}{2}(L(c)+a_{0}(c)+a_{n}(c)).
	\end{equation*}
\end{Proposition}

\subsection[Type A\_\{2n\}\^{}(2)]{Type $\boldsymbol{A_{2n}^{(2)}}$}
\label{subsecA(2n,2)}

We can realize the affine root system of type \smash{$A_{2n}^{(2)}$} as the
subsystem of the root system of type~\smash{$A_{2n-1}^{(1)}$} such that%
\begin{gather*} 
\begin{split}
& \alpha _{i}^{A_{2n}^{(2)}}=\alpha _{i}+\alpha _{2n-i},\quad i=1,\dots ,n-1 ,\qquad
\alpha _{n}^{A_{2n}^{(2)}}=2\alpha _{n} ,\qquad
\alpha _{0}^{A_{2n}^{(2)}}=\alpha _{0},\\
& \Lambda _{i}^{A_{2n}^{(2)}}=\Lambda _{i}+\Lambda _{2n-i},\quad i=1,\dots ,n-1, \qquad
\Lambda _{n}^{A_{2n}^{(2)}}=2\Lambda _{n}, \qquad
\Lambda _{0}^{A_{2n}^{(2)}}=\Lambda _{0}.
\end{split}
\end{gather*}
The affine Weyl group $W_{^{A_{2n}^{(2)}}}$ is equal to $W_{^{C_{n}^{(1)}}}$
and we thus have $W_{^{A_{2n}^{(2)}}}\cdot \varnothing =\mathcal{C}_{2n}^{s}$. For each $c$ in~\smash{$\mathcal{C}_{2n}^{s}$}, we can write
\begin{equation*}
\Lambda _{0}^{A_{2n}^{(2)}}-c\bigl(\Lambda
_{0}^{A_{2n}^{(2)}}\bigr)=\sum_{i=0}^{n}a_{i}^{_{A_{2n}^{(2)}}}(c)\alpha
_{i}^{A_{2n}^{(2)}},
\end{equation*}%
which gives by using (\ref{RelCt})%
\begin{equation*}
\Lambda _{0}-c(\Lambda _{0})=a_{0}^{A_{2n}^{(2)}}(c)\alpha
_{0}+2a_{n}^{A_{2n}^{(2)}}(c)\alpha
_{n}+\sum_{i=1}^{n-1}a_{i}^{A_{2n}^{(2)}}(c)(\alpha _{i}+\alpha _{2n-i})%
 .
\end{equation*}%
By comparing with (\ref{DelaAffineA}), we get%
\begin{gather*}
a_{i}^{A_{2n}^{(2)}}(c)=a_{i}(c)=a_{2n-i}(c)\quad \text{for }i=1,\dots ,n-1, \\
a_{n}^{A_{2n}^{(2)}}(c)=\frac{1}{2}a_{n}(c)\qquad \text{and}\qquad
a_{0}^{A_{2n}^{(2)}}(c)=a_{0}(c).
\end{gather*}

\begin{Proposition}
	\label{Prop_A(2n,2)}For any \smash{$c\in W_{^{A_{2n}^{(2)}}}^{0}=$} $\mathcal{C}
	_{2n}^{s}$, we have
	\begin{equation*}
	L_{A_{2n}^{(2)}}(c)=\frac{1}{2}(L(c)+a_{0}(c)).
	\end{equation*}
\end{Proposition}

\subsection[Type A\_\{2n\}\^{}\{'(2)\}]{Type $\boldsymbol{A_{2n}^{\prime (2)}}$}

We denote by \smash{$A_{2n}^{\prime (2)}$} the affine root system obtained by
transposing the Cartan matrix of type~\smash{$A_{2n}^{ (2)}$}. The types \smash{$A_{2n}^{(2)}$}
and \smash{$A_{2n}^{\prime (2)}$} coincide up to relabelling of the
nodes of their Dynkin diagram. Nevertheless, this relabelling does not fix
the affine Grassmannian elements. Equivalently, we could also consider the
orbit of the fundamental weight \raisebox{-2pt}{{\smash{$\Lambda _{n}^{A_{2n}^{(2)}}$}}}. We realize
the affine root system of type \smash{$A_{2n}^{\prime (2)}$} as the subsystem of the
root system of type \smash{$A_{2n-1}^{(1)}$} such that%
\begin{gather*}
\alpha _{i}^{A_{2n}^{\prime (2)}}=\alpha _{i}+\alpha _{2n-i},\quad i=1,\dots ,n-1,
\qquad
\alpha _{n}^{A_{2n}^{\prime (2)}}=\alpha _{n}, \qquad
\alpha _{0}^{A_{2n}^{\prime (2)}}=2\alpha _{0},\\
\Lambda _{i}^{A_{2n}^{\prime (2)}}=\Lambda _{i}+\Lambda _{2n-i},\quad i=1,\dots,n-1, \qquad
\Lambda _{n}^{A_{2n}^{\prime (2)}}=\Lambda _{n}, \qquad
\Lambda _{0}^{A_{2n}^{\prime (2)}}=2\Lambda _{0}.
\end{gather*}
The affine Weyl group \smash{$W_{^{A_{2n}^{(2)}}}$} is equal to \smash{$W_{^{C_{n}^{(1)}}}$}
and we have \smash{$W_{^{A_{2n}^{(2)}}}\cdot \varnothing =\mathcal{C}_{2n}^{s}$}. For
each $c$ in \smash{$\mathcal{C}_{2n}^{s}$}, we can write
\begin{equation*}
\Lambda _{0}^{A_{2n}^{\prime (2)}}-c\bigl(\Lambda _{0}^{A_{2n}^{\prime
		(2)}}\bigr)=\sum_{i=0}^{n}a_{i}^{_{A_{2n}^{\prime (2)}}}(c)\alpha
_{i}^{A_{2n}^{\prime (2)}},
\end{equation*}
which gives
\begin{gather*}
	2\Lambda _{0}-2c(\Lambda _{0}) =2a_{0}^{A_{2n}^{\prime (2)}}(c)\alpha
	_{0}+a_{n}^{A_{2n}^{\prime (2)}}(c)\alpha
	_{n}+\sum_{i=1}^{n-1}a_{i}^{A_{2n}^{\prime (2)}}(c)(\alpha _{i}+\alpha
	_{2n-i}),\qquad \text{i.e.,} \\
	\Lambda _{0}-c(\Lambda _{0}) =a_{0}^{A_{2n}^{\prime (2)}}(c)\alpha _{0}+
	\frac{1}{2}a_{n}^{A_{2n}^{\prime (2)}}(c)\alpha _{n}+\sum_{i=1}^{n-1}\frac{1%
	}{2 }a_{i}^{A_{2n}^{\prime (2)}}(c)(\alpha _{i}+\alpha _{2n-i}).
\end{gather*}
By comparing with (\ref{DelaAffineA}), we get%
\begin{gather*}
a_{i}^{A_{2n}^{\prime (2)}}(c)=2a_{i}(c)=2a_{2n-i}(c)\quad \text{for }i=1,\dots
,n-1 ,\\
a_{n}^{A_{2n}^{\prime (2)}}(c)=2a_{n}(c)\qquad \text{and} \qquad a_{0}^{A_{2n}^{\prime
		(2)}}(c)=a_{0}(c).
\end{gather*}

\begin{Proposition}
	\label{Prop_A(2n,2)prime} For any \smash{$c\in W_{^{A_{2n}^{\prime (2)}}}^{0}=\mathcal{C}_{2n}^{s}$}, we have
	\begin{equation*}
	L_{A_{2n}^{\prime (2)}}(c)=L(c)+a_{n}(c).
	\end{equation*}
\end{Proposition}

\subsection[Type B\_n\^{}(1)]{Type $\boldsymbol{B_{n}^{(1)}}$}

We will proceed in two steps by first embedding the affine root system of
type \smash{$B_{n}^{(1)}$} into the root system of type \smash{$D_{n+1}^{(2)}$} and next by
using Section~\ref{sucsecD^2}. We first write%
\begin{gather}
\alpha _{i}^{B_{n}^{(1)}}=\alpha _{i}^{D_{n+1}^{(2)}},\quad i=1,\dots ,n, \qquad
\alpha _{0}^{B_{n}^{(1)}}=2\alpha _{0}^{D_{n+1}^{(2)}}+\alpha
_{1}^{D_{n+1}^{(2)}},\nonumber\\
\Lambda _{i}^{B_{n}^{(1)}}=\Lambda _{i}^{D_{n+1}^{(2)}},\quad i=2,\dots ,n, \qquad
\Lambda _{0}^{B_{n}^{(1)}}=\Lambda _{0}^{D_{n+1}^{(2)}}, \qquad
\Lambda _{1}^{B_{n}^{(1)}}=\Lambda _{1}^{D_{n+1}^{(2)}}-\Lambda
_{0}^{D_{n+1}^{(2)}},
  \label{RelBA}
\end{gather}%
which indeed gives a root system of type \smash{$B_{n}^{(1)}$}. The affine Weyl
group \smash{$W_{B_{n}^{(1)}}$} can then be realized as the subgroup of \smash{$W_{D_{n+1}^{(2)}}=W_{C_{n}^{(1)}}$} such that
\begin{equation*}
W_{B_{n}^{(1)}}=\bigl\langle
s_{0}^{B_{n}^{(1)}}=s_{0}^{C_{n}^{(1)}}s_{1}^{C_{n}^{(1)}}s_{0}^{C_{n}^{(1)}},s_{1}^{C_{n}^{(1)}},\dots ,s_{n}^{C_{n}^{(1)}}\bigr\rangle
\end{equation*}%
or equivalently, \smash{$W_{B_{n}^{(1)}}$} is the subgroup of $\widetilde{\mathfrak{S%
}}_{n}$%
\begin{equation*}
W_{B_{n}^{(1)}}=\langle s_{0}s_{1}s_{2n-1}s_{0},s_{1}s_{2n-1},\dots
,s_{n-1}s_{n+1},s_{n}\rangle .
\end{equation*}%
Write $\mathcal{C}_{2n}^{s,p}$ for the subset of \smash{$\mathcal{C}_{2n}^{s}$} of
self-conjugate $2n$-cores with an even number of nodes on its main diagonal
(i.e., an even number of nodes with content equal to $0$).

\begin{Lemma}
	We have $W_{^{B_{n}^{(1)}}}\cdot \varnothing =\mathcal{C}_{2n}^{s,p}$, that is,
	the affine Grassmannian elements of type \smash{$B_{n}^{(1)}$} are parametrized by
	the self-conjugate $2n$-cores with an even diagonal.
\end{Lemma}

\begin{proof}
 The proof is similar to that of Lemma \ref%
{Lemma_selfconjugate}. One first checks that $\mathcal{C}_{2n}^{s,p}$ is
stable under the action of \smash{$W_{^{B_{n}^{(1)}}}$}. Next, one observes that for
any nonempty $c\in \mathcal{C}_{2n}^{s,p}$, there exists at least an integer
$i=0,\dots ,n$ such that all the $i$ nodes of $c$ are removable. Since $c$
belongs to $\mathcal{C}_{2n}^{s,p}$, this implies that \raisebox{-1pt}{\smash{$s_{i}^{B_{n}^{(1)}}%
\cdot c=c^{\prime }$}}, where $c^{\prime }\neq c$ belongs to $\mathcal{C}
_{2n}^{s}$ and has a number of nodes strictly less than $c$.
\end{proof}

For each $c$ in $\mathcal{C}_{2n}^{s,p}$, we can write%
\begin{equation*}
\Lambda _{0}^{B_{n}^{(1)}}-c\bigl(\Lambda
_{0}^{B_{n}^{(1)}}\bigr)=\sum_{i=0}^{n}a_{i}^{B_{n}^{(1)}}(c)\alpha
_{i}^{B_{n}^{(1)}},
\end{equation*}%
which gives by using (\ref{RelBA})%
\begin{equation*}
\Lambda _{0}^{D_{n+1}^{(2)}}-c\bigl(\Lambda
_{0}^{D_{n+1}^{(2)}}\bigr)=a_{0}^{B_{n}^{(1)}}(c)\bigl(\alpha
_{1}^{D_{n+1}^{(2)}}+2\alpha
_{0}^{D_{n+1}^{(2)}}\bigr)+\sum_{i=1}^{n}a_{i}^{B_{n}^{(1)}}(c)\alpha
_{i}^{D_{n+1}^{(2)}} .
\end{equation*}%
We get by using (\ref{A_i(Dn+1,2)})
\begin{gather*}
a_{i}^{B_{n}^{(1)}}(c)=a_{i}^{D_{n+1}^{(2)}}(c)=a_{i}(c)=a_{2n-i}(c)\quad \text{for }i=2,\dots ,n ,\\
a_{0}^{B_{n}^{(1)}}(c)=\frac{1}{2}a_{0}^{D_{n+1}^{(2)}}(c)=\frac{1}{2}
a_{0}(c), \qquad   a_{1}^{B_{n}^{(1)}}(c)=a_{1}^{D_{n+1}^{(2)}}(c)-\frac{1}{
	2}a_{0}^{D_{n+1}^{(2)}}(c)=a_{1}(c)-\frac{1}{2}a_{0}(c).
\end{gather*}
The following proposition can then be deduced from Proposition \ref{Prop_D(n+1,2)}.

\begin{Proposition}
	\label{Prop_B(n,1)} For any \smash{$c\in W_{^{B_{n}^{(1)}}}^{0}=\mathcal{C}%
	_{2n}^{s,p}$}, we have
	\begin{equation*}
	L_{B_{n}^{(1)}}(c)=L_{D_{n+1}^{(2)}}(c)-a_{0}^{D_{n+1}^{(2)}}(c)=\frac{1}{2}
	(L(c)-a_{0}(c)+a_{n}(c)).
	\end{equation*}
\end{Proposition}

\subsection[Type A\_\{2n-1\}\^{}(2)]{Type $\boldsymbol{A_{2n-1}^{(2)}}$}

Here again, we will proceed in two steps by first embedding the affine root
system of type \smash{$A_{2n-1}^{(2)}$} into the root system of type \smash{$A_{2n}^{(2)}$}
and next by using Section~\ref{subsecA(2n,2)}. We first write%
\begin{gather}
\alpha _{i}^{A_{2n-1}^{(2)}}=\alpha _{i}^{A_{2n}^{(2)}},\quad i=1,\dots ,n, \qquad
\alpha _{0}^{A_{2n-1}^{(2)}}=2\alpha _{0}^{A_{2n}^{(2)}}+\alpha
_{1}^{A_{2n}^{(2)}},\nonumber\\
\Lambda _{i}^{A_{2n-1}^{(2)}}=\Lambda _{i}^{A_{2n}^{(2)}},\quad i=2,\dots ,n, \qquad
\Lambda _{0}^{A_{2n-1}^{(2)}}=\Lambda _{0}^{A_{2n}^{(2)}}, \qquad
\Lambda _{1}^{A_{2n-1}^{(2)}}=\Lambda _{1}^{A_{2n}^{(2)}}-\Lambda
_{0}^{A_{2n}^{(2)}},
 \label{RelAA}
\end{gather}
which indeed gives a root system of type \smash{$A_{2n-1}^{(2)}$}. The affine Weyl
group \smash{$W_{A_{2n-1}^{(2)}}$} can then be realized as the subgroup of \smash{$W_{A_{2n}^{(2)}}=W_{C_{n}^{(1)}}$} such that
\begin{equation*}
W_{A_{2n-1}^{(2)}}=\bigl\langle
s_{0}^{A_{2n-1}^{(2)}}=s_{0}^{C_{n}^{(1)}}s_{1}^{C_{n}^{(1)}}s_{0}^{C_{n}^{(1)}},s_{1}^{C_{n}^{(1)}},\dots ,s_{n}^{C_{n}^{(1)}}\bigr\rangle =W_{B_{n}^{(1)}}.
\end{equation*}
We thus get the following lemma.

\begin{Lemma}
	We have \smash{$W_{^{A_{2n-1}^{(2)}}}\cdot \varnothing =\mathcal{C}_{2n}^{s,p}$}.
\end{Lemma}

For each $c$ in \smash{$\mathcal{C}_{2n}^{s,p}$}, we can write%
\begin{equation*}
\Lambda _{0}^{A_{2n-1}^{(2)}}-c\bigl(\Lambda
_{0}^{A_{2n-1}^{(2)}}\bigr)=\sum_{i=0}^{n}a_{i}^{A_{2n-1}^{(2)}}(c)\alpha
_{i}^{A_{2n-1}^{(2)}},
\end{equation*}%
which gives by using (\ref{RelAA})%
\begin{equation*}
\Lambda _{0}^{A_{2n}^{(2)}}-c\bigl(\Lambda
_{0}^{A_{2n}^{(2)}}\bigr)=a_{0}^{A_{2n-1}^{(2)}}(c)\bigl(\alpha
_{1}^{A_{2n}^{(2)}}+2\alpha
_{0}^{A_{2n}^{(2)}}\bigr)+\sum_{i=1}^{n}a_{i}^{A_{2n-1}^{(2)}}(c)\alpha
_{i}^{A_{2n}^{(2)}} .
\end{equation*}%
We get by using (\ref{A_i(Dn+1,2)})
\begin{gather*}
a_{i}^{A_{2n-1}^{(2)}}(c)=a_{i}^{A_{2n}^{(2)}}(c)=a_{i}(c)=a_{2n-i}(c)\quad \text{for }i=2,\dots ,n-1 ,\\
a_{n}^{A_{2n-1}^{(2)}}(c)=a_{n}^{A_{2n}^{(2)}}(c)=\frac{1}{2}a_{n}(c) ,\qquad
a_{0}^{A_{2n-1}^{(2)}}(c)=\frac{1}{2}a_{0}^{A_{2n}^{(2)}}(c)=\frac{1}{2}
a_{0}(c),\\
a_{1}^{A_{2n-1}^{(2)}}(c)=a_{1}^{A_{2n}^{(2)}}(c)-\frac{1}{2}a_{0}^{A_{2n}^{(2)}}(c)=a_{1}(c)-\frac{1}{2}a_{0}(c).
\end{gather*}
The following proposition can then be deduced from Proposition \ref{Prop_D(n+1,2)}.

\begin{Proposition}
	\label{Prop_A(2n-1,2)} For any \smash{$c\in W_{^{A_{2n-1}^{(2)}}}^{0}=\mathcal{C}%
	_{2n}^{s,p}$}, we have
	\begin{equation*}
	L_{A_{2n-1}^{(2)}}(c)=L_{A_{2n}^{(2)}}(c)-a_{0}^{A_{2n}^{(2)}}(c)=\frac{1}{2}
	(L(c)-a_{0}(c)).
	\end{equation*}
\end{Proposition}

\subsection[Type D\_n\^{}(1)]{Type $\boldsymbol{D_{n}^{(1)}}$}

Here again we use an embedding in the affine root system of type \smash{$D_{n+1}^{(2)}$} and write
\begin{gather}
\alpha _{i}^{D_{n}^{(1)}}=\alpha _{i}^{D_{n+1}^{(2)}},\quad i=1,\dots ,n-1, \qquad
\alpha _{0}^{D_{n}^{(1)}}=2\alpha _{0}^{D_{n+1}^{(2)}}+\alpha
_{1}^{D_{n+1}^{(2)}} ,\nonumber\\
\alpha _{n}^{D_{n}^{(1)}}=2\alpha _{n}^{D_{n+1}^{(2)}}+\alpha
_{n-1}^{D_{n+1}^{(2)}},\nonumber\\
\Lambda _{i}^{D_{n}^{(1)}}=\Lambda _{i}^{D_{n+1}^{(2)}},\quad i=2,\dots ,n-1, \qquad
\Lambda _{0}^{D_{n}^{(1)}}=\Lambda _{0}^{D_{n+1}^{(2)}},\Lambda
_{n}^{D_{n}^{(1)}}=\Lambda _{n}^{D_{n+1}^{(2)}} ,\nonumber\\
\Lambda _{1}^{D_{n}^{(1)}}=\Lambda _{1}^{D_{n+1}^{(2)}}-\Lambda
_{0}^{D_{n+1}^{(2)}},\qquad \Lambda _{n-1}^{D_{n}^{(1)}}=\Lambda
_{n-1}^{D_{n+1}^{(2)}}-\Lambda _{n}^{D_{n+1}^{(2)}},
 \label{RelDA}
\end{gather}
which gives a root system of type \smash{$D_{n}^{(1)}$}. The affine Weyl group \smash{$W_{D_{n}^{(1)}}$} can then be realized as the subgroup of
\smash{$W_{D_{n+1}^{(1)}}=W_{C_{n}^{(1)}}$} such that
\begin{equation*}
W_{D_{n}^{(1)}}=\bigl\langle
s_{0}^{D_{n}^{(1)}}=s_{0}^{C_{n}^{(1)}}s_{1}^{C_{n}^{(1)}}s_{0}^{C_{n}^{(1)}},s_{1}^{C_{n}^{(1)}},\dots ,s_{n-1}^{C_{n}^{(1)}},s_{n}^{D_{n}^{(1)}}=s_{n}^{C_{n}^{(1)}}s_{n-1}^{C_{n}^{(1)}}s_{n}^{C_{n}^{(1)}}\bigr\rangle
\end{equation*}%
or equivalently, \smash{$W_{D_{n}^{(1)}}$} is the subgroup of $\widetilde{\mathfrak{S%
}}_{n}$%
\begin{equation*}
W_{D_{n}^{(1)}}=\langle s_{0}s_{1}s_{2n-1}s_{0},s_{1}s_{2n-1},\dots
,s_{n-1}s_{n+1},s_{n}s_{n-1}s_{n+1}s_{n}\rangle .
\end{equation*}%
According to Table \ref{table2}, the affine Grassmannian
\[
W_{D_{n}^{(1)}}^{0}=W_{D_{n}^{(1)}}/W_{D_{n}}
\]
 is in one-to-one
correspondence with the sublattice lattice of $\mathbb{Z}^{n}$ of vectors $%
\beta =(\beta _{1},\dots ,\beta _{n})$ such that $\beta _{1}+\cdots +\beta
_{n}$ is even. On the other hand, we must have
\[
W_{D_{n}^{(1)}}^{0}\subset
W_{B_{n}^{(1)}}^{0}\qquad \text{for $W_{D_{n}^{(1)}}\subset W_{B_{n}^{(1)}}$}.
\]
 But by
using Table \ref{table2} again, the affine Grassmannian
\[
W_{B_{n}^{(1)}}^{0}=W_{B_{n}^{(1)}}/W_{B_{n}}
\]
 is again in one-to-one
correspondence with the vectors in $\mathbb{Z}^{n}$ whose sum of coordinates
is even. Therefore, we have
\[
W_{D_{n}^{(1)}}^{0}=W_{B_{n}^{(1)}}^{0}
\]
 and \raisebox{1pt}{\smash{$W_{D_{n}^{(1)}}^{0}$}}
 is yet parametrized by the elements of \smash{$\mathcal{C}_{2n}^{s,p}$}.

 Observe nevertheless that the length functions in \smash{$W_{D_{n}^{(1)}}^{0}$} and \smash{$W_{B_{n}^{(1)}}^{0}$} do not coincide.

For each $c$ in $\mathcal{C}_{2n}^{s,p}$, we can write%
\begin{equation*}
\Lambda _{0}^{D_{n}^{(1)}}-c\bigl(\Lambda
_{0}^{D_{n}^{(1)}}\bigr)=\sum_{i=0}^{n}a_{i}^{D_{n}^{(1)}}(c)\alpha
_{i}^{D_{n}^{(1)}},
\end{equation*}%
which gives by using (\ref{RelDA})%
\begin{gather*}
\Lambda _{0}^{D_{n+1}^{(2)}}-c\bigl(\Lambda_{0}^{D_{n+1}^{(2)}}\bigr)= a_{0}^{D_{n}^{(1)}}(c)\bigl(\alpha_{1}^{D_{n+1}^{(2)}}+2\alpha_{0}^{D_{n+1}^{(2)}}\bigr) \\
\hphantom{\Lambda _{0}^{D_{n+1}^{(2)}}-c(\Lambda_{0}^{D_{n+1}^{(2)}})=}{}
+\sum_{i=1}^{n-1}a_{i}^{D_{n}^{(1)}}(c)\alpha_{i}^{D_{n+1}^{(2)}}+a_{n}^{D_{n}^{(1)}}(c)\bigl(\alpha_{n-1}^{D_{n+1}^{(2)}}+2\alpha _{n}^{D_{n+1}^{(2)}}\bigr) .
\end{gather*}
We get by using (\ref{A_i(Dn+1,2)})
\begin{gather*}
a_{i}^{D_{n}^{(1)}}(c)=a_{i}^{D_{n+1}^{(2)}}(c)=a_{i}(c)=a_{2n-i}(c)\quad \text{for }i=2,\dots ,n-1 ,\\
a_{0}^{D_{n}^{(1)}}(c)=\frac{1}{2}a_{0}^{D_{n+1}^{(2)}}(c)=\frac{1}{2}
a_{0}(c), \qquad a_{1}^{D_{n}^{(1)}}(c)=a_{1}^{D_{n+1}^{(2)}}(c)-\frac{1}{
	2}a_{0}^{D_{n+1}^{(2)}}(c)=a_{1}(c)-\frac{1}{2}a_{0}(c) ,\\
a_{n}^{D_{n}^{(1)}}(c)=\frac{1}{2}a_{n}^{D_{n+1}^{(2)}}(c)=\frac{1}{2}
a_{n}(c),\qquad  a_{n-1}^{D_{n}^{(1)}}(c)=a_{n-1}^{D_{n+1}^{(2)}}(c)-
\frac{1}{2}a_{n}^{D_{n+1}^{(2)}}(c)=a_{n-1}(c)-\frac{1}{2}a_{n}(c).	
\end{gather*}
The following proposition can then be deduced from Proposition \ref{Prop_A(2n-1,2)}.

\begin{Proposition}
	\label{Prop_D(n,1)} For any \smash{$c\in W_{^{D_{n}^{(1)}}}^{0}=\mathcal{C}_{2n}^{s,p}$}, we have
	\begin{equation*}
	L_{D_{n}^{(1)}}(c)=L_{D_{n+1}^{(2)}}(c)-a_{0}^{D_{n+1}^{(2)}}(c)-a_{n}^{D_{n+1}^{(2)}}(c)=
	\frac{1}{2}(L(c)-a_{0}(c)-a_{n}(c)).
	\end{equation*}
\end{Proposition}

\subsection{Weighted boxes in cores}\label{subsec:weighted}

The previous formulas for the atomic lengths can be expressed succinctly by
defining a weight function $\pi $ on the relevant cores depending only on
the residues of the nodes. We then have
\begin{equation*}
L_{X_{n}^{(a)}}(c)=\sum_{\square \in c}\pi (\square ),
\end{equation*}%
where the map $\pi $ takes the same value $\pi _{i}$ on the set of nodes
having residue $i$. This is equivalent to associate a weight $\pi _{i}$ on
each node of the Dynkin diagram of the root system considered. When $\pi
_{i}=1$ for any $i$, the previous sum is just equal to the number of nodes
in $c$.
\begin{center}\renewcommand{\arraystretch}{1.32}
\begin{tabular}{|c|c|c|c|c|}
\hline
type & $A_{n-1}^{(1)}$ & $C_{n}^{(1)}$ & $B_{n}^{(1)}$ & $D_{n}^{(1)}$ \\
\hline
core set & $\mathcal{C}_{n}$ & $\mathcal{C}_{2n}^{s}$ & $\mathcal{C}%
_{2n}^{s,p}$ & $\mathcal{C}_{2n}^{s,p}$ \\ \hline
$\pi $ & $%
\begin{array}{c}
\pi _{i}=1 \\
\text{any }i%
\end{array}%
$ & $%
\begin{array}{c}
\pi _{i}=1 \\
\text{any }i%
\end{array}%
$ & $%
\begin{array}{c}
\pi _{0}=0,\pi _{n}=1 \\
\pi _{i}=\frac{1}{2},\text{o.t.w}%
\end{array}%
$ & $%
\begin{array}{c}
\pi _{0}=\pi _{n}=0 \\
\pi _{i}=\frac{1}{2},\text{o.t.w}%
\end{array}%
$ \\ \hline
\end{tabular}
\vspace{2mm}

\begin{tabular}{|c|c|c|c|c|}
\hline
type & $D_{n}^{(2)}$ & $A_{2n}^{(2)}$ & $A_{2n}^{\prime (2)}$ & $%
A_{2n-1}^{(2)}$ \\ \hline
core set & $\mathcal{C}_{2n}^{s}$ & $\mathcal{C}_{2n}^{s}$ & $\mathcal{C}%
_{2n}^{s}$ & $\mathcal{C}_{2n}^{s,p}$ \\ \hline
$\pi $ & $%
\begin{array}{c}
\pi _{0}=\pi _{n}=1 \\
\pi _{i}=\frac{1}{2},\text{o.t.w}%
\end{array}%
$ & $%
\begin{array}{c}
\pi _{0}=1 \\
\pi _{i}=\frac{1}{2},\text{o.t.w}%
\end{array}%
$ & $%
\begin{array}{c}
\pi _{0}=2 \\
\pi _{i}=1,\text{o.t.w}%
\end{array}%
$ & $%
\begin{array}{c}
\pi _{0}=0 \\
\pi _{i}=\frac{1}{2},\text{o.t.w}%
\end{array}%
$ \\ \hline
\end{tabular}
\end{center}

\subsection[Type D\_4\^{}(3)]{Type $\boldsymbol{D_{4}^{(3)}}$}

We realize the affine root system of type \smash{$D_{4}^{(3)}$} from the root system
of type \smash{$A_{5}^{(2)}$} by setting%
\begin{gather}
\alpha _{0}^{D_{4}^{(3)}}=\alpha _{0}^{A_{5}^{(2)}}, \qquad
\alpha _{1}^{D_{4}^{(3)}}=\alpha _{2}^{A_{5}^{(2)}}, \qquad
\alpha _{2}^{D_{4}^{(3)}}=\alpha _{1}^{A_{5}^{(2)}}+\alpha _{3}^{A_{5}^{(2)}},\nonumber\\
\Lambda _{0}^{D_{4}^{(3)}}=\Lambda _{0}^{A_{5}^{(2)}}, \qquad
\Lambda _{1}^{D_{4}^{(3)}}=\Lambda _{2}^{A_{5}^{(2)}}, \qquad
\Lambda _{2}^{D_{4}^{(3)}}=\Lambda _{1}^{A_{5}^{(2)}}+\Lambda
_{3}^{A_{5}^{(2)}}.\label{relaDG}
\end{gather}
The affine Weyl group \smash{$W_{D_{4}^{(3)}}$} can then be realized as the subgroup
of \smash{$W_{A_{5}^{(2)}}$} such that
\begin{equation*}
W_{D_{4}^{(3)}}=\bigl\langle
s_{0}^{D_{4}^{(3)}}=s_{0}^{A_{5}^{(2)}},s_{1}^{D_{4}^{(3)}}=s_{2}^{A_{5}^{(2)}},s_{2}^{D_{4}^{(3)}}=s_{1}^{A_{5}^{(2)}}s_{3}^{A_{5}^{(2)}}\bigr\rangle =\langle s_{0}s_{1}s_{5}s_{0},s_{2}s_{4},s_{1}s_{3}s_{5}\rangle .
\end{equation*}%
Let $\mathcal{C}_{6}^{g}$ be the set of $6$-cores associated to the elements
in $\mathbb{Z}^{6}$ of the form%
\begin{equation}\label{eq:d43}
(\beta _{1},\beta _{2},\beta _{1}-\beta _{2},\beta _{2}-\beta _{1},-\beta
_{2},-\beta _{1}).
\end{equation}%
Observe in particular that $\mathcal{C}_{6}^{g}\subset \mathcal{C}_{6}^{s,p}$%
.

\begin{Lemma}
	We have \smash{$W_{D_{4}^{(3)}}\cdot \varnothing =\mathcal{C}_{6}^{g}$}.
\end{Lemma}

\begin{proof}
	One checks easily that $\mathcal{C}_{6}^{g}$ is stable under the action of \smash{$W_{D_{4}^{(3)}}$}.

Moreover, when $c$ is not empty in $\mathcal{C}_{6}^{g}$,
	at least one of the actions of \smash{$s_{0}^{D_{4}^{(3)}}$} or \smash{$s_{1}^{D_{4}^{(3)}}$}
	on $c$ makes decrease its number of boxes.
\end{proof}

\begin{Remark}
Note that \smash{$s_3^{A_6^{(1)}}$} acts on $\mathcal{C}_{6}^{g}$ such that \smash{$s_3^{A_6^{(1)}}\mathcal{C}_{6}^{g}$} correspond to the set of self-conjugate $6$-cores such that the $2$-quotient is a $3$-core partition.
\end{Remark}

For each $c$ in $\mathcal{C}_{6}^{g}$, we can write%
\begin{equation*}
\Lambda _{0}^{D_{4}^{(3)}}-c\bigl(\Lambda
_{0}^{D_{4}^{(3)}}\bigr)=\sum_{i=0}^{2}a_{i}^{D_{4}^{(3)}}(c)\alpha
_{i}^{D_{4}^{(3)}},
\end{equation*}%
which gives by using (\ref{relaDG})%
\begin{equation*}
\Lambda _{0}^{A_{5}^{(2)}}-c\bigl(\Lambda
_{0}^{A_{5}^{(2)}}\bigr)=a_{0}^{D_{4}^{(3)}}(c)\alpha
_{0}^{A_{5}^{(2)}}+a_{1}^{D_{4}^{(3)}}(c)\alpha
_{2}^{A_{5}^{(2)}}+a_{2}^{D_{4}^{(3)}}(c)\bigl(\alpha _{1}^{A_{5}^{(2)}}+\alpha
_{3}^{A_{5}^{(2)}}\bigr).
\end{equation*}%
We get%
\begin{gather*}
a_{0}^{D_{4}^{(3)}}(c)=a_{0}^{A_{5}^{(2)}}(c)=\frac{1}{2}a_{0}(c), \qquad
a_{1}^{D_{4}^{(3)}}(c)=a_{2}^{A_{5}^{(2)}}(c)=a_{2}(c), \\
a_{2}^{D_{4}^{(3)}}(c)=a_{1}^{A_{5}^{(2)}}(c)=a_{3}^{A_{5}^{(2)}}(c)=a_{1}(c)-
\frac{1}{2}a_{0}(c)=\frac{1}{2}a_{3}(c).
\end{gather*}
The following proposition can then be deduced from Proposition \ref{Prop_A(2n-1,2)}.

\begin{Proposition}
\label{Prop_D(4,3)}	For any \smash{$c\in W_{^{D_{4}^{(3)}}}^{0}=\mathcal{C}_{6}^{g}$}, we have
	\begin{equation*}
	L_{D_{4}^{(3)}}(c)=L_{A_{5}^{(2)}}(c)-a_{3}^{A_{5}^{(2)}}(c)=\frac{1}{2}
	(L(c)-a_{0}(c)-a_{3}(c)).
	\end{equation*}
\end{Proposition}

\subsection[Type G\_2\^{}(1)]{Type $\boldsymbol{G_{2}^{(1)}}$}

We realize the affine root system of type \smash{$G_{2}^{(1)}$} from the root system
of type \smash{$A_{5}^{(2)}$} by setting%
\begin{gather}
\alpha _{0}^{G_{2}^{(1)}}=3\alpha _{0}^{A_{5}^{(2)}}, \qquad
\alpha _{1}^{G_{2}^{(1)}}=3\alpha _{2}^{A_{5}^{(2)}}, \qquad
\alpha _{2}^{G_{2}^{(1)}}=\alpha _{1}^{A_{5}^{(2)}}+\alpha _{3}^{A_{5}^{(2)}},\nonumber\\
\Lambda _{0}^{G_{2}^{(1)}}=3\Lambda _{0}^{A_{5}^{(2)}}, \qquad
\Lambda _{1}^{G_{2}^{(1)}}=3\Lambda _{2}^{A_{5}^{(2)}}, \qquad
\Lambda _{2}^{G_{2}^{(1)}}=\Lambda _{1}^{A_{5}^{(2)}}+\Lambda
_{3}^{A_{5}^{(2)}}.\label{relaG}
\end{gather}
The affine Weyl group \smash{$W_{G_{2}^{(1)}}$} is the same as \smash{$W_{D_{4}^{(3)}}$}
\begin{equation*}
W_{G_{2}^{(1)}}=\bigl\langle
s_{0}^{G_{2}^{(1)}}=s_{0}^{A_{5}^{(2)}},s_{1}^{G_{2}^{(1)}}=s_{2}^{A_{5}^{(2)}},s_{2}^{G_{2}^{(1)}}=s_{1}^{A_{5}^{(2)}}s_{3}^{A_{5}^{(2)}}\bigr\rangle =\langle s_{0}s_{1}s_{5}s_{0},s_{2}s_{4},s_{1}s_{3}s_{5}\rangle .
\end{equation*}
and we again have \smash{$W_{G_{2}^{(1)}}\cdot \varnothing =\mathcal{C}_{6}^{g}$}.

For each $c$ in $\mathcal{C}_{6}^{g}$, we can write%
\begin{equation*}
\Lambda _{0}^{G_{2}^{(1)}}-c\bigl(\Lambda
_{0}^{G_{2}^{(1)}}\bigr)=\sum_{i=0}^{2}a_{i}^{G_{2}^{(1)}}(c)\alpha
_{i}^{G_{2}^{(1)}},
\end{equation*}%
which gives by using (\ref{relaG})%
\begin{align*}
3\Lambda _{0}^{A_{5}^{(2)}}-3c\bigl(\Lambda
_{0}^{A_{5}^{(2)}}\bigr)&{}=3a_{0}^{G_{2}^{(1)}}(c)\alpha
_{0}^{A_{5}^{(2)}}+3a_{1}^{G_{2}^{(1)}}(c)\alpha
_{2}^{A_{5}^{(2)}}+a_{2}^{G_{2}^{(1)}}(c)\bigl(\alpha _{1}^{A_{5}^{(2)}}+\alpha
_{3}^{A_{5}^{(2)}}\bigr)\\
&{}= 3a_{0}^{A_{5}^{(2)}}(c)\alpha
_{0}^{A_{5}^{(2)}}+3a_{1}^{A_{5}^{(2)}}(c)\alpha
_{1}^{A_{5}^{(2)}}+3a_{2}^{A_{5}^{(2)}}(c)\alpha
_{2}^{A_{5}^{(2)}}+3a_{3}^{A_{5}^{(2)}}(c)\alpha _{3}^{A_{5}^{(2)}}.
\end{align*}
We get
\begin{gather*}
a_{0}^{G_{2}^{(1)}}(c)=a_{0}^{A_{5}^{(2)}}(c)=\frac{1}{2}a_{0}(c), \qquad
a_{1}^{G_{2}^{(1)}}(c)=a_{2}^{A_{5}^{(2)}}(c)=a_{2}(c), \\
a_{2}^{G_{2}^{(1)}}(c)=3a_{1}^{A_{5}^{(2)}}(c)=3a_{3}^{A_{5}^{(2)}}(c)=3a_{1}(c)-%
\frac{3}{2}a_{0}(c)=\frac{3}{2}a_{3}(c).
\end{gather*}
The following proposition can then be deduced from Proposition \ref
{Prop_A(2n-1,2)}.

\begin{Proposition}\label{Prop_G(2,1)}
	For any $c\in W_{^{G_{2}^{(1)}}}^{0}=\mathcal{C}_{6}^{g}$, we have
	\begin{equation*}
	L_{G_{2}^{(1)}}(c)=L_{A_{5}^{(2)}}(c)+a_{3}^{A_{5}^{(2)}}(c)=\frac{1}{2}
	(L(c)-a_{0}(c)+a_{3}(c)).
	\end{equation*}
\end{Proposition}

\section{Combinatorics of integer partitions}\label{sec:combi}

We have seen in the previous section how the affine Grassmannian elements are in one-to-one correspondence with some core partitions and how their atomic length could be expressed as modified weights of these partitions. As mentioned previously, the Nekrasov--Okounkov formula is derived by using a polynomial argument. Unfortunately, the formulas for the atomic lengths obtained in Section~\ref{subsec:weighted} due to these modified weights heavily depend on the rank of the affine root system considered when $\pi_i$ is not constant for any $i$: the weight of a fixed box changes with the rank $n$ considered. Therefore, this interpretation does not permit the direct use of the polynomial argument for all affine root systems. In this section, we bypass this obstruction by making explicit the one-to-one correspondence between the core partitions used in Section~\ref{sec:atomic} and integer partitions whose weight (number of boxes) exactly corresponds to the value of the atomic length.

We first recall a couple of results about the Littlewood decomposition, a bijection mapping an integer partition to its $n$-core and an $n$-tuple of partitions. Then we study the restriction of this application to the set of distinct partitions.
In what follows, we will map bijectively the elements $c$ of the subsets of cores arising in the previous section to elements $c'$ of subsets of partitions defined implicitly in terms of their abaci such that the atomic length of $c$ is equal to the weight of $c'$, i.e., $L(c)=\lvert c'\rvert$. This combinatorial construction relies on particular extremal representations of each affine Lie algebra that can be regarded as an analogue of the natural representations of the classical Lie algebras considered as matrix algebras.

\subsection{The Littlewood decomposition and bi-infinite binary words}

Recall that there is a natural correspondence between $\ccp$ and the set of bi-infinite words indexed by $\Z$ over the alphabet $\lbrace 0,1\rbrace$.
\begin{Definition}\label{def:psi}
Define
\[
\psi\colon \   \ccp   \to   \lbrace 0,1\rbrace^{\Z}, \qquad
   \lambda  \mapsto  (c_k)_{k\in\Z},
\]
such that
 \[c_k=\begin{cases}
 0 \quad \text{if } k\in\lbrace \lambda_i-i,\, i\in\N\rbrace,\\
 1 \quad \text{if } k\in \lbrace j-\lambda_j^{\mathrm{tr}}-1,\, j\in\N\rbrace.
 \end{cases}\]
\end{Definition}
Moreover, by definition of $\psi(\lambda)$, one has:
\begin{equation}\label{eq:balance}
\#\{k\leq-1,\, c_k=1\}= \#\{k\geq0,\, c_k=0\}=d_\lambda.
\end{equation}
 The application $\psi$ is a bijection from the set of integer partitions $\ccp$ and the subset of bi-infinite binary words
 \[\bigl\lbrace (c_k)\in\lbrace 0,1\rbrace^{\Z}\mid \#\{k\leq-1,\,c_k=1\}= \#\{k\geq0,\,c_k=0\}=d_\lambda\bigr\rbrace.\]
 Let $\partial \lambda$ be the border of the Ferrers diagram of $\lambda$. Each step on $\partial\lambda$ is either horizontal or vertical. The above correspondence amounts to encoding the walk along the border from the South-West to the North-East as depicted in Figure \ref{fig:word}: take ``$0$'' for a vertical step and ``$1$'' for a horizontal step. The resulting word is indexed by $\Z$. In order to keep this correspondence bijective, one needs to fix the index $0$. The choice within this framework is to set the letter of index $0$ to be the first step after the corner of the Durfee square, the largest square that can fit within the Ferrers diagram of $\lambda$. The sequence $\psi(\lambda)$ has many names across the literature, with slight variations on the alphabet, or the binary labels, such as Maya diagrams, edge sequences, Dirac sea, abacus. Note that the edge sequence defined in~\cite{LLMSmem} corresponds to the word obtained by labelling horizontal and vertical steps with letters ``$0$'' and letters ``$1$'', respectively.

\begin{figure}[th!]
\centering
\begin{tikzpicture}
  [
    dot/.style={circle,draw=black, fill,inner sep=1pt},
  ]

\foreach \x in {0,...,2}{
  \node[dot] at (\x,-4){ };
}

\foreach \x in {.2,1.2}
  \draw[->,thick,orange] (\x,-4) -- (\x+.6,-4);
\foreach \x in {3.2,4.2}
  \draw[->,thick,orange] (\x,-2) -- (\x+.6,-2);
\foreach \y in {1.8,.8}
  \draw[->,thick,orange] (5,-\y) -- (5,-\y+.6);
\draw[->,thick,orange] (2,-3.8) -- (2,-3.8+.6);
\draw[->,thick,orange] (3,-2.8) -- (3,-2.8+.6);
\draw[->,thick,orange] (2.2,-3) -- (2.2+.6,-3);

\node[dot] at (2,-3){};
\foreach \x in {3,...,5}
  \node[dot] at (\x,-2){};
\foreach \x in {1,...,4}
  \draw (\x,-.1) -- node[above,xshift=-0.4cm,yshift=1mm] {$\lambda_\x'$} (\x,+.1);

\node[above,xshift=0.5cm,yshift=1mm] at (4,0) {$\lambda_5'$};
\node[above,xshift=0.5cm,yshift=1mm] at (5,0) {NE};
\node[above,xshift=-4mm,yshift=1mm] at (0,0) {NW};

\foreach \y in {1,...,3}
  \draw (.1,-\y) -- node[above,xshift=-4mm,yshift=0.2cm] {$\lambda_\y$} (-.1,-\y);
\node[above,xshift=-4mm,yshift=2mm] at (0,-4) {$\lambda_4$};
\node[above,xshift=-4mm] at (0,-5) {SW};
\node at (0,-4.5) {0};
\node at (0,-5.5) {0};
\node at (2,-3.5) {0};
\node at (3,-2.5) {0};
\node at (5,-1.5) {0};
\node at (5,-0.5) {0};

\node at (0.5,-4) {1};
\node at (1.5,-4) {1};
\node at (4.5,-2) {1};
\node at (3.5,-2) {1};
\node at (5.5,-0) {1};
\node at (2.5,-3) {1};

\draw[dotted,gray] (0,-1)--(5,-1);
\draw[dotted,gray] (0,-2)--(3,-2);
\draw[dotted,gray] (0,-3)--(2,-3);

\draw[dotted,gray] (1,0)--(1,-4);
\draw[dotted,gray] (2,0)--(2,-3);
\draw[dotted,gray] (3,0)--(3,-2);
\draw[dotted,gray] (4,0)--(4,-2);

\draw[->,very thick, orange](1.5,-3.5)--(1.5,-0.5);
\draw[->,very thick, orange](1.5,-0.5)--(4.5,-0.5);

\node[dot] at (5,-1){};
\node[dot] at (5,0){};
\draw[->,thick,-latex] (0,-6) -- (0,-5);
\draw[thick] (0,-6) -- (0,1);
\draw[->,thick,-latex] (-1,0) -- (6,0);
\node[circle,draw=blue,fill=blue,inner sep=0pt,minimum size=5pt] at (3,-3){};

\end{tikzpicture}
\caption{$\partial\lambda$ and its binary correspondence for $\lambda=(5,5,3,2)$ with a hook.}
\label{fig:word}
\end{figure}

As illustrated in Figure \ref{fig:word}, the boxes of $\lambda$ are in bijective correspondence with letters of $\psi(\lambda)$.
\begin{Lemma}[{\cite[Lemma $2.1$]{Wmac}}]\label{lem:indices}
The map $\psi$ $($see Definition $\ref{def:psi})$ associates bijectively a box $s$ of hook length $h_s$ of the Ferrers diagram of $\lambda$ to a pair of indices $(i_s,j_s)\in\Z^2$ of the word $\psi(\lambda)$ such that
\begin{enumerate}\itemsep=0pt
\item[$(1)$] $i_s<j_s$,
\item[$(2)$] $c_{i_s}=1$, $c_{j_s}=0$,
\item[$(3)$] $j_s-i_s=h_s$,
\item[$(4)$] $s$ is a box above the main diagonal in the Ferrers diagram of $\lambda$ if and only if the number of letters {\rm ``$1$''} with negative index greater than $i_s$ is lower than the number of letters {\rm ``$0$''} with nonnegative index lower than $j_s$.
\end{enumerate}
\end{Lemma}

Hook lengths formulas are useful enumerative tools bridging combinatorics with other fields such as representation theory, probability, gauge theory or algebraic geometry. A much more recent identity is the Nekrasov--Okounkov formula. It was discovered independently by Nekrasov and Okounkov in their work on random partitions and Seiberg--Witten theory~\cite{NO}, by Westbury~\cite{We} in his work on universal characters for $\mathfrak{sl}_n$, and later by Han~\cite{Ha} based on one of the identities for affine type \smash{$A^{(1)}_{n-1}$}~in~\cite[Appendix $1$]{Mac} and a polynomial argument. This formula is commonly stated as follows:
\begin{equation*}
\sum_{\lambda\in\ccp}q^{\lvert \lambda\rvert}\prod_{h\in\cch(\lambda)}\biggl(1-\frac{z}{h^2}\biggr)=\prod_{k\geq 1}\bigl(1-q^k\bigr)^{z-1},
\end{equation*}
where $z$ is a fixed complex number.

Han's proof is based on the crucial observation that if we take $z=n^2$ in \eqref{NOdebut}, the products \smash{$\prod_{h\in\cch(\lambda)}\bigl(1-\frac{n^2}{h^2}\bigr)$} cancel whenever $\lambda$ does not belong to $\mathcal{C}_n$. Han's proof is to show that in this case, the equality \eqref{NOdebut} corresponds to a specialization of the Weyl--Kac denominator formula~in type \smash{$A_{n-1}^{(1)}$}. Since the equality holds for an infinite number of $n$, Han's proof is based on a~polynomial argument.

\begin{figure}[th]\centering
\begin{tabular}{ccc}
\begin{ytableau}
\none[\textcolor{red!60}{\rightarrow}]& & & & & &\none \\
\none[\textcolor{red!60}{\rightarrow}]&\none[\textcolor{red!60}{\rightarrow}]&&
\\
\none[\textcolor{red!60}{\rightarrow}]&\none[\textcolor{red!60}{\rightarrow}]&\none[\textcolor{red!60}{\rightarrow}]&\\
 \none &\none &\none &\none \\
 \none &\none &\none &\none \\
\end{ytableau}
&\qquad &
\begin{ytableau}
 *(green!60)+ &+ & + &+ & +\\
 - &*(green!60)+ &+ & \none \\
 - &- &*(green!60)+&\none \\
 -&\none &\none[\textcolor{green}{D(\lambda)}] &\none \\
 - &\none &\none &\none \\
\end{ytableau}\\
(a) && (b)
\vspace{4mm}\\
\begin{ytableau}
 *(green!60)+ &+ & + &*(yellow!60)+&+ & +\\
 - &*(green!60)+ &+ &*(yellow!60)+ & \none \\
 - &- &*(green!60)+&*(yellow!60)+ &\none \\
 -&\none &\none[\textcolor{green}{D(\lambda)}] &\none & \none \\
 - &\none &\none &\none & \none \\
\end{ytableau}
& \qquad &
\begin{ytableau}
 *(green!60)+ &+ & + &+&+ \\
 - &*(green!60)+ &+ & \none \\
 - &- &*(green!60)+&\none &\none[\textcolor{green}{D(\lambda)}] \\
 *(yellow!60)-&*(yellow!60)-&*(yellow!60)-&\none\\
 -&\none &\none &\none & \none \\
 - &\none &\none &\none & \none \\
\end{ytableau}\\
(c) && (d)
\end{tabular}
\caption{Distinct partition, a self-conjugate partition, a doubled distinct partition and its conjugate filled with $\varepsilon$:
(a)~shifted Young diagram of $\bar{\lambda}=(5,2,1)\in\ccd$, (b)~$\lambda=(5,3,3,1,1)\in\ccsc$, (c)~$\lambda=(6,4,4,1,1)\in\ccdd$, (d)~$\lambda=(5,3,3,3,1,1)\in\ccdd^{\mathrm{tr}}$.}
\label{fig:var}

\end{figure}

 The set of distinct partitions, denoted by $\ccd$, is the set of partitions such that no consecutive parts are equal (see~\cite{HX,Macbook,Sta99}). A distinct partition $\bar{\lambda}$ is identified with its \textit{shifted Young diagram}, which means the $i$-th row of the usual Young diagram is shifted by $i$ boxes to the right. The \textit{doubled distinct partition} of $\bar{\lambda}\in\ccd$, denoted by $\bar{\lambda}\bar{\lambda}$, is defined to be the usual partition whose Young diagram is obtained by adding $\bar{\lambda}_i$ boxes to the $i$-th column of the shifted Young diagram of $\bar{\lambda}$ for $1\leq i\leq \ell (\bar{\lambda})$.
 We denote the set of doubled distinct partitions by $\ccdd$.

  \textbf{Warning.} In Figure \ref{fig:var}\,(a), the shifted Young diagram of $\bar{\lambda}$ should have every row shifted by~$1$ to the right in order to obtain $\bar{\lambda}\bar{\lambda}$ from Figure~\ref{fig:var}\,(c).

 Following~\cite{HX}, the leftmost box of the $i$-th row of the shifted Young diagram of $\bar{\lambda}$ has coordinate $(i,i+1)$. The \textit{hook length} of a box of coordinate $(i,j)$ in the shifted diagram is the number of boxes strictly to the right, strictly below, the box itself plus $\bar{\lambda}_j$. Let us denote by~$\cch\bigl(\bar{\lambda}\bigr)$ the multiset of hook lengths in the shifted diagram for the strict partition $\bar{\lambda}$. One has then the following relation on multisets:
 \begin{equation}\label{eq:hookdd}
  \cch\bigl(\bar{\lambda}\bar{\lambda}\bigr)=\cch\bigl(\bar{\lambda}\bigr)\uplus\cch\bigl(\bar{\lambda}\bigr)\uplus\bigl\lbrace 2\bar{\lambda}_1,\dots, 2\bar{\lambda}_\ell\bigr\rbrace\setminus\bigl\lbrace \bar{\lambda}_1,\dots ,\bar{\lambda}_\ell\bigr\rbrace,
 \end{equation}
 where we use the symbol $\uplus$ for the union of multisets. Let us give an alternative definition of the set $\ccdd$, is that of all partitions $\lambda$ of Durfee square size $d$ such that ${\lambda_i=\lambda_i^{\mathrm{tr}}+1}$ for all~${i\in\lbrace 1,\dots,d\rbrace}$. We also define the set of conjugate of doubled distinct partitions $\ccdd^{\mathrm{tr}}:=\bigl\lbrace \lambda^{\mathrm{tr}}\mid \lambda\in\ccdd\bigr\rbrace$.
 These constructions maps $\bar{\lambda}=(1)\in\ccd$ to $\bar{\lambda}\bar{\lambda}=(2)\in \ccdd$ and to~${(1,1)\in\ccdd^{\mathrm{tr}}}$.

  Recall that a partition $\mu$ is \textit{self-conjugate} if its Ferrers diagram is symmetric along the main diagonal. Such a partition $\mu$ can also be constructed from an element $\bar{\lambda}\in\ccd$: $\mu$ is defined to be the usual partition whose Young diagram is obtained by adding $\bar{\lambda}_i-1$ boxes to the $i$-th column of the shifted Young diagram of $\bar{\lambda}$ for $1\leq i\leq \ell\bigl(\bar{\lambda}\bigr)$. One can go from a self-conjugate partition to a doubled distinct partition, respectively conjugate doubled distinct partition, by adding a~vertical strip, respectively a horizontal strip, of length of the size of the Durfee square (shaded in yellow in Figure \ref{fig:var}\,(c), respectively in Figure \ref{fig:var}\,(d)).

   For instance, in Figure \ref{fig:var}, take $\bar{\lambda}=(5,2,1)\in \ccd$, the corresponding element in the set of self-conjugate partitions $\ccsc$ $\lambda=(5,3,3,1,1)$ in Figure \ref{fig:var}\,(a) has its main diagonal $D(\lambda)$ shaded in green while in Figure \ref{fig:var}\,(c) $\bar{\lambda}\bar{\lambda}=(6,4,4,1,1)\in\ccdd$ has its main diagonal shaded in green as for the strip shaded in yellow, it corresponds to the boxes added to a self-conjugate partition to obtain a doubled distinct partition. The conjugate of a doubled distinct partition is also illustrated in Figure \ref{fig:var}\,(d). Note that if we take $\mu$ in one of the sets $\ccsc$, $\ccdd$, $\ccdd^{\mathrm{tr}}$ and $\bar{\lambda}\in\ccd$ the corresponding distinct partition, one has the following relation:
\begin{equation*}
 d_\mu=\ell\bigl(\bar{\lambda}\bigr).
 \end{equation*}
Note that the hook lengths on the main diagonal $D(\lambda)$ are $ 2 \bar{\lambda}_1,\dots,2\bar{\lambda}_\ell$. Moreover, the multiset of hook lengths of the $(l+1)$-th column (respectively row) of the Young diagram of the doubled distinct partition (respectively conjugate doubled distinct partition) coloured in yellow in~Figure~\ref{fig:var}\,(c) (respectively Figure~\ref{fig:var}\,(d)) is $\bigl\lbrace \bar{\lambda}_1,\dots,\bar{\lambda}_\ell\bigr\rbrace$.

\begin{figure}[h]
\centering
\begin{tabular}{ccc}
\begin{ytableau}
 *(green!60)9 &*(green!60)6 & *(green!60) 5&*(green!60)2 & *(green!60)1\\
 6 &*(green!60)3 &*(green!60)2 & \none \\
 5 &2 &*(green!60)1&\none \\
 2&\none &\none[\textcolor{green}{\bar{\lambda}}] &\none \\
 1 &\none &\none &\none \\
\end{ytableau}
&
\begin{ytableau}
 10 &*(green!60)7 & *(green!60)6 &*(green!60)5&*(green!60)2 & *(green!60)1\\
7 &4 &*(green!60)3 &*(green!60)2 & \none \\
 6 &3 & 2&*(green!60)1 &\none \\
 2&\none &\none[\textcolor{green}{\bar{\lambda}}] &\none & \none \\
 1 &\none &\none &\none & \none \\
\end{ytableau}
&
\begin{ytableau}
 *(green!60)10 &*(green!60)7 & *(green!60)6 &*(green!60)2&*(green!60)1 \\
 7 &*(green!60)4 &*(green!60)3 & \none \\
 6 &3 &*(green!60)2&\none[\textcolor{green}{\bar{\lambda}}] \\
5&2&1&\none\\
 2&\none &\none &\none & \none \\
 1 &\none &\none &\none & \none \\
\end{ytableau}\\
(a) $\lambda=(5,3,3,1,1)\in\ccsc$ & (b) $\lambda=(6,4,4,1,1)\in\ccdd$ & (c) $\lambda=(5,3,3,3,1,1)\in\ccdd^{\mathrm{tr}}$
\end{tabular}
\caption{A self-conjugate partition, a doubled distinct partition and its conjugate filled with its hook lengths.}

\end{figure}

Let us introduce here a signed statistic \textit{$\varepsilon_s$}, for a box $s$ of $\bar{\lambda}\bar{\lambda}$, is defined as $-$ if $s$ is strictly below the main diagonal of the Ferrers diagram of $\lambda$ and as $+$ otherwise, as depicted in Figure~\ref{fig:var}. This signed statistic already appears algebraically within the work of King~\cite{King} and combinatorially within the work of P\'etr\'eolle~\cite{MP}.
Moreover, the multiset $\cch\bigl(\bar{\lambda}\bigr)$ of hook lengths of~$\bar{\lambda}$ is equal to the multiset of hook lengths of $\bar{\lambda}\bar{\lambda}$ strictly above the main diagonal~$D(\lambda)$, which is~${\bigl\lbrace h_s,\, s\in\bar{\lambda}\bar{\lambda}\setminus D\bigl(\bar{\lambda}\bar{\lambda}\bigr),\,\varepsilon=+\bigr\rbrace=\cch\bigl(\bar{\lambda}\bigr)}$.

\begin{Remark}\label{rem:conju}
Let $\lambda$ be a partition and $\psi(\lambda)=(c_k)_{k\in\Z}$ be its corresponding word, as introduced in Definition \ref{def:psi}. Let $\lambda^{\mathrm{tr}}$ be the conjugate of $\lambda$ and $\psi\bigl(\lambda^{\mathrm{tr}}\bigr)=\bigl(c^{\mathrm{tr}}_i\bigr)_{i\in\Z}$. We have
\[\forall k\in\Z,\qquad
 c^{\mathrm{tr}}_k=1-c_{-k-1}.\]
\end{Remark}

Given the properties of symmetries of self-conjugate and doubled distinct partitions, they can alternatively be characterized by
\begin{equation}\label{eq:motdd}
\lambda\in\ccdd \iff \psi(\lambda)=(c_k)_{k\in\Z}, \qquad c_0=1  \quad \text{and} \quad  \forall k \in \N^*,\ c_{-k}=1-c_k,
\end{equation}
and
\begin{equation*}
\lambda\in\ccsc \iff \psi(\lambda)=(c_k)_{k\in\Z}, \qquad \forall k \in \N,\ c_{-k-1}=1-c_k.
\end{equation*}
From Remark \ref{rem:conju} and \eqref{eq:motdd}, the set $\ccdd^{\mathrm{tr}}$ also admits a similar characterization:
\begin{equation}\label{eq:motddprime}
\lambda\in\ccdd ^{\mathrm{tr}}\iff \psi(\lambda)=(c_k)_{k\in\Z},\qquad  c_{-1}=0  \quad \text{and} \quad  \forall k \in \N,\ c_{-k-2}=1-c_k.
\end{equation}

From the above relations and the definition of $\psi$ (see Definition~\ref{def:psi}), one has the following lemma.
\begin{Lemma}\label{lem:distinct}
Set $\bar{\lambda}=\bigl(\bar{\lambda_1},\dots,\bar{\lambda}_\ell\bigr)\in\ccd$ and let \smash{$\psi\bigl(\bar{\lambda}\bar{\lambda}\bigr)=(c_k)_{k\in\Z}$} and $\psi\bigl(\bar{\lambda}\bar{\lambda}^{\mathrm{tr}}\bigr)=\bigl(c^{\mathrm{tr}}_k\bigr)_{k\in\Z}$ and $\mu$ the self-conjugate partition corresponding to $\bar{\lambda}$. Set $\psi(\mu)=(c^{s}_k)_{k\in\Z}$. Then
\begin{align*}
&\lbrace k\in\N\mid c_k=0\rbrace=\bigl\lbrace \bar{\lambda}_i,\, i\in\lbrace 1,\dots,\ell\rbrace\bigr\rbrace,\\
&\bigl\lbrace k\in\N\mid c^{\mathrm{tr}}_k=0\bigr\rbrace=\big\lbrace \bar{\lambda}_i-1,\, i\in\lbrace 1,\dots,\ell\rbrace\big\rbrace,\\
&\lbrace k\in\N\mid c^{s}_k=0\rbrace=\big\lbrace \bar{\lambda}_i-1,\, i\in\lbrace 1,\dots,\ell\rbrace\big\rbrace.
\end{align*}
\end{Lemma}

Take a partition $\lambda$ and a strictly positive integer $n$. Obtaining what is called the $n$-quotient of $\lambda$ is straightforward from $\psi(\lambda)=(c_i)_{i\in\Z}$: we just look at subwords with indices congruent to the same values modulo $n$. In this case, the equality \eqref{eq:balance} is not necessarily verified. To be able to apply $\psi^{-1}$, the index has to be shifted by $\#\lbrace i\in\N \mid c_{ni+k}=1\rbrace-\#\lbrace i\in\N^* \mid c_{-ni+k}=1\rbrace$.

The sequence $10$ within these subwords are replaced iteratively by $01$ until the subwords are all the infinite sequence of ``$0$'''s before the infinite sequence of ``$1$'''s (in fact it consists in removing all rim hooks in $\lambda$ of length congruent to $0\pmod n$). Let $\omega$ be the partition corresponding to the word which has the subwords $\pmod n$ obtained after the removal of the~$10$~sequences. Note that $\omega$ is an $n$-core.

Now we recall the following classical map, often called the Littlewood decomposition (see for instance~\cite{GKS,HJ}).

\begin{Definition}\label{def:Phi}{\em
 Let $n \geq 2$ be an integer and consider
\[
\Phi_n\colon \  \mathcal{P}   \to   \mathcal{C}_{n} \times \mathcal{P}^n, \qquad
  \lambda   \mapsto   \bigl(\omega,\nu^{(0)},\dots,\nu^{(n-1)}\bigr),
\]
where, if we set $\psi(\lambda)=(c_i)_{i\in\Z}$, then for all $k\in\lbrace 0,\dots,n-1\rbrace$, one has
\[
\nu^{(k)}:=\psi^{-1}\bigl(\bigl(c_{n(i+m_k)+k}\bigr)_{i\in\Z}\bigr),
\]
 where $m_k=\#\lbrace i\in\N \mid c_{ni+k}=0\rbrace-\#\lbrace i\in\N^* \mid c_{-ni+k}=1\rbrace$.
 The tuple $\underline{\nu}=\bigl(\nu^{(0)},\dots,\nu^{(n-1)}\bigr)$ is called the $n$-quotient of $\lambda$ and is denoted by \textit{$\quot_n(\lambda)$}, while $\omega$ is the $n$-core of $\lambda$ denoted by~\textit{$\core_n(\lambda)$}.}
\end{Definition}
\begin{Proposition}
Let $n \geq 2$ be an integer. The application $\Phi_n$ is a bijection between $\ccp$ and~${\mathcal{C}_{n} \times \mathcal{P}^n}$.
\end{Proposition}
For example, if we take $\lambda = (4,4,3,2) \text{ and } n=3$, then $\psi(\lambda)=\dots \color{red}{0} \color{blue}{0} \color{green}{1} \color{red}{1} \color{blue}{0}\color{green}{1} \color{black}| \color{red}{0} \color{blue}{1} \color{green}{0} \color{red}{0} \color{blue}{1} \color{green}{1}\color{black}\dots$
\begin{align*}
\begin{array}{@{}rc|rcl}
\psi\left(\nu^{(0)}\right)=\dots \color{red} 001 \color{black}| \color{red}001\color{black}\dots& &\psi\left(w_{0}\right)=\dots \color{red} 000 \color{black}| \color{red}011\color{black}\dots,\\
 \psi\left(\nu^{(1)}\right)=\dots \color{blue} 000 \color{black}| \color{blue}111\color{black}\dots& \longmapsto& \psi\left(w_{1}\right)=\dots \color{blue} 000 \color{black}| \color{blue}111\color{black}\dots , \\
 \psi\left(\nu^{(2)}\right)=\dots \color{green} 011 \color{black}| \color{green}011\color{black}\dots & & \psi\left(w_{2}\right)=\dots \color{green} 001 \color{black}| \color{green}111\color{black}\dots .
\end{array}
\end{align*}
Thus
\[
\psi(\omega)=\dots \color{red}{0} \color{blue}{0} \color{green}{0} \color{red}{0} \color{blue}{0}\color{green}{1} \color{black}| \color{red}{0} \color{blue}{1} \color{green}{1} \color{red}{1} \color{blue}{1} \color{green}{1}\color{black}\dots
\]
and
\[
\quot_3(\lambda)=\bigl(\nu^{(0)},\nu^{(1)},\nu^{(2)}\bigr)=((1,1),\varnothing,(2)),\qquad \core_3(\lambda)= \omega=(1).
\]

Now we discuss the Littlewood decomposition for $\ccdd^{\mathrm{tr}}$. Let $n$ be a positive integer, take $\lambda\in \ccdd^{\mathrm{tr}}$, and set $\psi(\lambda)=(c_i)_{i \in\mathbb{Z}}\in \lbrace 0,1\rbrace^\Z$, as introduced in Definition \ref{def:psi}, and $(\omega,\underline{\nu})=(\core_n(\lambda),\quot_n(\lambda))$. Using \eqref{eq:motdd}, one has the equivalence (see, for instance,~\cite{Wmac}):
\begin{align*}
\lambda\in \ccdd^{\mathrm{tr}} \ &\iff \ \forall i_0 \in \lbrace 0,\dots,n-1\rbrace,\ \forall j \in \mathbb{N},\ c_{i_0+jn}=1-c_{-i_0-jn-2},\\
&\iff \  \forall i_0 \in \lbrace 0,\dots,n-1\rbrace,\ \forall j \in \mathbb{N},\ c_{i_0+jn}=1-c_{n-(i_0+2)-n(j-1)},\\
&\iff \ \forall i_0 \in \lbrace 0,\dots,n-1 \rbrace,\ \nu^{(i_0)}=\bigl(\nu^{(n-i_0-2)}\bigr)^{\mathrm{tr}}\quad \text{and}\quad \omega \in \ccdd^{\mathrm{tr}}_{(n)} .
\end{align*}

Therefore, $\lambda$ is uniquely defined if its $n$-core is known as well as the $\left\lfloor n/2\right\rfloor$ first elements of its $n$-quotient, which are partitions without any constraint. It implies that if $n$ is odd, there is a~one-to-one correspondence between a $\ccdd^{\mathrm{tr}}$ and a triplet made of one~$\ccdd^{\mathrm{tr}}$ $n$-core, an element of~$\ccdd^{\mathrm{tr}}$ and $(n-3)/2$ generic partitions. If $n$ is even, the Littlewood decomposition is a one to one correspondence between a self-conjugate partition and a quadruplet made of one~$\ccdd^{\mathrm{tr}}$ $n$-core, $(n-2)/2$ generic partitions, a partition of~$\ccdd^{\mathrm{tr}}$ and a self-conjugate partition $\mu=\nu^{((n-1)/2)}$.
Hence the restriction of the Littlewood decomposition when applied to elements of $\ccdd^{\mathrm{tr}}$ is as follows.
\begin{Lemma}[{\cite[Lemma 2.8]{Wmac}}]\label{lem:DDprimeLittlewood}
Let $n$ be a positive integer. The Littlewood decomposition $\Phi_n$ $($see Definition $\ref{def:Phi})$ maps a conjugate doubled distinct partition $\lambda$ to $\bigl(\omega,\nu^{(0)},\dots,\nu^{(n-1)}\bigr)=(\omega,\underline{\nu})$ such that:
\begin{align*}
&(DD'1)\quad \text{the first component } \omega \text{ is a $\ccdd^{\mathrm{tr}}$ $n$-core and }\nu^{(0)},\dots,\nu^{(n-1)} \text{are partitions},\\
&(DD'2) \quad \forall j \in \lbrace 0,\dots,\lfloor n/2 \rfloor-1\rbrace,\ \nu^{(j)}=\bigl(\nu^{(n-2-j)}\bigr)^{\mathrm{tr}},\ \nu^{(n-1)}\in\ccdd^{\mathrm{tr}},
\\
&\quad \quad \quad \quad \text{and if $n$ is even, } \nu^{(n/2-1)}=\bigl(\nu^{(n/2-1)}\bigr)^{\mathrm{tr}}\in\ccsc, \\
&(DD'3) \quad |\lambda|=\begin{cases}\displaystyle|\omega|+2n\sum_{i=0}^{(n-3)/2} \bigl\lvert\nu^{(i)}\bigr\rvert +n\bigl\lvert \nu^{(n-1)}\bigr\rvert&
 \text{if n is odd},\\
\displaystyle|\omega|+2n\sum_{i=0}^{n/2-2} \bigl\lvert\nu^{(i)}\bigr\rvert+n\bigl\lvert \nu^{(n-1)}\bigr\rvert+n\bigl\lvert \nu^{(n/2-1)}\bigr\rvert & \text{if n is even},
\end{cases}\\
&(DD'4)\quad \mathcal{H}_n(\lambda)=n\mathcal{H}(\underline{\nu}),
\intertext{where a set $S$,}
&\quad nS:=\{ns,s\in S\}\notag
\intertext{and}
&\quad \mathcal{H}(\underline{\nu}):=\bigcup\limits_{i=0}^{n/2-1} \bigl(\cch \bigl(\nu^{(i)}\bigr)
\uplus\cch\bigl(\nu^{(n-2-i)}\bigr)\bigr).
\end{align*}
\end{Lemma}

\subsection{Cores and multicharges}

Let us first start by defining the multicharge associated to a core.

\begin{Definition}
	\label{defgks} Let $n\geq 2$ be an integer and consider
	\begin{equation*}
	\phi _{n}\colon \ \mathcal{C}_{n}   \rightarrow   \Z^{n}, \qquad
	  \omega   \mapsto   (m_{0},\dots ,m_{n-1}),
	\end{equation*}
	with $m_{i}:=\min \{k\in \Z\mid c_{kn+i}=1\}=\min \bigl\{k-\lambda _{k}^{\mathrm{tr}}-1,\, k\in \N\mid k-\lambda _{k}^{\mathrm{tr}}-1\equiv i\pmod{n}\bigr\}$.
\end{Definition}

We then get the following theorem.

\begin{Theorem}[{\cite[Theorem $2.10$]{Johnson}} and {\cite[Bijection $2$]{GKS}}]\label{thm:gk}
    Let $\omega $ be an $n$-core and $\psi (\omega )=( c_{k}) _{k\in \Z}$
	be its corresponding word $($see Definition~$\ref{def:psi})$. The map $\phi _{n}$
	is a bijection from $\cC_{n}$ to {$\Z_{0}^{n}:=\bigl\{(m_{i})_{0\leq i\leq n-1}\in %
	\Z^{n}\mid \sum_{i=0}^{n-1}m_{i}=0\bigr\}$}. Moreover, we have
	\begin{equation*}
	\lvert \omega \rvert =\frac{n}{2}\sum_{i=0}^{n-1}m_{i}^{2}+%
	\sum_{i=0}^{n-1}im_{i}. 
	\end{equation*}
\end{Theorem}

Given an $n$-core $\omega $, recall that $a_{i}$ counts the number of nodes
(or boxes) in the Young diagram of $c$ with residues $i$. One can also
define the vector $(\beta _{1},\dots ,\beta _{n})\in \mathbb{Z}^{n}$.
\begin{equation}
\beta _{i}:=a_{(i-1)\func{mod}n}-a_{i\func{mod}n} \label{eq:beta_a}
\end{equation}
for any $i=1,\dots ,n-1$.
It can be proven (see, for instance,~\cite{Rost}) that the two vectors $%
(\beta _{1},\dots ,\beta _{n})$ and $(a_{0},\dots ,a_{n-1})$ are related
according to the following formulas:
\begin{equation*}
a_{0}=\frac{1}{2} \Vert \beta  \Vert _{2}^{2}=\frac{1}{2}\bigl(\beta
_{1}^{2}+\cdots +\beta _{n}^{2}\bigr)
\end{equation*}
and for any $i=1,\dots ,n-1$
\begin{equation*}
a_{i}=a_{0}-\beta _{1}-\cdots -\beta _{i}.
\end{equation*}
In fact, both vectors $(\beta _{1},\dots ,\beta _{n})$ and $
(m_{0},\dots ,m_{n-1})$ are related by%
\begin{equation*}
m_{i}=\beta _{n-i}\qquad \text{for any }i=0,\dots ,n-1.
\end{equation*}

\begin{Remark}
	Observe also that when $\omega $ is regarded as an element of the affine
	Grassmannian of type \smash{$A_{n-1}^{(1)}$}, the vector $\beta $ is the one appearing
	in~(\ref{betac}).
\end{Remark}

In the case of a self-conjugate $n$-core, we get the following corollary of
Theorem \ref{thm:gk} (see \mbox{\cite[Sections~$7$ and~$8$]{GKS}} and~\cite[Section $%
4.1$]{Wmac} for details).

\begin{Corollary}\label{cor:gksc} Let $\omega \in \cC_{n}$ be an $n$-core. Then $\omega \in \cC_{n}^{s}$ $($i.e., is a self-conjugate $n$-core$)$ if and only if
    $\phi_{n}(\omega )=(m_{0},\dots ,m_{n-1})\in \Z^{n}$ $($see Definition $\ref{defgks})$ satisfies $m_{n-1-i}=-m_{i}$ for any~${i=0,\dots ,\lfloor n/2\rfloor}$. In
	particular, if $n$ is odd, we have $m_{\lfloor n/2\rfloor }=-m_{\lfloor
		n/2\rfloor }=0$. Moreover, the map from $\cC_{n}^{s}$ to $\Z^{\lfloor
		n/2\rfloor }$ that sends $\omega $ to $\bigl(m_{0},\dots ,m_{\lfloor n/2\rfloor
		-1}\bigr)$ is a bijection and we have the equality
	\begin{equation*}
	\lvert \omega \rvert =n\sum_{i=0}^{\lfloor n/2\rfloor
		-1}m_{i}^{2}+\sum_{i=0}^{\lfloor n/2\rfloor
		-1}(2i-n+1)m_{i}=n\sum_{i=1}^{\lfloor n/2\rfloor }\beta
	_{i}^{2}-\sum_{i=1}^{\lfloor n/2\rfloor }(2i-1)\beta _{i}.
	\end{equation*}
\end{Corollary}

Observe that in the right-hand side of the above equality seen as a polynomial in the variables~$(m_i)_{0\leq i \leq \lfloor (n-1)/2\rfloor }$,
 the coefficients of the $m_{i}$'s have a parity opposite to the leading coefficient.
Moreover, none of the coefficients of the monomial is divisible by the
leading coefficient. A~particular important case for the paper is that of
the core $c$ in \smash{$\mathcal{C}_{2n}^{s}$}. We then get ${\phi
_{2n}(c)=(m_{0},m_{1},\dots ,m_{2n-1})}$ with $m_{2n-1-i}=-m_{i}$ for any $%
i=0,\dots ,n-1$ and by Proposition~\ref{Prop_C(n,1)}
\begin{equation}
\lvert c\rvert
=2n\sum_{i=0}^{n-1}m_{i}^{2}+%
\sum_{i=0}^{n-1}(2i-2n+1)m_{i}=L_{A_{2n-1}^{(1)}}(c)=L_{C_{n}^{(1)}}(c).
\label{eq:Lc}
\end{equation}

We get a similar corollary for the set of ``doubled distinct $n$ cores'': $\cC%
_{n}^{dd}:=\cC_{n}\cap \ccdd$.

\begin{Corollary}
	\label{cor:gkdd}Let $\omega \in \cC_{n}$ be an $n$-core.Then $\omega $
	belongs to $\cC_{n}^{dd}$ if and only if $\phi _{n}(\omega )=(m_{0},\dots
	,m_{n-1})\in \Z^{n}$ satisfies $m_{0}=0$ and $m_{n-i}=-m_{i}$ for all $i\in
	\{1,\dots ,\lfloor (n-1)/2\rfloor \}$. Moreover, the map from~$\cC_{n}^{dd}$
	to $\Z^{\lfloor (n-1)/2\rfloor }$ that sends $\omega $ to $\bigl(m_{1},\dots
	,m_{\lfloor (n-1)/2\rfloor }\bigr)$ is a bijection and we have the equality
	\begin{equation}
	\lvert \omega \rvert =n\sum_{i=1}^{\lfloor (n-1)/2\rfloor
	}m_{i}^{2}+\sum_{i=1}^{\lfloor (n-1)/2\rfloor }(2i-n)m_{i}. \label{eq:gkdd}
	\end{equation}
\end{Corollary}

Note that in the right-hand side of the above equality seen as a polynomial in $(m_i)_{0\leq i \leq \lfloor (n-1)/2\rfloor }$
 the coefficients of the $m_{i}$'s share the same parity as the leading coefficient.

Since the elements of $\cC_{n}^{dd}:=\cC_{n}\cap \ccdd$ are not
self-conjugate in general, it is relevant to consider \smash{$\cC_{n}^{dd,\mathrm{tr}}=\bigl(\cC_{n}^{dd}\bigr) ^{\mathrm{tr}}$}, the conjugate set of doubled distinct $%
n $ cores. Let~${\phi _{n}\bigl(\omega ^{\mathrm{tr}}\bigr)=(-m_{n-1},\dots ,-m_{0})}$,
with $\omega \in \cC_{n}^{dd}$. Then by setting $m_{i}^{\prime
}=-m_{n-i-1}$, $i=0,\dots n-1$, one gets~${\phi _{n}\bigl(\omega ^{\mathrm{tr}%
}\bigr)\!=\bigl(m_{0}^{\prime },\dots ,m_{n-1}^{\prime }\bigr)}$ with $m_{n-1}^{\prime }=0$
and $m_{i}^{\prime }=m_{i+1}$ for all $i\in \{0,\dots ,\lfloor
(n-1)/2\rfloor -1\}$. Hence, the map from~\smash{$\cC_{n}^{dd,\mathrm{tr}}$}~to $
\Z^{\lfloor (n-1)/2\rfloor }$ that sends $\omega ^{\mathrm{tr}}$ to \smash{$\bigl(m_{0}^{\prime },\dots ,m_{\lfloor (n-1)/2\rfloor -1}^{\prime
}\bigr)=\bigl(m_{1},\dots ,m_{\lfloor (n-1)/2\rfloor }^{\prime }\bigr)$} is a~bijection and
we have this time the equality
\begin{equation*}
\lvert \omega \rvert =n\sum_{i=0}^{\lfloor (n-1)/2\rfloor -1}\bigl(m_{i}^{\prime
}\bigr)^{2}+\sum_{i=0}^{\lfloor (n-1)/2\rfloor -1}(2i-n+2)m_{i}^{\prime }.
\end{equation*}

\subsection{Distinct core partitions and abaci}

Recall that we already introduced the sets of core-partitions $\cC_{n}^{s}$,
$\cC_{n}^{dd}:=\cC_{n}\cap \ccdd$ and $\smash{\cC_{n}^{dd,\mathrm{tr}}}=\smash{\bigl(\cC_{n}^{dd}\bigr) ^{\mathrm{tr}}}$.
In order to get simple expressions of the
atomic lengths as a number of boxes in all types, we will need to consider
new sets of partitions obtained by doubling distinct partitions. They will
not be $2n$-cores in general. To do that, let us introduce the following set
of partitions. Write
\begin{itemize}\itemsep=0pt
	\item $\cC_{n}^{d}$ for the set of distinct partitions $\bar{\lambda}$ such
	that $n\not\in \cch\bigl(\bar{\lambda}\bigr)=\bigl\{h_{s},\, s\in \lambda ,\bar{\lambda}\bar{%
		\lambda}\mid \varepsilon =+\bigr\}$,
	
	\item $\cC_{n}^{d,r}$ be the subset of $\cC_{n}^{d}$ such that $n/2$ is not
	a part of $\bar{\lambda}$,
	
	\item $\cC_{n}^{d,\mathrm{tr}}$ be the set of distinct partitions such that $%
	n\not\in \cch\bigl(\bar{\lambda}^{\mathrm{tr}}\bigr)=\bigl\{h_{s},\, s\in \bar{\lambda}\bar{%
		\lambda}\mid \varepsilon =-\bigr\}$,
	
	\item $\cC_{n}^{d,\mathrm{tr},r}$ be the subset of $\cC_{n}^{d,\mathrm{tr}}$
	such that $n/2\not\in \bar{\lambda}$. Note that the last condition is equivalent to the fact that the main diagonal of \smash{$\bigl(\bar{\lambda}\bar{\lambda}\bigr)^{\mathrm{tr}}$} does not contain a box whose hook length is equal to $n$.
\end{itemize}

By abuse of notation, we will refer to the sets $\cC_{n}^{d,\mathrm{tr}}$
and $\cC_{n}^{d,\mathrm{tr},r}$ as sets of distinct core partitions.
The following lemma obtained from \eqref{eq:hookdd}, shows a simple relation between the sets $%
\cC_{n}^{d,r}$ and $\cC_{n}^{dd}$.

\begin{Lemma}
	\label{lem:dd} Let $n$ be an integer greater or equal to $2$. Take $\bar{%
		\lambda}\in\ccd$. Then $\bar{\lambda}$ belongs to $\cC_{n}^{d,r}$ if and
	only if $\bar{\lambda}\bar{\lambda}$ belongs to $\cC_{n}^{dd}$, the set of
	doubled distinct $n$-core partitions.
\end{Lemma}

\textbf{Warning.} when $\bar{\lambda}$ belongs to $\cC_{n}^{d,%
	\mathrm{tr}}$, its doubled version $\bar{\lambda}\bar{\lambda}$ does not
belong to $\cC_{n}^{dd,\mathrm{tr}}$ in general.

As a consequence of the previous lemma, one derives the following proposition

\begin{Proposition}
	\label{prop:ddbeta} Let $n$ be an integer greater or equal to $2$. Take $%
	\bar{\lambda}\in\cC_{2n+a}^{d,r}$ and $a\in\lbrace 1,2\rbrace$. Set $%
	\phi_{2n+a}\bigl(\bar{\lambda}\bar{\lambda}\bigr)=(m_0,\dots,m_{2n+a-1})$ $($see
	Definition $\ref{defgks})$. Then $\bar{\lambda}$ is in bijective
	correspondence with $\phi_{2n+a}\bigl(\bar{\lambda}\bar{\lambda}\bigr)$ and
	\begin{equation*}
	\bigl\lbrace \bar{\lambda}_1,\dots,\bar{\lambda}_\ell\bigr\rbrace=
	\bigcup_{i=1}^{2n+a-1} \bigl\lbrace (2n+a)k+i,\, 0\leq k\leq m_i-1\mid
	m_i>0\bigr\rbrace.
	\end{equation*}
\end{Proposition}

\begin{proof}
	Let $a$ be either $1$ or $2$. Now let describe the connection between the
	word corresponding to $\bar{\lambda}\bar{\lambda}\in\cC^{dd}_{2n+a}$ and $%
	\bar{\lambda}$. Set $\phi_{2n+a}\bigl(\bar{\lambda}\bar{\lambda}\bigr)=(m_0,m_1,\dots,m_n,m_{n+1},m_{n+2},\dots,m_{2n+a-1})$ (see Definition \ref%
	{defgks}) and $\psi\bigl(\bar{\lambda}\bar{\lambda}\bigr)=(c_k)_{k\in\Z}$ (see
	Definition \ref{def:psi}), then by definition of $\phi_{2n+a}$, $%
	(2n+a)m_i+i=\max\lbrace (k+1)g+i\mid c_{kg+i}=0\rbrace$ for any $0\leq i\leq
	2n+a-1$.
%
 Take $1\leq i\leq n$. If $%
	\max(m_{2n+a-i},m_i)>0$, then Lemma \ref{lem:distinct} guarantees that there
	exists $k\in\lbrace 1,\dots,\ell\rbrace$ such that $\bar{\lambda}%
	_k=\max((2n+a)(m_{2n+a-i}-1)+2n+a-i,(2n+a)(m_{i}-1)+i)$. Moreover, if $m_i=0$%
	, by~Lemma~\ref{lem:distinct}, this is then equivalent to
	\begin{equation*}
	\bigl\lbrace \bar{\lambda}_k,\, k\in\lbrace 1,\dots,\ell\rbrace \mid \bar{%
		\lambda}_k\equiv \pm i\pmod{2n+a}\bigr\rbrace=\varnothing.
	\end{equation*}
	Therefore, we have the set equality
	\begin{equation*}
	\bigl\lbrace \bar{\lambda}_1,\dots,\bar{\lambda}_\ell\rbrace=
	\bigcup_{i=1}^{2n+a-1} \lbrace (2n+a)k+i,\, 0\leq k\leq m_i-1\mid
	m_i>0\bigr\rbrace.\qedhere
\tag*{\qed}
\end{equation*}
\renewcommand{\qed}{}
\end{proof}




As already mentioned, when $\bar{\lambda}$ belongs to~$\cC_{n}^{d,\mathrm{tr%
}}$, its doubled version $\bar{\lambda}\bar{\lambda}$ does not belong to~\smash{$\cC%
_{n}^{dd,\mathrm{tr}}$}. Nevertheless, the core and the quotient of $\bar{%
	\lambda}\bar{\lambda}$ take a very particular form as explained in the
following proposition which can be deduced from Lemma \ref%
{lem:DDprimeLittlewood}.

\begin{Proposition}
	\label{prop:ddprime} Let $n$ be an integer greater or equal to $2$. Take $%
	\bar{\lambda}\in \ccd$ and set $\Phi _{n}\bigl(\bar{\lambda}\bar{\lambda}^{%
		\mathrm{tr}}\bigr)=\smash{\bigl( \omega ,\nu ^{(0)},\dots ,\nu ^{(n-1)}\bigr)} $ $($see
	Definition $\ref{def:Phi})$. Then $\bar{\lambda}$ belongs to \smash{$\cC_{n}^{d,%
		\mathrm{tr},r}$} if and only if
	\begin{itemize}\itemsep=0pt
		\item $\nu^{(j)}=\varnothing$ for any $0\leq j\leq n-2$,
		
		\item let $m=\max\bigl(\lbrace 1\rbrace\cup\bigl\lbrace\bigl(n+\bar{\lambda}_i\bigr)/n \mid
		\bar{\lambda}_i\equiv 0\pmod{n}\bigr\rbrace\bigr)$. Then $\nu^{(n-1)}$ is the rectangle
		${(m-1)\times m}$ partition.
	\end{itemize}
	
	Moreover, in the odd case, setting $\phi _{2n-1}(\omega )=(m_{0},\dots
	,m_{n-2},-m_{n-2},\dots ,-m_{0},0)$ one gets%
	\begin{equation}
	\bigl\lvert \bar{\lambda}\bigr\rvert =\frac{2n-1}{2}\left(
	m^{2}+\sum_{i=0}^{n-2}m_{i}^{2}\right) +\left( -\frac{2n-1}{2}%
	m+\sum_{i=0}^{n-2}\biggl(i-n+\frac{1}{2}\biggr)m_{i}\right) . \label{eq:gkddnuimpair}
	\end{equation}%
	In the even case, by setting $\phi _{2n}(\omega )=(m_{0},\dots
	,m_{n-2},0,-m_{n-2},\dots ,-m_{0},0)$ one gets%
	\begin{equation}
	\bigl\lvert \bar{\lambda}\bigr\rvert =n\left( m^{2}+\sum_{i=0}^{n-2}m_{i}^{2}\right)
	+\left( -nm+\sum_{i=0}^{n-2}(i-n+1)n_{i}\right) .
	\label{eq:gkddnupairbvtquot}
	\end{equation}%
\end{Proposition}

\begin{proof}
	Consider $\bar{\lambda}\in\cC_{n}^{d,\mathrm{tr},r}$, the Littlewood
	decomposition $\Phi_n\bigl(\bar{\lambda}\bar{\lambda}^\mathrm{tr}%
	\bigr)=\bigl(\omega,\nu^{(0)},\dots,\nu^{(n-1)}\bigr)$ (see Definition \ref{def:Phi})
    and $\psi\bigl(\bar{\lambda}\bar{\lambda}^\mathrm{tr}\bigr)=(c_k)_{k\in\Z}$
	its corresponding word (see Definition \ref{def:psi}). Since $n/2\not\in\bar{%
		\lambda}$, the main diagonal $D\bigl(\bar{\lambda}\bar{\lambda}^%
	\mathrm{tr}\bigr)$ of $\bigl(\bar{\lambda}\bar{\lambda}\bigr)^\mathrm{tr}$ has no box of hook length $n$.
	If there exists $i\in\lbrace 0,\dots,n-2\rbrace$ such that $\nu^{(i)}\neq
	\varnothing$, there exists at least one box in $\nu^{(k)}$ whose hook length
	is equal to $1$. By Lemma \ref{lem:indices} and by property $(DD^{\prime }4)$
	from Lemma \ref{lem:DDprimeLittlewood}, there exists $m\in\Z$ such that $%
	c_{mn+i}=1$ and $c_{(m+1)n+i}=0$. Note that the hook length of the box
	$(mn+i,(m+1)n+i)$ is equal to $n$. If $%
	i\not\equiv-i-2\pmod{n}$, that is, if $i\neq n/2$, by \eqref{eq:motddprime},
	one gets that ${c_{-(m+1)n-i-2}=1}$ and ${c_{-mn-i-2}=0}$. Moreover, one of the
	two boxes corresponding to the pairs of indices $(mn+i,(m+1)n+i)$ and $%
	(-(m+1)n-i-2,-mn-i-2)$ is strictly above the main diagonal $D\bigl(\bar{\lambda}\bar{\lambda}^\mathrm{tr}\bigr)$,
	which contradicts the fact that $\smash{\bar{\lambda}\in\cC_{n}^{d,\mathrm{tr},r}}$.
	If $i=n/2$, which implies that $n$ is even in this case, using the same
	arguments, $\smash{\bar{\lambda}\in\cC_{n}^{d,\mathrm{tr},r}}$ implies that there is
	no hook of length $n$ apart maybe on the diagonal~$D\bigl(\bar{\lambda}\bar{\lambda}^\mathrm{tr}\bigr)$. Since that $%
	n/2\not\in \bar{\lambda}$, we have that $\nu^{(n/2)}=\varnothing$. Therefore, \smash{$%
	\bar{\lambda}\in\cC_{n}^{d,\mathrm{tr},r}$} implies that~\smash{${\nu^{(i)}=\varnothing}$}
	for any $0\leq j\leq n-2$.
	
	Now we prove the second part of the statement. Let $d$ be the size of the
	Durfee square of~$\nu^{(n-1)}$. Using the same arguments as above, if $%
	\nu^{(n-1)}$ has a hook length equal to $1$ in any row but the $d+1$-th row,
	the properties of symmetry along the main diagonal and the property~$(DD^{\prime }4)$
    from Lemma \ref{lem:DDprimeLittlewood} imply that $%
	\nu^{(n-1)}$ contains a unique box whose hook length is equal to $1$.
	Therefore, $\nu^{(n-1)}$ is a rectangle of size $m\times(m-1)$ with $%
	m=\max\bigl(\lbrace 1\rbrace\cup\bigl\lbrace\bigl(n+\bar{\lambda}_i\bigr)/n \mid\bar{\lambda%
	}_i\equiv 0\pmod{n}\bigr\rbrace\bigr)$ or equivalently $m=\max(k+1\mid
	c_{kn+n-1}= 0)$. Note that the word to $\nu^{(n-1)}$ is
	\begin{equation*} \label{eq:nu0typeb}
	\psi\bigl(\nu^{(n-1)}\bigr)=\dots 0\underbrace{1\dots 1}_{m-1}\underbrace{0\mid 0\dots
		0}_{m}1\dots.
	\end{equation*}
	
	Equalities \eqref{eq:gkddnuimpair} and \eqref{eq:gkddnupairbvtquot} are
	derived using \eqref{eq:gkdd} and Lemma \ref{lem:DDprimeLittlewood} $%
	(DD^{\prime }3)$.
\end{proof}


Using the same arguments as for Proposition \ref{prop:ddprime}, we obtain
similarly for $\cC_{n}^{d,\mathrm{tr}}$.

\begin{Proposition}
\label{prop:ddprimepair} Let $n$ be an integer greater or equal to $2$. Take
$\bar{\lambda}\in\ccd$ and set the Littlewood decomposition $\Phi_n\bigl(\bigl(\bar{
\lambda}\bar{\lambda}\bigr)^\mathrm{tr}\bigr)=\bigl(\omega,\nu^{(0)},\dots,\nu^{(n-1)}
\bigr)$. Then $\bar{\lambda}$ belongs to $\cC_{n}^{d,\mathrm{tr}} $ if and
only if
\begin{itemize}\itemsep=0pt
\item $\nu^{(j)}=\varnothing$ for any $0\leq j\leq n-2\setminus\bigl\lbrace\frac{n}{
2}\bigr\rbrace$,

\item let $m=\max\bigl(\bigl\lbrace 1\rbrace\cup\lbrace\bigl(n+\bar{\lambda}_i\bigr)/n \mid%
\bar{\lambda}_i\equiv 0\pmod{n}\bigr\rbrace\bigr)$. Then $\nu^{(n-1)}$ is the rectangle partition $(m-1)\times m$. Moreover, when $n$ is even, then \smash{$\nu^{(n/2-1)}$} is the square $m^{\prime }\times m^{\prime }$ with
\[
m^{\prime}=\max\bigl(\lbrace 0\rbrace\cup\bigl\lbrace \bar{\lambda}_i/n \mid\bar{\lambda}_i\equiv n/2\pmod{n}\bigr\rbrace\bigr)
\]
or equivalently $m^{\prime
}=\max(k+1\mid c_{kn+n/2-1}= 0)$.
\end{itemize}

In addition, by setting $\phi_{2n-2}(\omega)=(m_0,\dots,m_{n-3},0,-m_{n-3},\dots,-m_0,0)$, we get
\begin{equation*}
\bigl\lvert \bar{\lambda}\bigr\rvert=(2n-2)\left(\sum_{i=0}^{n-3}m_i^2+m^2+m^{\prime
2}\right)+\left(-(2n-2)m+\sum_{i=0}^{n-2}2(i-n+1)m_i\right).
\end{equation*}
\end{Proposition}


%

The forthcoming subsections use the same arguments. We first use the result
from Section \ref{sec:atomic} giving the correspondence between the affine
Grassmannian elements and some subsets of self-conjugate $(2n)$-cores. Next,
by considering the bijection described in Theorem \ref{thm:gk}, we exhibit
the correspondence between the affine Grassmannian elements and some subsets
of integer partitions so that the atomic length of any affine Grassmannian
element coincides with the weight (i.e., the number of boxes) of the
corresponding partition.

\subsection[Type C\_n\^{}(1)]{Type $\boldsymbol{C_{n}^{(1)}}$}

This is the easiest case since for any $c\in W_{^{C_{n}^{(1)}}}^{0}$ identified to \smash{$\mathcal{C}_{2n}^{s}$}, we have
\[
L_{C_{n}^{(1)}}(c)=L_{A_{2n-1}^{(1)}}(c)= \vert c \vert ,
\] i.e., the
atomic length is given by the number of boxes of $c$ regarded as a
self-conjugate $2n$-core. Therefore, as already seen in (\ref{eq:Lc})%
\begin{equation*}
L_{C_{n}^{(1)}}(c)=2n\sum_{i=0}^{n-1}m_{i}^{2}+%
\sum_{i=0}^{n-1}(2i-2n+1)m_{i}=2n\sum_{i=0}^{n-1}\beta
_{i}^{2}-\sum_{i=0}^{n-1}(2i-1)\beta _{i}. 
\end{equation*}%
Observe also that we have $a_{0}=\sum_{i=0}^{n-1}\beta _{i}^{2}$ and $%
a_{n}=a_{0}-\beta _{1}-\cdots -\beta _{n}$ from which it becomes easy to
compute the atomic length in any type from the results of Section \ref{sec:atomic}.

\subsection[Type D\_\{n+1\}\^{}(2)]{Type $\boldsymbol{D_{n+1}^{(2)}}$}

Recall from Proposition \ref{Prop_D(n+1,2)}, that if \smash{$c\in
W_{^{D_{n+1}^{(2)}}}^{0}=\mathcal{C}_{2n}^{s}$}, we get that
\begin{align*}
L_{D_{n+1}^{(2)}}(c)& =\frac{1}{2}\left( (2n+2)\sum_{i=1}^{n}\beta
_{i}^{2}-\sum_{i=1}^{n}2i\beta _{i}\right) \\
& =(n+1)\sum_{i=1}^{n}\beta
_{i}^{2}-\sum_{i=1}^{n}i\beta _{i}
=(n+1)\sum_{i=0}^{n-1}m_{i}^{2}+\sum_{i=0}^{n-1}(i-n)m_{i}.
\end{align*}
Set \smash{$c\in W_{^{D_{n+1}^{(2)}}}^{0}=\mathcal{C}_{2n}^{s}$} and $\phi
_{2n}(c)=(m_{0},\dots ,m_{n-1},-m_{n-1},\dots ,-m_{0})$ (see Definition \ref{defgks}). Set
\[
c^{\prime }=\phi _{2n+2}^{-1}(0,m_{0},\dots
,m_{n-1},0,-m_{n-1},\dots ,-m_{0}).
\]
 One sees by using the equality
 \[
L_{D_{n+1}^{(2)}}(c)=\frac{1}{2}L_{C_{n}^{(1)}}(c^{\prime })
\] that $c$ could be replaced by $c^{\prime }$ in order to get an atomic length
counted by the number of boxes in a partition. In fact, we have the
following stronger statement.

\begin{Proposition}
	The core $c^{\prime }$ belongs to \smash{$\cC_{2n+2}^{dd}$} and there
	exists $\bar{\lambda}\in \cC_{2n+2}^{d,r}$ such that \smash{$c^{\prime }=\bar{%
		\lambda}\bar{\lambda}$}. Moreover, the map which associates to any \smash{$c\in \cC%
	_{2n}^{s}$} the so obtained distinct partition \smash{$\bar{\lambda}\in \cC%
	_{2n+2}^{d,r}$} is a bijection such that \smash{$L_{D_{n+1}^{(2)}}(c)=\bigl\lvert \bar{\lambda}\bigr\rvert $}.
\end{Proposition}

\begin{proof}
	The first sentence of the proposition is a consequence of Corollary \ref%
	{cor:gkdd} (even case) and Lemma \ref{lem:dd}. Now the map which sends $%
	c\in \cC_{2n}^{s}$ to \smash{$\bar{\lambda}\in \cC_{2n+2}^{d,r}$} is bijective by
	composition. The~equality $L(c)=\bigl\lvert \bar{\lambda}\bigr\rvert $ derives from
	Proposition \ref{prop:ddbeta}, \eqref{eq:gkdd} and \eqref{eq:Lc} setting $%
	\beta _{i}=m_{2n+2-i}$.
\end{proof}

\begin{Example}
	Take $n=2$ and
\[
c=(3,2,1)\in W_{^{D_{3}^{(2)}}}^{0}=\mathcal{C}_{4}^{s}
\]
 and
	$\phi_{4}(c)=(1,-1,1,-1)$. Then we obtain $c'=\phi^{-1}_{6}(0,1,-1,0,1,-1)=(5,3,1,1)$
	and $\bar{\lambda}=(4,1)$. We have \smash{$L_{D_{n+1}^{(2)}}(c)=\lvert c^{\prime
	}\rvert=5$}.
\end{Example}

\begin{Remark}
 Note that in this case the relevant partitions arise from a natural realization of the affine Weyl group. Indeed the group \smash{$W_{^{D_{n+1}^{(2)}}}$} can then be seen as the subgroup of $\widetilde{
	\mathfrak{S}}_{2n+2}$ such that%
\begin{equation*}
W_{^{D_{n+1}^{(2)}}}=\bigl\langle s_{i}^{D_{n+1}^{(2)}},\, i=0,\dots ,n\bigr\rangle
\end{equation*}
	with
\[
s_{i}^{D_{n+1}^{(2)}}=s_{i+1}s_{2n+2-i},\quad i=1,\dots ,n-1 ,\qquad
s_{n}^{D_{n+1}^{(2)}}=s_{n+2}s_{n+1}s_{n+2}, \qquad
s_{0}^{D_{n+1}^{(2)}}=s_0s_1 s_{0}.
\]
With this presentation, one has
\[
W_{^{D_{n+1}^{(2)}}}\cdot \varnothing =\cC_{2n+2}^{dd}
\]
 with the previous action of
$\widetilde{ \mathfrak{S}}_{2n+2}$ on the~$(2n+2)$-cores. In fact, this will be also the case for all the affine types considered later in this paper whose Dynkin diagram is a line (i.e., which do not contain a type $D$ subdiagram).
\end{Remark}

\subsection[Type A\_\{2n\}\^{}(2)]{Type $\boldsymbol{A_{2n}^{(2)}}$}

Recall from Proposition \ref{Prop_A(2n,2)}, that if \smash{$c\in
W_{^{A_{2n}^{(2)}}}^{0}=$} \smash{$\mathcal{C}_{2n}^{s}$}, then by using \eqref{eq:beta_a}, one derives:
\begin{equation*}
L_{A_{2n}^{(2)}}(c)=\frac{1}{2}\left( (2n+1)\sum_{i=1}^{n}\beta
_{i}^{2}-\sum_{i=1}^{n}(2i-1)\beta _{i}\right) .
\end{equation*}%
One can observe that in this case the coefficient of the monomial $\beta _{n}
$ is $-(2n+1)/2$ which happens to be the opposite of the coefficient of $%
\beta _{n}^{2}$. Set $\phi _{2n}(c)=(m_{0},\dots ,m_{n-1},-m_{n-1},\dots
,-m_{0})$ and $c^{\prime }=\phi _{2n+1}^{-1}(0,m_{0},\dots
,m_{n-1},-m_{n-1},\dots ,-m_{0})$ (see Definition \ref{defgks}). Using %
\eqref{eq:gkdd}, $c^{\prime }$ is a~doubled distinct $(2n+1)$-core.
Therefore, by Lemma \ref{lem:dd}, there exists a unique \smash{$\bar{\lambda}\in \cC_{2n+1}^{d,r}$} such that~\smash{$c^{\prime }=\bar{\lambda}\bar{\lambda}$}. Moreover,
we have that
\begin{equation*}
L_{A_{2n}^{(2)}}(c)=\bigl\lvert \bar{\lambda}\bigr\rvert .
\end{equation*}

\begin{Example}
	Take $n=2$ and
\[
c=(3,2,1)\in W_{^{A_{4}^{\prime (2)}}}^{0}=\mathcal{C}%
	_{4}^{s},\qquad   \phi_{4}(c)=(1,-1,1,-1).
\] Then $c^{\prime
		-1}_{5}(0,1,-1,1,-1)=(4,3,1)$ and $\bar{\lambda}=(3,1)$. We have $	L(c^{\prime })=8$.
\end{Example}

\subsection[Type A\_\{2n\}\^{}\{'(2)\}]{Type $\boldsymbol{A_{2n}^{\prime (2)}}$}

Recall from Proposition \ref{Prop_A(2n,2)prime}, that if \smash{$c\in
W_{^{A_{2n}^{(2)}}}^{0}=\mathcal{C}_{2n}^{s}$}, then by using \eqref{eq:beta_a}, one derives
\begin{equation*}
L_{A_{2n}^{\prime (2)}}(c)=(2n+1)\sum_{i=1}^{n}\beta
_{i}^{2}-\sum_{i=1}^{n}2i\beta _{i}.
\end{equation*}
Set
\begin{gather*}
\phi _{2n}(c)=(m_{0},\dots ,m_{n-1},-m_{n-1},\dots ,-m_{0}),\\   c^{\prime }=\phi _{2n+1}^{-1}(m_{0},\dots ,m_{n-1},0,-m_{n-1},\dots ,-m_{0}).
\end{gather*}
The map sending \smash{$c\in \cC_{n}^{s}$} to \smash{$c^{\prime }\in \cC_{2n+1}^{s}$} is a bijection. Moreover, by Corollary \ref{cor:gksc}, we have
\begin{equation*}
L_{A_{2n}^{\prime (2)}}(c)=\lvert c^{\prime }\rvert.
\end{equation*}

\begin{Example}
	Take $n=2$ and
\[
c=(3,2,1)\in W_{^{A_{4}^{\prime (2)}}}^{0}=\mathcal{C}%
	_{4}^{s}, \qquad  \phi _{4}(c)=(1,-1,1,-1).
\] Then $c^{\prime }=\phi
	_{5}^{-1}(1,-1,0,1,-1)=(4,2,1,1)$. We have $L(c)=\vert c^{\prime
	} \vert =8$.
\end{Example}

%
%

\subsection[Type B\_n\^{}(1)]{Type $\boldsymbol{B_{n}^{(1)}}$}

For all the remaining types in this section, the affine Grassmannian
elements are in bijection with the subset \smash{$\cC_{2n}^{s,p}$} of self-conjugate
$(2n)$-cores with an even diagonal. So consider \smash{$c\in \cC_{2n}^{s,p}$} and set  $%
\phi _{2n}(c)=(m_{0},\dots ,m_{n-1},-m_{n-1},\dots ,-m_{0})$. A core has an
even diagonal if and only if the number of horizontal steps after the corner
of the Durfee square is even. This is equivalent to the fact that the number
of letters ``$0$'' of positive index in the corresponding word $\psi (c)$ is
even. Therefore, \smash{$c\in \cC_{2n}^{s,p}$} is equivalent to
\begin{equation*}
\sum_{i=0}^{n-1}\lvert m_{i}\rvert =\sum_{i=0}^{n-1}\lvert \beta _{i}\rvert
\equiv 0\pmod{2}.
\end{equation*}%
One derives from Proposition \ref{Prop_B(n,1)} that for any \smash{$c\in
W_{^{B_{n}^{(1)}}}^{0}=\mathcal{C}_{2n}^{s,p}$}
\begin{equation*}
L_{B_{n}^{(1)}}(c)=\frac{1}{2}\left( 2n\sum_{i=1}^{n}\beta
_{i}^{2}-\sum_{i=1}^{n}2i\beta _{i}\right) =n\sum_{i=1}^{n}\beta
_{i}^{2}-\sum_{i=1}^{n}i\beta _{i}.
\end{equation*}%
Consider \smash{$\bar{\lambda}\in \cC_{2n}^{d,\mathrm{tr},r}$} and set $\Phi _{2n}\bigl(\bigl(%
\bar{\lambda}\bar{\lambda}\bigr)^{\mathrm{tr}}\bigr)=\bigl( \omega ,\nu ^{(0)},\dots
,\nu ^{(n-1)}\bigr)$. By Proposition \ref{prop:ddprime}, we obtain  $\phi
_{2n}(\omega )=(m_{0},\dots ,m_{n-2},0,-m_{n-2},\dots ,-m_{0},0)$ and \smash{$\nu
^{(n-1)}$} is a rectangle $(m-1)\times m$ with~${m\in \N^{\ast }}$. This
permits to define a bijection $f_{2n}$ associating to each $\bar{\lambda}$
its corresponding vector $(m,m_{0},\dots ,m_{n-2})\in \N^{\ast }\times \Z%
^{n-1}$.
We can now consider the map $g\colon \N^{\ast }\times \Z^{n-1} \rightarrow  \Z^{2n} $ defined as follows:
\begin{equation*}
 (m,m_{0},\dots ,m_{n-2})  \mapsto
\begin{cases}
\displaystyle (m,m_{0},\dots ,m_{n-2},-m_{n-2},\dots ,-m_{0},-m)\text{ if }m\equiv
\sum_{i=0}^{n-2}m_{i}\pmod{2}, \\
(-m+1,m_{0},\dots ,m_{n-2},-m_{n-2},\dots ,-m_{0},m-1)\text{ otherwise,}
\end{cases}
\end{equation*}%
which is a bijection from $\N^{\ast }\times \Z^{n-1}$ to
\begin{equation*}
\left\{(m_{0},\dots ,m_{n-1},-m_{n-1},\dots ,-m_{0})\in \Z^{2n}\mid
\sum_{i=0}^{n-1}\lvert m_{i}\rvert \equiv 0\pmod{2}\right\}.
\end{equation*}%
So, the map $\phi _{2n}^{-1}\circ g\circ f_{2n}$ is the desired bijection
from $\cC_{2n}^{d,\mathrm{tr},r}$ to $\cC_{2n}^{s,p}$. Moreover, by
Proposition~\ref{prop:ddprime}, by setting \smash{$\bar{\lambda}=\bigl(\phi
_{2n}^{-1}\circ g\circ f_{2n}\bigr)^{-1}(c)$}, we get
\begin{equation*}
L_{B_{n}^{(1)}}(c)=\bigl\lvert \bar{\lambda}\bigr\rvert .
\end{equation*}

\begin{Example}
	Take $n=3$ and
\[
c=(10,6,6,4,3,3,1,1,1,1)\in W_{^{B_{3}^{(1)}}}^{0}=\mathcal{C}_{6}^{s,p}, \qquad \phi_{6}(c)=(1,-1,-2,2,1,-1).
\]
 Then $g^{-1}\circ\phi_{6}(c)=(1,-1,-2)$. Then $c^{\prime }=(10,6,6,3,3,3,3,1,1,1,1)$ and $\bar{\lambda}=(10,5,4)$. We have \smash{$L_{B_{3}^{(1)}}(c)=\bigl\lvert \bar{\lambda}\bigr\rvert=19$}.
\end{Example}

\subsection[Type A\_\{2n-1\}\^{}(2)]{Type $\boldsymbol{A_{2n-1}^{(2)}}$}

Recall from Proposition \ref{Prop_A(2n-1,2)} that if \smash{$c\in
W_{^{A_{2n-1}^{(2)}}}^{0}=\mathcal{C}_{2n}^{s,p}$}, we have
\begin{equation*}
L_{A_{2n-1}^{(2)}}(c)=\frac{1}{2} \left((2n-1)\sum_{i=1}^n\beta_i^2-\sum_{i=1}^n (2i-1)\beta_i\right).
\end{equation*}

As in the previous subsection, by Proposition \ref{prop:ddprime}, \smash{$\bar{
	\lambda}\in\cC_{2n-1}^{d,\mathrm{tr},r}$} is mapped bijectively to the vector
$(m,m_0,\dots,m_{n-2})\in\N^*\times \Z^{n-1}$. Let us denote this bijection
by $f_{2n-1}$.

Therefore, $c$ is mapped bijectively to \smash{$\bar{\lambda}\in\cC_{2n-1}^{d,%
	\mathrm{tr},r}$} with \smash{$\bar{\lambda}=\bigl(\phi_{2n}^{-1}\circ g\circ
f_{2n-1}\bigr)^{-1}(c)$}. Moreover,
\begin{equation*}
L_{A_{2n-1}^{(2)}}(c)=\bigl\lvert \bar{\lambda}\bigr\rvert.
\end{equation*}

\begin{Example}
	Take $n=3$ and
\[
c=(10,6,6,4,3,3,1,1,1,1)\in W_{^{A_{5}^{(2)}}}^{0}=\mathcal{C%
	}_{6}^{s,p}, \qquad \phi_{6}(c)=(1,-1,-2,2,1,-1).
\]
 Then $g^{-1}\circ%
	\phi_{6}(c)=(1,-1,-2)$. Then $c^{\prime }=(8,5,5,3,3,3,1,1,1)$ and $\bar{%
		\lambda}=(8,4,3)$. We have \smash{$L_{A_{5}^{(2)}}(c)=\bigl\lvert \bar{\lambda}\bigr\rvert=15$}.
\end{Example}

\subsection[Type D\_n\^{}(1)]{Type $\boldsymbol{D_{n}^{(1)}}$}

Recall from Proposition \ref{Prop_D(n,1)} that if \smash{$c\in
W_{^{D_{n}^{(1)}}}^{0}=\mathcal{C}_{2n}^{s,p}$}, we have then by using \eqref{eq:beta_a}
\begin{equation*}
L_{D_{n}^{(1)}}(c)=(n-1)\sum_{i=1}^{n}\beta
_{i}^{2}-\sum_{i=1}^{n}(i-1)\beta _{i}.
\end{equation*}

%
%
%

By Proposition \ref{prop:ddprimepair}, $\bar{\lambda}\in\cC_{2n}^{d,\mathrm{%
		tr}}$ is mapped bijectively to the vector $(m,m_0,\dots,m_{n-3},m^{\prime
})\in\N^*\times \Z^{n-1}\times \N$. Let us denote this bijection by $%
f_{2n-2} $. Therefore, $c$ is mapped bijectively to $\bar{\lambda}\in\cC%
_{2n-2}^{d,\mathrm{tr}}$ with \smash{$\bar{\lambda}=\bigl(\phi_{2n}^{-1}\circ g\circ
f_{2n-2}\bigr)^{-1}(c)$}. Moreover,
\begin{equation*}
L_{D_{n}^{(1)}}(c)=\bigl\lvert \bar{\lambda}\bigr\rvert.
\end{equation*}

\subsection[Type D\_4\^{}(3)]{Type $\boldsymbol{D_{4}^{(3)}}$}

By Proposition \ref{Prop_D(4,3)}, for any
\[
c\in W_{^{D_{4}^{(3)}}}^{0},
\]
 note that \smash{$L_{D_{4}^{(3)}}(c)$} has the same expression as the atom\-ic length of an element of type \smash{$D_3^{(1)}$}, which does not exist. Nevertheless the subset~\smash{$\cC_{2n-2}^{d,\mathrm{tr}}$}
is well defined for $n=3$: it is the subset of distinct partitions such that
$4\not\in \cch\bigl(\bar{\lambda}^{\mathrm{tr}}\bigr)=\{h_{s},\,{s\in \bar{\lambda}\bar{\lambda}},\,{\varepsilon =-}\}$.
	Moreover, the equation \eqref{eq:d43} implies that $\beta_3=\beta_1-\beta_2$.
	By Proposition~\ref{prop:ddprimepair},
\[
c\in W_{^{D_{4}^{(3)}}}^{0}=\mathcal{C}_{6}^{g}
\]
 is mapped bijectively to
	\smash{$\bar{\lambda}\in\cC_{4}^{d,\mathrm{tr}}$} with \smash{$\bar{\lambda}=\bigl(\phi_{6}^{-1}\circ g\circ
f_{4}\bigr)^{-1}(c)$} such that $m=\pm m'\pm m_0-\varphi$, where $m$, $m'$ and $m_0$ are as defined in Proposition \ref{prop:ddprimepair} and $\varphi\in\lbrace 0,1\rbrace$. Moreover,
\begin{equation*}
L_{D_{4}^{(3)}}(c)=\bigl\lvert \bar{\lambda}\bigr\rvert.
\end{equation*}

\subsection[Type G\_2\^{}(1)]{Type $G_{2}^{(1)}$}

Similarly to the previous subsection, from Proposition \ref{Prop_G(2,1)}), note that this time for any
\[
c\in W_{^{G_{2}^{(1)}}}^{0},
\]
\raisebox{1pt}{\smash{$L_{G_{2}^{(1)}}(c)$}} has the same expression as the atomic length of an element type \smash{$B_3^{(1)}$} with an additional restriction.

By Proposition \ref{prop:ddprime},
\[
c\in W_{^{G_{2}^{(1)}}}^{0}
\] is mapped bijectively to \smash{$\bar{\lambda}=\bigl(\phi
_{6}^{-1}\circ g\circ f_{6}\bigr)^{-1}(c)$} with the additional condition that $m$ is equal to $m_0+m_1$ if $m\equiv m_0+m_1 \pmod{2}$ or to $-m_0-m_1-1$ otherwise, where $m$, $m_0$, $m_1$ are defined as in Proposition \ref{prop:ddprime}. Moreover,
\begin{equation*}
L_{G_{2}^{(1)}}(c)=\bigl\lvert \bar{\lambda}\bigr\rvert.
\end{equation*}

\section[More on n-cores and their n-abaci]{More on $\boldsymbol{n}$-cores and their $\boldsymbol{n}$-abaci}\label{sec:abacus}

The goal of this section is first to use the connection between the affine Grassmannian elements and the combinatorial models described in Section \ref{sec:combi} to compute $\varepsilon(c)$ as its appears in Theorem~\ref{thm:macdo}. We next reinterpret the dominant weights appearing in the Weyl characters of the decomposition obtained in Theorem \ref{thm:macdo}. These results can be deduced from those obtained in~\cite{Wmac} and we omit the proofs for short.

Take $n$ a positive integer. Let us introduce for $\lambda\in\ccp$ and $\bar{\lambda}\in\ccd$
\begin{align*}
& H_n(\lambda):=\lbrace s\in\lambda\mid h_s<n\rbrace,\\
& H_n\bigl(\bar{\lambda}\bigr):=\bigl\lbrace h_s \in\cch\bigl(\bar{\lambda}\bigr)\mid h_s<n\bigr\rbrace,\\
& H_n\bigl(\bar{\lambda}^\mathrm{tr}\bigr):=\bigl\lbrace h_s \in\cch\bigl(\bar{\lambda}^{\mathrm{tr}}\bigr)\mid h_s<n\bigr\rbrace,\\
& H_n\bigl(\bar{\lambda}\cup D\bigl(\bar{\lambda}\bar{\lambda}^\mathrm{tr}\bigr)\bigr):=\bigl\lbrace h_s \in\cch\bigl(\bar{\lambda}\bigr)\uplus\bigl\lbrace 2\bar{\lambda}_1,\dots,2\bar{\lambda}_\ell\bigr\rbrace\mid h_s<n\bigr\rbrace,\\
& H_n\bigl(\bar{\lambda}^{\mathrm{tr}}\cup D\bigl(\bar{\lambda}\bar{\lambda}^\mathrm{tr}\bigr)\bigr):=\bigl\lbrace h_s \in\cch\bigl(\bar{\lambda}^{\mathrm{tr}}\bigr)\uplus\bigl\lbrace 2\bar{\lambda}_1,\dots,2\bar{\lambda}_\ell\bigr\rbrace\mid h_s<n\bigr\rbrace.
\end{align*}
We can in fact specify the results of the previous section and use the sets above to compute the signature of the affine Grassmannian elements as detailed in Table~\ref{tableau}.
\begin{table}[th]\centering\renewcommand{\arraystretch}{1.32}
\caption{Table of affine types with their corresponding partition family.}\label{tableau}
\begin{tabular}{|c|c|c|}
\hline
$T$ & $\cC_T$ & $\varepsilon(c)$ \\
\hline
\rule{0pt}{15pt} $A^{(1)}_{n-1}$ & $\cC_n$ &$(-1)^{\#H_n(\lambda)}$\\
\hline
\rule{0pt}{15pt} $B^{(1)}_{n}$& $\mathcal{C}%
_{2n}^{d,\mathrm{tr},r}$ &$(-1)^{\#H_{2n}(\bar{\lambda}^{\mathrm{tr}}\cup D(\bar{\lambda}\bar{\lambda}^\mathrm{tr}))+\ell(\bar{\lambda})}$ \\
\hline
\rule{0pt}{15pt}$A^{(2)}_{2n-1}$ & $\cC_{2n-1}^{d,\mathrm{tr},r}$ &$(-1)^{\#H_{2n-1}(\bar{\lambda)}^{\mathrm{tr}}}$\\
\hline
\rule{0pt}{15pt}$C^{(1)}_{n}$& $\mathcal{C}_{2n}^{s}$ & $(-1)^{\#H_{2n}(\bar{\lambda})}$\\
\hline
\rule{0pt}{15pt}$D^{(2)}_{n+1}$& $\cC_{2n+2}^{d,r}$ &$(-1)^{\#H_{2n+2}(\bar{\lambda}\cup D(\bar{\lambda}\bar{\lambda}^\mathrm{tr}))+\ell(\bar{\lambda})}$\\
\hline
\rule{0pt}{15pt}$A^{(2)}_{2n}$& $\cC_{2n+1}^s$ &$(-1)^{\#H_{2n+1}(\bar{\lambda})}$\\
\hline
\rule{0pt}{15pt}$A^{\prime (2)}_{2n}$& $\cC_{2n+1}^{d,r}$ &$(-1)^{\#H_{2n+1}(\bar{\lambda}\cup D(\bar{\lambda}\bar{\lambda}^\mathrm{tr}))+\ell(\bar{\lambda})}$ \\
\hline

\rule{0pt}{15pt}$D^{(1)}_{n}$& $\mathcal{C}_{2n-2}^{d,\mathrm{tr}}$ & $(-1)^{\#H_{2n-2}(\bar{\lambda}^{\mathrm{tr}}\cup D(\bar{\lambda}\bar{\lambda}^\mathrm{tr}))}$\\
\hline
\end{tabular}
\end{table}

\subsection{Hook length product}

We recall here the following definition that provides a bridge between the abacus model appearing in~\cite{Wmac} and the affine Grassmannian elements.

\begin{Definition}
Let $n$ and $g$ be two positive integers such that $n\leq g$. Set $\lambda \in \ccp$ and $\psi(\lambda)=(c_k)_{k\in\Z}$ its corresponding binary word. For $i\in\lbrace 0,\dots,g-1\rbrace$, define $m_i:=\max\lbrace ({k+1})g+i\mid c_{kg+i}=0\rbrace$. Let $\sigma\colon\lbrace 1,\dots,g\rbrace\rightarrow \lbrace 0,\dots,g-1\rbrace$ be the unique bijection such that $\beta_{\sigma(1)}>\dots>\beta_{\sigma(g)}$. The vector $\mathbf{v}:=\bigl(\beta_{\sigma(1)},\dots,\beta_{\sigma(n)}\bigr)$ is called the $V_{g,n}$-coding corresponding~to~$\lambda$.
\end{Definition}

By definition of the elements of $\cC_T$ in Table \ref{tableau} and Proposition \ref{prop:ddbeta}, one has the following lemma.

\begin{Lemma}
Take $T$ in Table $\ref{tableau}$ such that \smash{$T\neq A_{n-1}^{(1)}$} and $\bar{\lambda}\in\cC_T=\cC_{g}^{a}$. Let $\mu$ be the corresponding partition in either $\ccdd$ $($if $a$ is $d,r)$, $\ccdd^{\mathrm{tr}}$ $($if $a$ is $d,\mathrm{tr},r$ or $a$ is $d,\mathrm{tr})$ or $\ccsc$ $($in the other cases$)$. Let $\mathbf{v}:=\bigl(\beta_{\sigma(1)},\dots,\beta_{\sigma(n)}\bigr)$ be the $V_{g,n}$-coding corresponding to $\mu$. Then for any $1\leq i\leq n$ such that $v_i>g$, we have
\begin{itemize}\itemsep=0pt
\item if $\mu$ is a doubled distinct partition, $v_i=\max\bigl\lbrace \bar{\lambda}_k+g,\, k\geq 1,\,\bar{\lambda}_k\equiv \sigma(i)\pmod{g} \bigr\rbrace$
\item otherwise, $v_i=\max\lbrace \bar{\lambda}_k-1,\, k\geq 1,\,\bar{\lambda}_k-1+g\equiv \sigma(i)\pmod{g} \rbrace$.

\end{itemize}
\end{Lemma}

Recall from Theorem \ref{thm:macdo} that if $c\in W_{a}^{0}$, we can set $c=ut_{\nu}$ with $\nu\in M^{*}\cap(-P_{+})$ and $u\in W^{\nu}$ for
any element of the affine Grassmannian $W_{a}^{0}$. Moreover, each weight $%
-\eta^{\vee}\nu+u^{-1}(\mathring{\rho})-\mathring{\rho}$ belongs to $P_{+}$,
the set of dominant weights for $\mathring{A}$. The purpose of the following theorem is to fill the gaps between our reformulation of the Weyl--Kac formula as a sum over elements of the affine Grassmannian and the formulation of the same formula within the language of multivariate series and $\theta$ function (see~\cite[Proposition $6.1$]{RS}), where the sum is over $\Z^n$. Given Remark \ref{rem:reseau}, we set for any affine algebra $\mathfrak{g}$ of rank $n$ from Table \ref{table2}:
\[\tilde{\eta}^{\vee}=\begin{cases} \eta^{\vee} &\text{if } \mathfrak{g}\neq C_n^{(1)},\\
2 \eta^{\vee} &\text{if } \mathfrak{g}= C_n^{(1)}.\end{cases}\]

Note that $\tilde{\eta}^{\vee}$ corresponds to the parameter $g$ in~\cite{Mac,Wmac}.

\begin{Theorem}
With the previous notation the one-to-one correspondence described in Section~$\ref{sec:combi}$ associates to any affine Grassmannian element $c$ its corresponding partition $\lambda$ of $\cC_T$ in Table~$\ref{tableau}$ so that by writing $\mathbf{v}=(v_i)_{1\leq i\leq n}$ for the $V_{\tilde{\eta}^{\vee},n}$ coding of $\lambda$, we get:
\[-\tilde{\eta}^{\vee}\nu+u^{-1}(\mathring{\rho})-\mathring{\rho}=\bigl(v_i+i-\tilde{\eta}^{\vee}\bigr)_{1\leq i\leq n}.\]
\end{Theorem}

Recall that the proof of the Nekrasov--Okounkov formula by Han~\cite{Han08} heavily relies on a~specialization of the Weyl--Kac in type \smash{$A_{n-1}^{(1)}$} and a polynomial argument. This approach can be extended to derive ($u$-analogues of) Nekrasov--Okounkov formulas for any classical non-exceptional affine type $T$. In order to do so from Theorem \ref{thm:macdo}, one needs to be able to convert the characters on the right-hand side of Theorem \ref{thm:macdo} (or more precisely their associated $u$-dimension formulas) as a product over hook lengths of elements of $\cC_T$. Actually one can reformulate the results from~\cite[Section $4.2$]{Wmac} for a set of Laurent variables $(X_k)_{k\in \Z}$ as the bridge between hook length products on the one hand and products over weights of the form $-\eta ^{\vee }\nu +u^{-1}(\mathring{\rho})-\mathring{\rho}$ associated with an affine Grassmannian element as in Theorem \ref{thm:macdo}.

\begin{Theorem}\label{thm:enumhook}
Let $(X_k)_{k\in\Z}$ be formal variables.
Consider an affine type $T$ of rank $n$ in Table~$\ref{tableau}$ and any element of the affine Grassmannian $c\in W_{a}^{0}$ such that $c=ut_{\nu}$ with $\nu\in Q\cap(-P_{+})$ and~$u\in W^{\nu}$. Set $\mathbf{r}=(r_i)_{1\leq i\leq n}$ the vector such that $r_i=\bigl(-\tilde{\eta}^{\vee}\nu+u^{-1}(\mathring{\rho})-\mathring{\rho}\bigr)_i$ for all $i$. Let $\lambda\in\cC_T$ and $\cch(\lambda)_T$ its corresponding elements in Tables~$\ref{tableau}$ and~$\ref{tableau2}$, respectively. Set
\[
\alpha_{T}(i)=\#\bigl\lbrace h_s\in \cch(\lambda)_T , \, h_s=\tilde{\eta}^{\vee}-i\bigr\rbrace
\]
 and where
 \[
 \alpha'_{T}(i)=\alpha_T(i)+\#\bigl\lbrace h\in \bar{\lambda},\, 2h=\tilde{\eta}^{\vee}-i\bigr\rbrace+\#\bigl\lbrace h\in \bar{\lambda},\, h=\tilde{\eta}^{\vee}-i\bigr\rbrace.
\]
Then we have
\begin{itemize}\itemsep=0pt
\item if \smash{$T=A_{n-1}^{(1)}$},
\begin{equation}\label{eq:hookproductA}
\prod_{h\in\cch(\lambda)_T}\frac{X_{h-\tilde{\eta}^{\vee}}X_{h+\tilde{\eta}^{\vee}}}{X_h^2}=\prod_{i=1}^{\tilde{\eta}^{\vee}-1}\left(\frac{X_{-i}}{X_i}\right)^{\alpha_{T}(i)}\Delta_T(\mathbf{r}),
\end{equation}
\item if \smash{$T\in C_n^{(1)}\cup D_{n+1}^{(2)}\cup A_{2n}^{ (2)}\cup A_{2n}^{\prime (2)}$},
\begin{equation*}
\prod_{h\in\bar{\lambda}}\frac{X_{2h-\tilde{\eta}^{\vee}}X_{h-\tilde{\eta}^{\vee}}}{X_{2h}X_h}
\prod_{h\in\cch(\lambda)_T}\frac{X_{h-\tilde{\eta}^{\vee}}X_{h+\tilde{\eta}^{\vee}}}{X_h^2}=\prod_{i=1}^{\tilde{\eta}^{\vee}-1}\left(\frac{X_{-i}}{X_i}\right)^{\alpha'_{T}(i)}\Delta_T(\mathbf{r}),
\end{equation*}
\item if \smash{$T\in B_n^{(1)}\cup A_{2n-1}^{(2)}$},
\begin{equation*}
\prod_{h\in\bar{\lambda}}\frac{X_{2h-\tilde{\eta}^{\vee}}X_{h+\tilde{\eta}^{\vee}}}{X_{2h}X_h}\prod_{h\in\cch(\lambda)_T}\frac{X_{h-\tilde{\eta}^{\vee}}X_{h+\tilde{\eta}^{\vee}}}{X_h^2}=\prod_{i=1}^{\tilde{\eta}^{\vee}-1}\left(\frac{X_{-i}}{X_i}\right)^{\alpha'_{T}(i)}\Delta_T(\mathbf{r}),
\end{equation*}
\item if \smash{$T\in D_n^{(1)}$},
\begin{equation*}
\prod_{h\in\bar{\lambda}}\frac{X_{2h+\tilde{\eta}^{\vee}}X_{h+\tilde{\eta}^{\vee}}}{X_{2h}X_h}\prod_{h\in\cch(\lambda)_T} \frac{X_{h-\tilde{\eta}^{\vee}}X_{h+\tilde{\eta}^{\vee}}}{X_h^2}=\prod_{i=1}^{\tilde{\eta}^{\vee}-1} \left(\frac{X_{-i}}{X_i}\right)^{\alpha'_{T}(i)}\Delta_T(\mathbf{r}).
\end{equation*}
\end{itemize}
\end{Theorem}

\begin{table}[th]\centering\renewcommand{\arraystretch}{1.55}
\caption{Table of affine types with their corresponding hook length product.}\label{tableau2}
	\begin{tabular}{|c|c|c|}
		\hline
		$T$ & $\cch(\lambda)_T$ & $\Delta_T(\mathbf{r})$ \\
		\hline
		\rule{0pt}{17pt} $A^{(1)}_{n-1}$ & $\cch(\lambda)$ &
		$\displaystyle\prod_{1\leq i<j\leq n} \frac{X_{r_i-r_j}}{X_{j-i}}$\\
		\hline
		\rule{0pt}{17pt} $B^{(1)}_{n}$& $\cch\bigl(\bar{\lambda}^{\mathrm{tr}}\bigr)$ &$\displaystyle\prod_{1\leq i<j\leq n} \frac{X_{r_i-r_j}X_{r_i+r_j}}{X_{j-i}X_{\tilde{\eta}^{\vee}+2-i-j}}$ \\
		\hline
		\rule{0pt}{17pt}$A^{(2)}_{2n-1}$ & $\cch\bigl(\bar{\lambda}^{\mathrm{tr}}\bigr)$ &$\displaystyle \prod_{i=1}^{n}\frac{X_{r_i}}{X_i}\Delta_{B_n^{(1)}}(\mathbf{r})$\\
		\hline
		\rule{0pt}{17pt}$C^{(1)}_{n}$& $\cch\bigl(\bar{\lambda}^{\mathrm{tr}}\bigr)$ &$\displaystyle \prod_{i=1}^{n}\frac{X_{r_i}}{X_i}\displaystyle\prod_{1\leq i<j\leq n} \frac{X_{r_i-r_j}X_{r_i+r_j}}{X_{j-i}X_{\tilde{\eta}^{\vee}-i-j}}$\\
		\hline
		\rule{0pt}{17pt}$D^{(2)}_{n+1}$& $\cch\bigl(\bar{\lambda}^{\mathrm{tr}}\bigr)$ &$\displaystyle\prod_{1\leq i<j\leq n} \frac{X_{r_i-r_j}X_{r_i+r_j}}{X_{j-i}X_{\tilde{\eta}^{\vee}+1-i-j}}$\\
		\hline
		\rule{0pt}{17pt}$A^{(2)}_{2n}$& $\cch\bigl(\bar{\lambda}^{\mathrm{tr}}\bigr)$ & $\displaystyle \prod_{i=1}^{n}\frac{X_{r_i}}{X_i}\Delta_{D_{n+1}^{(2)}}(\mathbf{r})$\\
		\hline
		\rule{0pt}{17pt}$A^{\prime (2)}_{2n}$& $\cch\bigl(\bar{\lambda}^{\mathrm{tr}}\bigr)$ &$\displaystyle\prod_{1\leq i<j\leq n} \frac{X_{r_i-r_j}X_{r_i+r_j}}{X_{j-i}X_{\tilde{\eta}^{\vee}-i-j}}$ \\
		\hline
		
		\rule{0pt}{15pt}$D^{(1)}_{n}$& $\displaystyle \cch\bigl(\bar{\lambda}^{\mathrm{tr}}\bigr)$ & $\displaystyle \prod_{i=1}^{n}\frac{X_{2r_i}X_i}{X_{2i}X_{r_i} }\Delta_{B_n^{(1)}}(\mathbf{r})$\\
		\hline
	\end{tabular}
\end{table}

\begin{Remark}
The case \smash{$T=A_{n-1}^{(1)}$} in the above theorem is a reformulation of~\cite[Theorem 1]{HD}. The hook length products in the above table are obtained from the principal specialization of the Schur functions appearing in Theorem \ref{thm:macdo}. Usually, this principal specialization in Schur functions is rather expressed in terms of hook-contents in Young diagrams (see~\cite[Theorem $15.3$]{Stanleyhook}). Thus one can ask about the connection between hook-content formulas and the hook length formulas of $n$-cores. Following the notation from Theorem \ref{thm:enumhook}, the hook-content reformulation of \eqref{eq:hookproductA} which corresponds to~\cite[Theorem $17$]{HD} is
\[
\prod_{1\leq i<j\leq n} \frac{X_{r_i-r_j}}{X_{j-i}}=\prod_{s\in\mu} \frac{X_{\tilde{\eta}^{\vee}+c_s}}{X_{h_s}},
\]
where $\mu$ is the partition whose parts are the integers $(r_i-r_n)$ with $1\leq i\leq n$ and $c_s$ is the content of the box $s\in\mu$.
One can show that similar results hold for other types.
\end{Remark}

\begin{Remark}
In the left-hand side products of Theorem \ref{thm:enumhook}, one can note that if we take $(X_k)_{k\in\Z}=(k)_{k\in\Z}$, the product of hook lengths is equal to $0$ for any partition $\lambda\in\ccp\setminus \cC_T$. For instance, in type \smash{$A_{n-1}^{(1)}$}, the product over all the hook lengths $\bigl(1-n^2/h^2\bigr)$ of a partition $\lambda$ is equal to $0$ if and only if $\lambda$ is not an $n$-core.
\end{Remark}

The remaining steps leading to Nekrasov--Okounkov type identities are to consider Theorem~\ref{thm:macdo} as an equality of multivariate power series setting $x_i={\rm e}^{-\varepsilon_i}$ for $1\leq i\leq n$. Then we specialize the set of variables $(x_i)$ so that we can use the Weyl $u$-dimension formula, then we prove that the Nekrasov--Okounkov type identities hold for an infinite number of ranks and conclude using a polynomial argument. Theorem \ref{thm:enumhook} is the cornerstone to provide the connection between the Weyl character formula on the one hand and the hook length products on the other.\looseness=1

Let $u$ be a formal variable. The specializations of the variables $x_i$'s as well as the corresponding hook length product and the specialization of the Laurent variables used in Theorem~\ref{thm:enumhook} are summarized in the following table, where $\varepsilon(c)$ is the signature of the corresponding affine Grassmannian element (see Table~\ref{tableau}).

\begin{table}[th]\centering\renewcommand{\arraystretch}{1.52}
\caption{Table of affine types with their corresponding specializations and hook length product.}\label{tableau3}

\begin{tabular}{|c|c|c|c|}
\hline
$T$ & $(x_i)_{1\leq i\leq n}$ & $X_k$ & hook length for the group character $s_{-\tilde{\eta}^{\vee }\nu +u^{-1}(\mathring{\rho})-\mathring{\rho}}$\\
\hline
\rule{0pt}{15pt} $A^{(1)}_{n-1}$ & $\bigl(u^{i-1}\bigr)$ & $1-u^k$&
$\varepsilon(c)\displaystyle\prod_{h\in\cch(\lambda)_T} \frac{(1-u^{h+\tilde{\eta}^{\vee}})(1-u^{h-\tilde{\eta}^{\vee}})}{(1-u^h)^2}$\\
\hline
\rule{0pt}{15pt} $B^{(1)}_{n}$& $\bigl(u^{2i-1}\bigr)$ & $1-u^{2k}$&
\makecell{
$\varepsilon(c)\displaystyle u^{-\tilde{\eta}^{\vee}\ell\bigl(\bar{\lambda}\bigr)}\prod_{h\in\cch(\lambda)_T} \frac{(1-u^{2(h-\tilde{\eta}^{\vee})})(1-u^{2(h+\tilde{\eta}^{\vee})})}{(1-u^{2h})^2}$\vspace{1mm}\\ $\displaystyle\times\prod_{h\in\bar{\lambda}}\frac{(1+u^{2h-\tilde{\eta}^{\vee}})(1-u^{2(h+\tilde{\eta}^{\vee})})(1-u^{2h+\tilde{\eta}^{\vee}})}{(1-u^{4h})(1-u^{2h})}$}
\\
\hline
\rule{0pt}{15pt}$A^{(2)}_{2n-1}$ & $\bigl(u^{i}\bigr)$ & $1-u^{k}$&\makecell{
$\varepsilon(c)(-1)^{\ell\bigl(\bar{\lambda}\bigr)}\displaystyle u^{-\ell\bigl(\bar{\lambda}\bigr)\tilde{\eta}^{\vee}/2}\prod_{h\in\cch(\lambda)_T} \frac{(1-u^{h-\tilde{\eta}^{\vee}})(1-u^{h+\tilde{\eta}^{\vee}})}{(1-u^{h})^2}$\vspace{1mm}\\
$\displaystyle \times
\prod_{h\in\bar{\lambda}}\frac{(1-u^{h-\tilde{\eta}^{\vee}/2})(1-u^{h+\tilde{\eta}^{\vee}}) (1+u^{h+\tilde{\eta}^{\vee}/2})}{(1-u^{2h})(1-u^{h})(1+u^{h-\tilde{\eta}^{\vee}/2})}$}
\\
\hline
\rule{0pt}{15pt}$C^{(1)}_{n}$& $\bigl(u^{i}\bigr)$ & $1-u^{k}$&\makecell{
$\displaystyle \varepsilon(c)u^{\ell\bigl(\bar{\lambda}\bigr)\tilde{\eta}^{\vee}/2}
\prod_{h\in\cch(\lambda)_T}\frac{(1-u^{h-\tilde{\eta}^{\vee}})(1-u^{h+\tilde{\eta}^{\vee}})}{(1-u^{h})^2}
$\vspace{1mm}\\  $\displaystyle\times\prod_{h\in\bar{\lambda}}\frac{(1-u^{h-\tilde{\eta}^{\vee}/2})(1-u^{h-\tilde{\eta}^{\vee}})(1+u^{h+\tilde{\eta}^{\vee}/2})}{(1-u^{2h})(1-u^{h})}$}\\
\hline
\rule{0pt}{15pt}$D^{(2)}_{n+1}$& $\bigl(u^{2i-1}\bigr)$ & $1-u^{2k}$& \makecell{
$\varepsilon(c)(-1)^{\ell\bigl(\bar{\lambda}\bigr)}\displaystyle \prod_{h\in\cch(\lambda)_T}\frac{(1-u^{2(h-\tilde{\eta}^{\vee})})(1-u^{2(h+\tilde{\eta}^{\vee})})}{(1-u^{2h})^2}$\vspace{1mm}\\
$\displaystyle \times \prod_{h\in\bar{\lambda}}\frac{(1+u^{2h-\tilde{\eta}^{\vee}})(1-u^{2h+\tilde{\eta}^{\vee}})}{(1-u^{4h})}$}
\\
\hline
\rule{0pt}{15pt}$A^{(2)}_{2n}$& $\bigl(u^{2i-1}\bigr)$ & $1-u^{2k}$&\makecell{
$\varepsilon(c)\displaystyle u^{\tilde{\eta}^{\vee}\ell\bigl(\bar{\lambda}\bigr)} \prod_{h\in\cch(\lambda)_T}\frac{(1-u^{2(h-\tilde{\eta}^{\vee})})(1-u^{2(h+\tilde{\eta}^{\vee})})}{(1-u^{2h})^2}$\vspace{1mm}\\ $\displaystyle\times\prod_{h\in\bar{\lambda}}\frac{(1+u^{2h-\tilde{\eta}^{\vee}})(1-u^{2(h-\tilde{\eta}^{\vee})})(1-u^{2h+\tilde{\eta}^{\vee}})}{(1-u^{4h})(1-u^{2h})}$}
\\
\hline
\rule{0pt}{15pt}$D^{(1)}_{n}$& $\bigl(u^{2i-1}\bigr)$ & $1-u^{2k}$&\makecell{
$\varepsilon(c)(-1)^{\ell\bigl(\bar{\lambda}\bigr)}\displaystyle u^{-2\tilde{\eta}^{\vee}\ell\bigl(\bar{\lambda}\bigr)}\prod_{h\in\cch(\lambda)_T} \frac{(1-u^{2(h-\tilde{\eta}^{\vee})})(1-u^{2(h+\tilde{\eta}^{\vee})})}{(1-u^{2h})^2}
$
\vspace{1mm}\\
$\displaystyle\times \prod_{h\in\bar{\lambda}}\frac{(1+u^{2h-\tilde{\eta}^{\vee}})(1-u^{2(h+\tilde{\eta}^{\vee})}) (1-u^{2h+\tilde{\eta}^{\vee}})}{(1-u^{4h})(1-u^{2h})}$}
\\
\hline
\end{tabular}
\end{table}

To get the $u$-deformation of the Nekrasov--Okounkov formula as stated in \eqref{Hande} (and its generalizations in any affine classical type), we start with Theorem \ref{thm:macdo} and specialize each group character appearing in the right-hand side as prescribed in the second column of Table \ref{tableau3}. We obtain the right-hand side of the equalities of Theorem \ref{thm:enumhook}.
Then, one gets the hook length formulas of the fourth column by using Theorem \ref{thm:enumhook}. Observe the signatures $\epsilon(c)$ of the affine Grassmannian elements simplify. In affine type $A$ this yields \eqref{Hande}. For the other types, we refer to~\cite[Section 5.3]{Wmac} for the precise formulas.

 \textbf{Warning.} Theorem \ref{thm:macdo} uses the \textit{dual} atomic length $L^{\vee}$. This subtlety forces to each affine root system to be replaced by its dual (obtained by transposing its generalized Cartan matrix) between Table~\ref{tableau} and Tables \ref{tableau2} and \ref{tableau3}.

\subsection{Root systems and hook length}

The integers $\alpha_T(i)$ and $\alpha'_T(i)$ in Table \ref{tableau2} can be expressed as a refinement of the cardinality of affine inversion sets. Here we detail the case of type \smash{$A_{n-1}^{(1)}$}. Similar interpretations can be derived for the other affine classical types. Our goal in this subsection is to suggest a direct approach between the integer partitions of Section \ref{sec:combi} and the affine root systems, or equivalently between the sublattices of Section \ref{sec:4} and the integer partitions introduced in Section \ref{sec:combi}.

Take \smash{$w\in\widetilde{\mathfrak{S}}_{n}$}. Recall that $R(w):=\big\lbrace\alpha\in R_+\mid w^{-1}(\alpha)\in-R_+\bigr\rbrace$ is the inversion set of $w$ with{\samepage
\[R_+=\bigl\lbrace \alpha+k\delta\mid \text{where } k\in \N \text{ if } \alpha\in R_+^{0},  \ \text{and} \
k\in \N^* \text{ if } \alpha\in -R_+^{0}\bigr\rbrace,\]
where $R_+^{0}$ is the set of positive roots of $\mathfrak{sl}_n$.}

One can also represent any element of the affine symmetric group $\widetilde{\mathfrak{S}}_{n}$ in its window notation $w=[w(1),\dots ,w(n)]$ such that%
\begin{equation*}
w(1)+\cdots +w(n)=\frac{n(n+1)}{2}
\end{equation*}%
and for any $i\neq j$%
\begin{equation*}
w(i)\not\equiv w(j)\func{mod}n.
\end{equation*}%
Recall we have a surjective group morphism%
\begin{equation*}
\pi \colon \
\widetilde{\mathfrak{S}}_{n}\rightarrow \mathfrak{S}_{n} ,\qquad
w=[w(1),\dots ,w(n)]\mapsto \pi (w)=(\overline{w}(1),\dots ,\overline{w}%
(n)),
\end{equation*}
where for any $i=1,\dots ,n$, $1\leq \overline{w}(i)\leq n$ is such that $%
\overline{w}(i)=w(i)\func{mod}n$.

The window notation of a core $c=[c(1),\dots ,c(n)]$ seen as an affine
Grassmannian element and its associated vector $(\beta _{1},\dots ,\beta
_{n})$ are related by the formulas%
\begin{equation*}
\beta _{i}=\frac{c(i)-\overline{c}(i)}{n},\qquad i=1,\dots ,n
\end{equation*}%
and we then have
\begin{equation*}
c(1)<\cdots <c(n).
\end{equation*}%
More generally $\ker \pi $ is isomorphic to the root lattice $Q$ of type $%
A_{n-1}$ seen as an abelian group.

We use an embedding of the root system of $A_\infty$ into the one of \smash{$A_{n-1}^{(1)}$}. To do so, let us consider~${E=\bigoplus_{k\in\Z} \Z\varepsilon_k}$ which can be seen as the set of sequences $u=(u_k)_{k\in\Z}$ whose support is finite. Set $F=\bigoplus_{i=1}^n \Z\varepsilon_i\oplus \Z\delta$ of rank $n+1$. Let us introduce the $\Z$-linear surjective map $\theta_n$ defined as follows:
\[
\theta_n\colon \ E   \to   F,\qquad
\varepsilon_k  \mapsto  \varepsilon_i-a\delta ,
\]
where $k=i+an$ with $i\in\lbrace 1,\dots,n\rbrace$.

Remark first that $\theta_n(\varepsilon_0-\varepsilon_n)=\delta$.

Moreover, $u=(u_k)_{k\in\Z}$ belongs to $\ker \theta_n$ if
\begin{align}
&\forall i\in\lbrace 1,\dots,n\rbrace,\qquad\sum_{a=-\infty}^{\infty} u_{i+an}=0,\label{condition1}\\
&\sum_{a=-\infty}^{\infty} -a(u_{1+an}+\dots+u_{n+an})=0.\label{condition2}
\end{align}

The affine symmetric group $\widetilde{\mathfrak{S}}_{n}$ acts naturally on $E$:
\[w\dot (u_k)_{k\in\Z}=\bigl(u_{w(k)}\bigr)_{k\in\Z}.\]
\begin{Lemma}\label{lem:stability}
The set $\ker \theta_n$ is stable under the action of $\widetilde{\mathfrak{S}}_{n}$.
\end{Lemma}
The proof of the above lemma follows from \eqref{condition1} and \eqref{condition2} above, the fact that elements of~$\widetilde{\mathfrak{S}}_{n}$ are $n$ periodic and
\[\lbrace w(1),\dots ,w(n)\rbrace\equiv\lbrace 0,\dots,n-1\rbrace\pmod{n}.\]

We have the following proposition.

\begin{Proposition}
The action of $\widetilde{\mathfrak{S}}_{n}$ on $E$ induces an action of $\widetilde{\mathfrak{S}}_{n}$ on $F$ by setting for any~$y\in F$:
\[w(y)=w(x), \qquad\text{for any } x\in \theta_n^{-1}\lbrace y\rbrace.\]
Set $Q:=\sum_{i=1}^{n}\Z(\varepsilon_i-\varepsilon_{i+1})\oplus\Z\delta$. The restriction of this action coincides with the action of~$\widetilde{\mathfrak{S}}_{n}$~on the root lattice of \smash{$A_{n-1}^{(1)}$}.
\end{Proposition}
\begin{proof}
The proof of the independence of the choice of $w\in\theta_n^{-1}\lbrace y\rbrace$ follows by Lemma \ref{lem:stability} and the surjectivity of $\theta_n$.

In order to check that the restriction of $\widetilde{\mathfrak{S}}_{n}$ to the lattice $Q$ actually corresponds to the affine action, we have first to check that
\[ w(\delta)=w(\varepsilon_0-\varepsilon_n)=\varepsilon_{w(0)}-\varepsilon_{w(n)}=\delta.\]
It remains to check that the generators $s_0,s_1,\dots,s_{n-1}$ acts as expected on $\delta$ and $\varepsilon_i-\varepsilon_{i+1}$ for any $1\leq i\leq n-1$.
\end{proof}

Now consider $w$ in $\widetilde{\mathfrak{S}}^0_n:=\widetilde{\mathfrak{S}}_{n}/\mathfrak{S}_{n}$. Take $(a,b)\in\lbrace 1,\dots n\rbrace^2$ such that $a\neq b$ and a positive integer $k$ such that $\alpha=\varepsilon_a-\varepsilon_b+k\delta\in R_+$. Recall that the height of $\alpha$ $ht(\alpha)$ is $a-b+kn$ and set $ht_n(\alpha)=a-b$, the height of $\alpha$ modulo $\delta$. Take $(i,j)\in\lbrace 1,\dots,n\rbrace^2$ the unique elements such that $a\equiv w(i)\pmod{n}$ and $b\equiv w(j)\pmod{n}$. Define $(m_a,m_b)\in\Z^2$ such that $a=w(i)+m_an$ and $b=w(j)+m_bn$. Thus
\[\alpha=\varepsilon_{w(i)}-\varepsilon_{w(j)}+(k+m_b-m_a)\delta.\]

Set $k'=k+m_b-m_a$. Our hypothesis is that $\alpha$ is a positive root. It implies that $k=k'+m_a-m_b$ is greater than or equal to $\mathds{1}_{a>b}(a,b)$.
Moreover, $\alpha\in R(w)$ if and only if
\[w^{-1}(\alpha)=\varepsilon_i-\varepsilon_j+k'\delta\in -R_+.\]
The above condition is equivalent to the fact that $\varepsilon_i-\varepsilon_j+(k+m_b-m_a)\delta$ belongs to $-R_+$. It~implies that $m_b-m_a\leq -\mathds{1}_{a>b}(a,b)$. Therefore, $a-b+(m_b-m_a)n\leq-\mathds{1}_{a>b}(a,b)+a-b<0$. Thus $w(i)<w(j)$. Given that $w\in\widetilde{\mathfrak{S}}^0_n$, this is equivalent to $i<j$.

Combining these conditions, it follows that $k$ satisfies the inequality
\begin{equation}\label{eq:ineg}
\mathds{1}_{a>b}(a,b)\leq k\leq m_a-m_b-1.
\end{equation}

 Let $n \geq 2$ be an integer and consider
\[
J_n\colon \   \widetilde{\mathfrak{S}}^0_n   \to   \lbrace 0,1\rbrace^{\Z}, \qquad
  w   \mapsto   (c_k)_{k\in\Z},
\]
where $c_k=1$ if $k\geq w(i)-1$ and $0$, with $i\in\lbrace 1,\dots,n\rbrace$ such that $k\equiv w(i)-1\pmod{n}$.

Note that for any $w\in\widetilde{\mathfrak{S}}^0_n$, the element $\psi^{-1}(J_n(w))$ belongs to $\mathcal{C}_n$, where $\psi$ is the bijection defined in Definition \ref{def:psi}. It can be seen that this map is actually a bijection and it is illustrated in the Example \ref{ex:exfin}.

Moreover, by definition of $V_{n,n}$-coding, $c(n-i)=v_i+1$ for any $i\in\lbrace 1,\dots,n\rbrace$.

We are now ready to state the following proposition.
\begin{Proposition}\label{prop:inversion}
Set $w\in\widetilde{\mathfrak{S}}^0_n$ and $\lambda=\Psi^{-1}(J_n(w))$ the corresponding $n$-core partition. Define the set $R_i(w):=\bigl\lbrace \alpha\in R(w),\, ht_n\bigl(w^{-1}(\alpha)\bigr)=n-i\bigr\rbrace$.
Then for any $i\in\lbrace 1,\dots,n-1\rbrace$, we have
\[\#R_i(w)=\#\lbrace h\in \cch(\lambda),\, h=n-i\rbrace.\]
\end{Proposition}

\begin{proof}
Set $w\in\widetilde{\mathfrak{S}}^0_n$, $J_n(w)=(c_k)_{k\in\Z}$ and $\lambda=\Psi^{-1}(J_n(w))$ the corresponding $n$-core partition. Set $\alpha\in R_i(w)$ and $(a,b)\in\lbrace 1,\dots,n\rbrace^2$ such that $\alpha=\varepsilon_a-\varepsilon_b+k\delta$
and $w^{-1}(a)-w^{-1}(b)\equiv n-i\pmod{n}$. Take $(c,d)\in\lbrace 1,\dots,n\rbrace$ the unique elements such that $a\equiv w(c)\pmod{n}$ and $b\equiv w(d)\pmod{n}$. Take $(m_a,m_b)\in\Z^2$ such that $a=w(c)+m_an$ and $b=w(d)+m_bn$. Thus
\[\alpha=\varepsilon_{w(c)}-\varepsilon_{w(d)}+(k+m_b-m_a)\delta.\]
Hence, by definition of $J_n$, we have
\begin{align*}
&\lbrace k\in \Z,\, c_{k}=1,\, k\equiv w(c)-1\pmod{n}\rbrace=\lbrace a-1+qn,\, q\geq -m_a\rbrace,\\
&\lbrace k\in \Z,\, c_{k}=0,\, k\equiv w(d)-1\pmod{n}\rbrace=\lbrace b-1+qn,\, q< -m_b\rbrace.
\end{align*}

By Lemma \ref{lem:indices}, we conclude the proof by noticing that the set of boxes $s\in\lambda$ such that $h_s=n-i$ such that the indices $i_s$ and $j_s$, are congruent to $w(c)-1$ respectively $w(d)-1$ $\pmod{n}$ corresponds to the following set:
\[\lbrace (a-1+kn,b-1+(k+\mathds{1}_{a>b}(a,b))n),\, -m_a\leq k \leq -m_b-1-\mathds{1}_{a>b}(a,b)\rbrace.
\]

Therefore, setting $k''=k+m_a+\mathds{1}_{a>b}(a,b)$, $k''$ follows the inequalities \eqref{eq:ineg}, which concludes the proof.
\end{proof}

\begin{Remark}
As a corollary of Proposition \ref{prop:inversion}, one derives the well-known result (see, for instance,~\cite[Proposition $1.3$]{LLMSSZ} and~\cite[Theorem $7$]{LapointeMorse}):
\[\#R(w)=\lbrace h\in\cch(\lambda),\, h<n\rbrace.\]
\end{Remark}

\begin{Example}\label{ex:exfin}
Take $n=3$ and $w\in\widetilde{\mathfrak{S}}^0_3$ such that $w(1)=-3$, $w(2)=1$ and $w(3)=8$. Then $w^{-1}(1)=2$, $w^{-1}(2)=-3$ and $w^{-1}(3)=7$ and
$\lambda=\psi^{-1}(J_3(w))=(5,3,1,1)$. Its corresponding binary word is drawn in Figure~\ref{fig:word111}.

\begin{figure}[th]
\centering
\begin{tikzpicture}
  [
    dot/.style={circle,draw=black, fill,inner sep=1pt},
  ]

\foreach \x in {0,...,1}{
  \node[dot] at (\x,-4){ };
}

\foreach \x in {.2}
  \draw[->,thick,orange] (\x,-4) -- (\x+.6,-4);
\foreach \x in {1.2}
  \draw[->,thick,orange] (\x,-2) -- (\x+.6,-2);

\foreach \x in {3.2,4.2}
  \draw[->,thick,orange] (\x,-1) -- (\x+.6,-1);
\foreach \y in {.8}
  \draw[->,thick,orange] (5,-\y) -- (5,-\y+.6);
\foreach \y in {1.8}
  \draw[->,thick,orange] (3,-\y) -- (3,-\y+.6);

\draw[->,thick,orange] (1,-3.8) -- (1,-3.8+.6);
\draw[->,thick,orange] (1,-2.8) -- (1,-2.8+.6);
\draw[->,thick,orange] (2.2,-2) -- (2.2+.6,-2);

\node[dot] at (2,-2){};
\node[dot] at (3,-2){};
\node[dot] at (1,-4){};
\node[dot] at (1,-3){};
\foreach \x in {3,...,5}
  \node[dot] at (\x,-1){};
%
\node[above,xshift=0.5cm,yshift=1mm] at (5,0) {NE};
\node[above,xshift=-4mm,yshift=1mm] at (0,0) {NW};

\foreach \y in {1,...,3}
  \draw (.1,-\y) -- node[above,xshift=-4mm,yshift=0.2cm] {$\lambda_\y$} (-.1,-\y);
\node[above,xshift=-4mm,yshift=2mm] at (0,-4) {$\lambda_4$};
\node[above,xshift=-4mm] at (0,-5) {SW};
\node at (0,-4.5) {0};
\node at (0,-5.5) {0};
\node at (1,-3.5) {0};
\node at (1,-2.5) {0};
\node at (3,-1.5) {0};
\node at (5,-0.5) {0};

\node at (0.5,-4) {1};
\node at (1.5,-2) {1};
\node at (4.5,-1) {1};
\node at (3.5,-1) {1};
\node at (5.5,-0) {1};
\node at (2.5,-2) {1};

\draw[dotted,gray] (0,-1)--(5,-1);
\draw[dotted,gray] (0,-2)--(3,-2);
\draw[dotted,gray] (0,-3)--(1,-3);

\draw[dotted,gray] (1,0)--(1,-4);
\draw[dotted,gray] (2,0)--(2,-2);
\draw[dotted,gray] (3,0)--(3,-2);
\draw[dotted,gray] (4,0)--(4,-1);


\node[dot] at (5,-1){};
\node[dot] at (5,0){};
\draw[->,thick,-latex] (0,-6) -- (0,-5);
\draw[thick] (0,-6) -- (0,1);
\draw[->,thick,-latex] (-1,0) -- (6,0);
\node[circle,draw=blue,fill=blue,inner sep=0pt,minimum size=5pt] at (2,-2){};

\end{tikzpicture}
\caption{$\lambda=(5,3,1,1)$ and its corresponding binary word.}\label{fig:word111}
\end{figure}

We have that
\begin{gather*}
R_1(w)=\lbrace \varepsilon_3-\varepsilon_2+k\delta,\, 1\leq k\leq 3\rbrace,\\
R_2(w)=\lbrace \varepsilon_1-\varepsilon_2+k\delta,\, 0\leq k\leq 1\rbrace\cup \lbrace \varepsilon_3-\varepsilon_1+\delta\rbrace.
\end{gather*}
Moreover, $\#R_2(w)=\#\lbrace s\in\lambda,\, h_s=1\rbrace$ and $\#R_1(w)=\#\lbrace s\in\lambda,\, h_s=2\rbrace$.
\end{Example}

\subsection*{Acknowledgements}
The authors would like to thank the referees as well as the editors for their precious comments and suggestions to help this paper to be much clearer.
The authors would also like to thank Olivier Brunat, Nathan Chapelier-Laget, Thomas Gerber, Igor Haladjian, Guo-Niu Han and Ole Warnaar for their interest and their precious remarks as well as Philippe Nadeau for suggesting the reformulation of Theorem~\ref{thm:enumhook}.
Both authors are supported by the Agence Nationale de la Recherche funding ANR CORTIPOM 21-CE40-001.


\pdfbookmark[1]{References}{ref}
\LastPageEnding

\end{document}